\documentclass[10pt]{article}
\usepackage{hyperref}
\usepackage{amsfonts}
\usepackage{mathrsfs}
\usepackage{array}
\usepackage{amsmath,amsthm,amsfonts,amssymb,epsfig}
\usepackage{soul}
\usepackage{color}
\usepackage{stmaryrd}
\usepackage{hyperref}
\soulregister\ref7
\soulregister\bf7

\newtheorem{thm}{Theorem}[section]

\newtheorem{lem}[thm]{Lemma}

\theoremstyle{definition}

\theoremstyle{remark}

\numberwithin{equation}{section}
\usepackage[framemethod=tikz]{mdframed}
\usepackage{xcolor}
\newmdenv[
linecolor=yellow, 
backgroundcolor=yellow!100, 
innerlinewidth=0pt, 
outerlinewidth=0pt, 
leftmargin=0pt, 
rightmargin=0pt
]{highlighted}

\newcommand\be{\begin{equation}}
\newcommand\ee{\end{equation}}
\newcommand\bes{\begin{eqnarray}}
\newcommand\ees{\end{eqnarray}}
\newcommand\bess{\begin{eqnarray*}}
\newcommand\eess{\end{eqnarray*}}

\def\theequation{\arabic{section}.\arabic{equation}}
\begin{document}
	
\title{Carleman estimate for full-discrete approximations of the complex Ginzburg-Landau equation with dynamic boundary conditions and applications to controllability}
	
\author{ Xu Zhu$^1$,\quad Wenwen Zhou$^1$,\quad  Bin Wu$^{1,2}$\footnote{Corresponding author. email: binwu@nuist.edu.cn}
   \\
$ ^1$School of Mathematics and Statistics\\ Nanjing University of
Information Science
and Technology\\
Nanjing 210044, China\\
$ ^2$Center for Applied Mathematics of Jiangsu Province \\
Nanjing University of Information Science and
Technology \\
Nanjing 210044, China\\
}

\maketitle

\begin{abstract}

In this paper, we investigate Carleman estimate and controllability result for the fully-discrete approximations of a one-dimensional Ginzburg-Landau equation with dynamic boundary conditions.  We first establish a new discrete Carleman estimate for the corresponding adjoint system. Based on this Carleman estimate, we obtain a relaxed observability inequality for the adjoint system, and then a controllability result for the fully-discrete Ginzburg-Landau equation with dynamic boundary conditions.


\end{abstract}


\section{Introduction}

\subsection{Motivation}

Let $T>0$, $\Omega=(0,1)$ and $\omega$ be a nonempty subset of $\Omega$.
Given the parameters $\alpha>0$, $c\in \mathbb{R} $ and $\beta, \gamma \in \mathbb{R} \backslash \{0\}$, we consider the following linear Ginzburg-Landau equation with dynamic boundary conditions [\ref{Correa-NARWA-2018}]:
\begin{equation} \label{eq1}
	\left\{
\begin{aligned}
&\partial_{t}y-(\alpha + i\beta)\partial_{x}^{2}y+(c+i\gamma)y= {\bf 1}_{\omega}v,
	&&(x,t)\in \Omega \times (0,T),\\
&\partial_t {y}_{\Gamma_{0}}(t)-(\alpha+i\beta)\partial_{x}y(0,t)+(c+i\gamma)y_{\Gamma_0}= 0,
	&&t\in (0,T),\\
&\partial_t {y}_{\Gamma_{1}}(t)+(\alpha+i\beta)\partial_{x}y(1,t)+(c+i\gamma)y_{\Gamma_1}=0,
	&&t\in (0,T),\\
&y_{\Gamma_{0}}(t)=y(0,t),\quad y_{\Gamma_{1}}(t)=y(1,t),
	&&t\in (0,T),\\
&y(x,0) = y^{0}(x), &&x\in \Omega,\\
&y_{\Gamma_{0}}(0)=y_{\Gamma_{0}}^{0},\quad y_{\Gamma_{1}}(0)=y_{\Gamma_{1}}^{0},
	&&x\in \Omega.
\end{aligned}
\right.
\end{equation}
Here $(y_{\Gamma_{0}},y,y_{\Gamma_{1}}) = ( y_{\Gamma_{0}}(t),y(x,t),y_{\Gamma_{1}}(t))$ is the state of the system, $(y^{0}_{\Gamma_{0}}, y^{0}, y^{0}_{\Gamma_{1}})$ is the initial condition, and the control function $v$ is acting on an open subinterval $\omega \subset \Omega$ through ${\bf 1}_{\omega}(x)$, the characteristic function of $\omega$.

Recently, a variety of mathematical models with dynamic boundary conditions have attracted much attention owing to their extensive applications, such as in population dynamics [\ref{Farkas-MBE-2010}], heat transport in different media [\ref{Gal-DCDS-2008}] and oscillations under strong friction [\ref{Bothe-SPTHA-2005}]. We refer to [\ref{Carreno-arXiv-2023}] for the well-posedness of system (\ref{eq1}) under suitable assumptions in the class of $( y_{\Gamma_0}, y,y_{\Gamma_1})\in C([0,T];\mathbb H^1)$, where
\begin{align*}
\mathbb H^{1}:=\{(u_{\Gamma_0}, u,  u_{\Gamma_1})\in  \mathbb R\times L^2(\Omega)\times \mathbb R; \ u\in H^1(\Omega), u(0)=u_{\Gamma_0}, u(1)=u_{\Gamma_1}\}.
\end{align*}

The null controllability for (\ref{eq1}) is formulated as follows: for any initial state $( y^0_{\Gamma_{0}}, y^0, y^0_{\Gamma_{1}})$ $\in \mathbb H^1$, one can find a control $v\in L^2(\omega\times(0,T))$ such that the associated solution $(y_{\Gamma_0}, y, y_{\Gamma_1})$ of (\ref{eq1}) satisfies 
\begin{align*}
y(x,T)=0,\quad {x}\in\Omega,\quad {\rm and}\quad \left(y_{\Gamma_0}(t), y_{\Gamma_1}(t)\right)=(0,0),\quad t\in (0,T). 
\end{align*}
In [\ref{Carreno-arXiv-2023}], the null controllability for a cubic Ginzburg-Landau equation with dynamic boundary conditions was studied. More precisely, consider the following cubic Ginzburg-Landau equation with dynamic
boundary conditions
\begin{equation} \label{lxbb}
	\left\{
\begin{aligned}
&\partial_{t}y-(\alpha+ i\beta)\Delta y+(c+i\gamma )|y|^2y = {\bf 1}_{\omega}v,
	&&(x,t)\in \Omega \times (0,T),\\
&\partial_{t}y_{\Gamma}+(\alpha+ i\beta)\partial_{\nu}y-(\alpha+i\beta)\Delta_{\Gamma}y_{\Gamma} +(c+i\gamma )|y_\Gamma|^2y_{\Gamma}=0,
	&&(x,t)\in\Gamma \times (0,T),\\
&y=y_{\Gamma},
	&&(x,t)\in\Gamma \times (0,T),\\
&(y(0),y_{\Gamma}(0)) = (y^0, y^0_{\Gamma}),
	&&(x,t)\in\Omega\times\Gamma,
\end{aligned}
\right.
\end{equation}
where $\Omega\subset \mathbb R^{d}$ $ (d\geq 1)$ is a bounded domain with boundary $\Gamma$ of class $C^2$, $\partial_{\nu} y$ is the
normal derivative associated to the outward normal $\nu$ of $\Omega$ and $\Delta_\Gamma$ denotes the Laplace-Beltrami operator on $\Gamma$. Firstly the authors investigated the null controllability for a linearized system by the duality approach and an appropriate Carleman estimate. Further,  the local null
controllability of the nonlinear system was obtained by an inverse mapping theorem. From proposition 4.3 in [\ref{Carreno-arXiv-2023}] the null controllability for (\ref{eq1}) could be deduced directly.   There are also some known controllability and observability results about the deterministic and stochastic
 complex Ginzburg-Landau equations in continuous settings [\ref{Fu-JFA-2009}, \ref{Fu-JDE-2017}, \ref{Rosier-Poincare-2009}].

The main objective of this paper is to investigate discrete Carleman estimate and controllability for the full discrete approximations of system (\ref{eq1}). It is well known that the null controllability result does not necessarily remain valid in the discrete setting. For example, a specific 2-D counterexamples was proposed for which a spatial discrete parabolic equation is not approximately controllable for discrete parameter $\Delta x$ in [\ref{Zuazua-SIREV-2005}]. Thus, we are interested in the $\phi(\Delta x)$-null controllability of the system (\ref{eq1}). More precisely, we want to obtain uniformly bounded controls such that the norm of the semi-discrete solution at time $T$, $y(T)$, is approximately of the size $\sqrt{\phi(\Delta x)}$ with a function $\phi$ that tends to zero exponentially  as space
discretization parameter tends to zero.

In order to formulate our discrete controllability, let us consider $M, N\in \mathbb{N}^*$, the space and time-discretization parameters $\Delta x=\frac{1}{M+1}$ and $\triangle t=\frac{T}{N}$. Furthermore, we introduce
an equidistant spatial mesh of $(0, 1)$, $0=x_0<x_1<\cdots<x_{M+1}=1$, and an equidistant time mesh of $(0, T)$, $0=t^0<t^1<\cdots<t^N=T$. Let $x_{j} = jh$, $t^{n}=n\triangle t$ and 
\begin{align*}
g_j=\left\{\begin{array}{ll} y^0_{\Gamma_0}, &j=0,\\
                             y^0(x_j), & 1\leq j\leq M,\\
                             y^{0}_{\Gamma_1},& j=M+1,
\end{array}\right.
\quad {\rm and}\quad  y^n_j=y(x_j,t_n),\quad 1\leq j\leq M, 1\leq n\leq N.
\end{align*} 
Then the full-discrete approximations of (\ref{eq1}) can read as follows
\begin{equation}\label{1.3}
	\left\{
	\begin{array}{ll}
		\frac{y^{n+1}_j-y^{n}_j}{\Delta t}-(\alpha+i\beta)\frac{y^{n+1}_{j+1}-2y^{n+1}_j+y^{n+1}_{j-1}}{(\Delta x)^2}\\
\hspace{4.3cm}+(c+i\gamma)y_j^{n+1}={\bf 1}_{\omega} v^{n+1}_j, & 1\leq j\leq M,\ 0\leq n\leq N-1, \\
\frac{y^{n+1}_0-y^{n}_0}{\Delta t}-(\alpha+i\beta)\frac{y^{n+1}_{1}-y^{n+1}_0}{\Delta x}+(c+i\gamma)y_0^{n+1}=0, & 0\leq n\leq N-1,\\
		\frac{y^{n+1}_{M+1}-y^{n}_{M+1}}{\Delta t}+(\alpha+i\beta)\frac{y^{n+1}_{M}-y^{n+1}_{M-1}}{\Delta x}+(c+i\gamma)y_{M+1}^{n+1}=0, & 0\leq n\leq N-1,\\
		y^{0}_{j}=g_j,  &0\leq j\leq M+1.
	\end{array}
	\right.
\end{equation}

 Now, we state our discrete null controllability problem in this paper. To formulate it, we introduce the notation $$g=(g_0, g_1\cdots,g_{M+1})^{\rm T}\quad{\rm and}\quad y^n=(y^n_0, y^n_1,\cdots,y^n_{M+1})^{\rm T}, \ n=0,1,\cdots,N$$
to denote discrete functions.
 \vspace{2mm}
 
\noindent{\bf Discrete $\phi(\Delta x)$-null controllability problem.} \ For any initial data $g\in L^2(\overline{\mathcal M})$, find control $v\in {L}^2(\omega\times \mathcal N^*)$ such that the corresponding solution $y$ satisfies 
 $$ \left\|y^{N}\right\|_{{L}^2(\overline{\mathcal M})}\leq \sqrt{\phi(\Delta x)}\left\|g\right\|_{{L}^2(\overline{\mathcal M})}$$
with a function $\phi(\Delta x)$, where the notations of discrete domains $\mathcal M$, $\mathcal N^*$,  discrete function spaces $L^2(\overline{\mathcal M})$ and ${L}^2(\omega\times \mathcal N^*)$, and discrete norm $\|\cdot\|_{{L}^2(\overline{\mathcal M})}$ will be presented in Section 1.2. 

\vspace{2mm}

There has been much work about the study of controllability for various partial differential equations with dynamic boundary conditions. Such sort of
works mainly focus on the continuous cases. We refer to [\ref{Khoutaibi-DCDS-S-2022}] for parabolic equations, [\ref{Mercado-SICON-2023}] for Schr\"{o}dinger equations,  [\ref{Baroun-JDCS-2023}, \ref{Baroun-arXiv-2024}] for stochastic parabolic equations and [\ref{Carreno-arXiv-2023}] for  Ginzburg-Landau equations.  Moreover, as for Carleman estimates for complex partial differential equations, we refer to [\ref{Fu-JFA-2009}] for parabolic equation with a complex principal part, [\ref{Baudouin-IP-2002}] for Schr\"{o}dinger equations and [\ref{Fu-SICON-2017}, \ref{Fu-JDE-2017}] for complex Ginzburg–Landau equations.

There are rich references on discrete Carleman estimates for various partial differential equations. In [\ref{Baudouin-SICON-2023}], a discrete Carleman estimate for the spatial semi-discrete hyperbolic equation was established and subsequently applied it to obtain the stability of an inverse problem to recover a potential term in this discrete hyperbolic equation. Discrete Carleman estimate for one-dimensional elliptic equation has been developed in [\ref{Boyer-JMPA-2010}]  to establish a relaxed observability for the associated semi-discrete parabolic equation. As for discrete Carleman estimates for parabolic equations, we refer to [\ref{Boyer-ESAIM:COCV-2005}] for time semi-discrete case, [\ref{Boyer-POINCARE-AN-2014}, \ref{Nguyen-MCRF-2014}] for spatial semi-discrete case,  and [\ref{Carreno-ACM-2023}] for full-discrete case.  These discrete Carleman estimates were applied to obtain relaxed controllability results for the associated discrete parabolic equations. In [\ref{Cerpa-JMPA-2022}] and [\ref{Hernandez-arXiv-2021}], discrete Carleman estimates for fourth-order differential operators were established. Moreover, there are some works about discrete Carleman estimates for stochastic differential equations, e.g. [\ref{Wu-IP-2024},\ref{Zhao-arXiv-2024}] for stochastic second-order semi-discrete parabolic equations and [\ref{Wang-arXiv-2024}] for stochastic fourth-order semi-discrete parabolic equations. To the best of our knowledge, there is only one paper to discuss Carleman estimates and controllability for discrete differential equations with dynamic boundary conditions [\ref{Lecaros-JDE-2023}], in which  Carleman estimate for a fully discrete approximations of the parabolic equation with dynamic boundary conditions was established and then applied to obtain $\phi(\Delta x)$-null controllability. 

However, to the author's knowledge, discrete Carleman estimates and controllability for the full discrete approximations of Ginzburg-Landau equation with dynamic boundary conditions have not been studied thoroughly yet. In order to handle dynamic boundary conditions, we need to carefully combine all boundary terms in proving the Carleman estimates. It is well-known that Ginzburg-Landau equation is closely related to both parabolic equations and Schr\"{o}dinger equations. In order to reveal the distinctions between our Carleman estimate and the ones for these equations, we also decouple $\alpha$ and $\beta$ from the related constants appearing in the proof. It is observed that our results are consistent with those of parabolic equations. However, our method is not applicable to the Schrödinger equations. The reason is that we need to handle the second-order difference operator arising from discretization, which leads to that $\alpha$ could not vanish. 
Our proofs follow the ideas proposed for the continuous Ginzburg-Landau equations in [\ref{Fu-JFA-2009} ] and [\ref{Rosier-Poincare-2009} ]. The primary challenge in developing discrete Carleman estimates lies in computing multiple discrete difference operators on the Carleman weight functions.


\subsection{Discrete settings}
To facilitate problem solving, we introduce some notations.
Let $M,N \in \mathbb{N^{*}}$ and $T>0$.
We use the notation $\llbracket a, b \rrbracket = [a,b] \cap \mathbb{N}$ for any real numbers $a<b$ and set the space- and time-discretization parameters $\Delta x=\frac{1}{M+1}$ and $\Delta t=\frac{T}{N}$.
Let us also define the uniform meshes in space and time as follows:
\begin{align*}
\mathcal K_{\Delta x}=\{x_j:=jh;\  i\in \llbracket 0, M+1 \rrbracket\}\quad {\rm and}\quad \mathcal K_{\Delta t}=\{t^n:=n\Delta t;\  n\in \llbracket 0, N \rrbracket\}.
\end{align*}
Moreover, we denote the discrete approximation of a function $r=r(x,t)$ at a grid point $(x_{j},t^n)\in \mathcal K_{\Delta x}\times \mathcal K_{\Delta t}$  as $r_{j}^{n}:=r(x_{j},t^n)$ and further set
\begin{align*}
r^n=(r^n_0, r^n_1, \cdots, r^n_{M+1})^{\rm T},\quad n=0,1,\cdots, N.
\end{align*}


\subsubsection{Space discretization}
Now, we give the notations for the spatial discretization.
For $\mathcal{W} \subseteq \mathcal{K}_{\Delta x}$, with the translation operators, we can build two dual sets given by
$$\mathcal{W}'={\bf s}_{-}(\mathcal{W}) \cap {\bf s}_{+}(\mathcal{W}),
\quad
\mathcal{W}^{*}={\bf s}_{-}(\mathcal{W}) \cup {\bf s}_{+}(\mathcal{W}),$$
where ${\bf s}_{+}$ and ${\bf s}_{-}$ are the translation operators defined by
\begin{align*}
{\bf s}_{\pm}(\mathcal W)=\left\{x_j \pm \frac{\Delta x}{2};\ x_j\in \mathcal W\right\}.
\end{align*}
The dual of a dual set will be denoted as $\overline{\mathcal{W}}:=(\mathcal{W}^{*})^{*}$ and $\widetilde{\mathcal{W}}:=(\mathcal{W}')'$.
Then, any subset $\mathcal{W} \subset \mathcal{K}_{\Delta x}$ that verifies $\widetilde{\overline{\mathcal{W}}}=\mathcal{W}$ will be called regular mash. Furthermore, we define its boundary by $\partial\mathcal{W}:=\overline{\mathcal{W}} \backslash \mathcal{W}$. We denote the interior of the meshes in space by $\mathcal M:= \widetilde {\mathcal K}_{\Delta x}$.

We denote by $\mathbb{C}(\mathcal W)$ the set of complex-valued function
defined in $\mathcal W \subset \mathcal K_{\Delta x}$. Then for $u\in C(\mathcal W)$, we define the average and the difference operator as
$${\bf A}_{x}u=\frac{{\bf s}_{+}u+{\bf s}_{-}u}{2},
\quad
{\bf D}_{x}u=\frac{{\bf s}_{+}u-{\bf s}_{-}u}{\Delta x},$$
where ${\bf s}_{\pm}u(x_j)=u(x_j\pm\frac{\Delta x}{2})$ for $x_j\in \mathcal W$. For $u, v\in \mathbb{C}(\overline {\mathcal W})$, we define the integral of $u$ on $\mathcal{W}$ and $\partial \mathcal W$
$$\int\limits_{\mathcal{W}}u
=\sum_{x_j \in \mathcal{W}}u(x_j)\Delta x,\quad \int\limits_{\partial\mathcal{W}}u=\sum\limits_{x_j \in \partial\mathcal{W}}u(x_j)$$
and the following $L^{2}$-inner product
$$(u,v)_{L^{2}(\mathcal{W})}
=\int\limits_{\mathcal{W}}uv^*
=h\sum_{x_j \in \mathcal{W}}u(x_j)v^*(x_j),$$
where $v^*$ denotes the conjugate of $v$. The associated norm is denoted by $\|u\|_{L^2(\mathcal W)}$, i.e.
\begin{align*}
\|u\|_{L^2(\mathcal W)}=\left(\sum\limits_{x_j \in \mathcal{W}}|u(x_j)|^2\Delta x\right)^{\frac{1}{2}}.
\end{align*}
 Analogous definitions are for $\|u\|_{L^2(\partial\mathcal W)}$ and $\|u\|_{L^{\infty}(\mathcal{W})}$ in terms of
$$\|u\|_{L^2(\partial\mathcal W)}=\left(\sum\limits_{x_j \in \partial\mathcal{W}}|u(x_j)|^2\right)^{\frac{1}{2}},\quad \|u\|_{L^{\infty}(\mathcal{W})}=\sup \limits_{x_j \in \mathcal{W}}|u(x_j)|.$$
In addition, we define the ${L}^{2}$-norm for $u\in \mathbb{C}(\overline{\mathcal W})$ by
$$||u||^{2}_{{L}^{2}(\overline{\mathcal{W}})}:=||u||^{2}_{L^{2}(\mathcal{W})}+||u||^{2}_{L^{2}(\partial\mathcal{W})}.$$

To introduce the boundary conditions, for $x_j \in \partial\mathcal{W}$ we define the outward normal at $x_j$,
\begin{equation}\label{1-1.4}
n_x(x_j):=
	\left\{
\begin{aligned}
&1,
	&&{\bf s}_{-}(x_j) \in \mathcal{W}^{*} \; \text{and} \; {\bf s}_{+}(x_j) \notin \mathcal{W}^{*},\\
-&1,
	&&{\bf s}_{-}(x_j) \notin \mathcal{W}^{*} \; \text{and} \; {\bf s}_{+}(x_j) \in \mathcal{W}^{*},\\
&0,
	&&\text{otherwise}.\\
\end{aligned}
\right.
\end{equation}
and the trace operator for $u \in \mathbb{C}(\mathcal{W}^{*})$ as
\begin{equation}\label{1-1.5}
{\bf tr}(u)(x_j):=
	\left\{
\begin{aligned}
&{\bf s}_{-}u(x_j),
	&&n_x(x_j)=1,\\
&{\bf s}_{+}u(x_j),
	&&n_x(x_j)=-1,\\
&0,
	&&n(x_j)=0.\\
\end{aligned}
\right.
\end{equation}

\subsubsection{Time discretization}
Now, we will introduce some notations to define the discretization of the time variable.
We define $\mathcal{N}:=\{t^{n};\ n \in \llbracket 1, N \rrbracket \}$ as the primal mesh in time and $\mathcal{N}^{*}:=\{t^{n+\frac{1}{2}};\ n \in \llbracket 0, N-1 \rrbracket \}$ as the dual mesh in time.
Also we write $\partial \mathcal N=\{0, T\}$, $\overline{\mathcal{N}}=\mathcal{N} \cup \{0\}$ and $\overline{\mathcal{N}^{*}}=\mathcal{N}^{*} \cup \{T+\frac{\triangle t}{2}\}$.

We define the time discrete derivative of a function $y:[0,T] \rightarrow \mathbb{R}$ sampled on $\overline{\mathcal{N}}$ as follows:
$${\bf D}_{t}y:=\frac{{\bf t}^{+}y-{\bf t}^{-}y}{\triangle t},$$
where  ${\bf t}^{\pm}y(t^n):=y(t^n \pm \frac{\triangle t}{2})$ for $t^n \in \mathcal{N}$.
Similarly to the discrete spatial variable, we define the outward normal for $t^n\in \partial \mathcal N$ as
\begin{equation}\label{1.6}
n_t(t^n):=
	\left\{
\begin{aligned}
&1,
	&&t^n=T,\\
-&1,
	&&t^n=0.
\end{aligned}
\right.
\end{equation}
 For $u\in \mathbb{C}(\overline {\mathcal N})$, we define the discrete integral of $u$ with respect to time
$$\int\limits_{\mathcal{N}}\hspace{-1.05em}-\ u
=\sum_{t^n \in \mathcal{N}}u(t^n)\Delta t,\quad \int\limits_{\partial\mathcal{N}}\hspace{-1.05em}-\ u
=\sum_{t^n \in \partial\mathcal{N}}u(t^n).$$

Combining space and time discrete integrals, for $u\in \mathbb{C}(\mathcal K_{\Delta x}\times\mathcal K_{\Delta t})$ we define 
\begin{align*}
&\iint\limits_{\mathcal M\times\mathcal N}u=\int\limits_{\mathcal{N}}\hspace{-1.05em}-\int\limits_{\mathcal M}u=\sum_{t^n \in \mathcal{N}}\sum_{x_j \in \mathcal{M}}u(x_j,t^n)\Delta x\Delta t,\\
&\iint\limits_{\partial\mathcal M\times\mathcal N}u=\int\limits_{\mathcal{N}}\hspace{-1.05em}-\int\limits_{\partial\mathcal M}u=\sum_{t^n \in \mathcal{N}}\sum_{x_j \in \partial\mathcal{M}}u(x_j,t^n)\Delta t,\\
&\iint\limits_{\mathcal M\times\partial\mathcal N}u=\int\limits_{\partial\mathcal{N}}\hspace{-1.05em}-\int\limits_{\mathcal M}u=\sum_{t^n \in \partial\mathcal{N}}\sum_{x_j \in \mathcal{M}}u(x_j,t^n)\Delta x,\\
&\iint\limits_{\partial\mathcal M\times\partial\mathcal N}u=\int\limits_{\partial\mathcal{N}}\hspace{-1.05em}-\int\limits_{\partial\mathcal M}u=\sum_{t^n \in \partial\mathcal{N}}\sum_{x_j \in \partial\mathcal{M}}u(x_j,t^n),\\
&\iint\limits_{\{x_j\}\times\mathcal N}u=\int\limits_{\mathcal{N}}\hspace{-1.05em}-u(x_j,t^n)=\sum_{t^n \in \mathcal{N}}u(x_j,t^n)\Delta t,\quad x_j\in \partial \mathcal M
\end{align*}
and the $L^2$-inner product in $\mathcal M\times \mathcal N$
\begin{align*}(u,v)_{L^{2}(\mathcal{M}\times\mathcal N)}
=\iint\limits_{\mathcal{M}\times\mathcal N}uv^*
=\sum_{t^n\in \mathcal N}\sum_{x_j \in \mathcal{M}}u(x_j,t^n)v^*(x_j,t^n)\Delta x\Delta t.
\end{align*}
The associated norm is denoted by $\|u\|_{L^2(\mathcal M\times\mathcal N)}$, i.e.
\begin{align*}
\|u\|_{L^2(\mathcal M\times \mathcal N)}=\left(\sum\limits_{t^n \in \mathcal{N}}\sum\limits_{x_j \in \mathcal{M}}|u(x_j, t^n)|^2\Delta x\Delta t\right)^{\frac{1}{2}}.
\end{align*}
Analogous definitions are for $(\cdot,\cdot)_{L^2({\{0\}\times\mathcal N})}$ and $(\cdot,\cdot)_{L^2({\{1\}\times\mathcal N})}$, which represent the $L^2$-inner products on the boundary.

\subsection{Main results}

For our discrete controllability problem, we can use a suitable discrete framework to process it.
In this paper, we use a backward difference for time derivative and a centered second order difference for the space derivative.
Thus, the full-discrete approximations \eqref{1.3} can be written as
\begin{equation} \label{1.4}
	\left\{
\begin{aligned}
&{\bf D}_{t}y-(\alpha + i\beta){\bf D}_{x}^{2}{\bf t}^{+}(y)+(c+i\gamma){\bf t}^+(y)=1_{\omega}v,
	&&(x,t) \in \mathcal{M} \times \mathcal{N}^{*},\\
&{\bf D}_{t}y(0,t)-(\alpha+i\beta){\bf D}_{x}{\bf t}^{+}({\bf s}_{+}y)(0,t)+(c+i\gamma){\bf t}^+(y)(0, t)=0,
	&&t \in \mathcal{N}^{*},\\
&{\bf D}_{t}y(1,t) + (\alpha+i\beta){\bf D}_{x}{\bf t}^{+}({\bf s}_{-}y)(1,t)+(c+i\gamma){\bf t}^+(y)(1, t)=0,
	&&t\in \mathcal{N}^{*},\\
&y(x,0)=g(x),
	&&x\in \overline{\mathcal{M}},
\end{aligned}
\right.
\end{equation}
where $g\in L^2(\overline{\mathcal{M}})$ stands for the initial data, $v$ is the discrete control function acting on a subset $\omega \subset \Omega$. 

Due to the so-called Hilbert Uniqueness Method [\ref{Lions-SIREV-1988}], it is well-known that the controllability result of (\ref{1.4})
can be proved by observability estimate for the following adjoint system:
\begin{equation} \label{dofc}
	\left\{
\begin{aligned}
&-{\bf D}_{t}q-(\alpha - i\beta){\bf D}_{x}^{2}{\bf t}^{-}(q)+(c-i\gamma){\bf t}^-(q)=0,
	&&(x,t) \in \mathcal{M} \times \mathcal{N},\\
&-{\bf D}_{t}q(0,t)-(\alpha-i\beta){\bf D}_{x}{\bf t}^{-}({\bf s}_{+}q)(0,t)+(c-i\gamma){\bf t}^-(q)(0,t) = 0,
	&&t \in \mathcal{N},\\
&-{\bf D}_{t}q(1,t) + (\alpha-i\beta){\bf D}_{x}{\bf t}^{-}({\bf s}_{-}q)(1,t)+(c-i\gamma){\bf t}^-(q)(1,t) = 0,
	&&t \in \mathcal{N},\\
&q(x,T+{\Delta t}/{2})=q_{T}(x),
	&&x \in \overline{\mathcal{M}}.
\end{aligned}
\right.
\end{equation}

To establish the observability inequality, we need to obtain a discrete Carleman estimate inequality for the solution $q\in \mathbb{C}(\overline {\mathcal M}\times\overline{\mathcal N ^*})$ of \eqref{dofc}.
To do this, we first introduce the following weight function in Carleman estimate.
Let $\omega_0$ be an open subset of $\Omega$ such that ${\omega}_{0} \subset \subset{\omega}$.
Then, we can find a function $\psi:\overline{\Omega}\to\mathbb{R}$ satisfying
\begin{align} \label{psi}
\psi>0 \; \text{in} \; \Omega,
\quad
|\nabla \psi|>c_0 \; \text{in} \; \Omega \setminus \overline{\omega_{0}} \;
\quad\text{and} \quad
\psi_{x}(0)>0, \; \psi_{x}(1)<0
\end{align}
with some $c_0>0$. Further, we introduce
\begin{align} \label{varphi}
\varphi(x)= e^{\lambda \psi(x)} - e^{\lambda K}<0,
\quad
\phi(x)= e^{\lambda \psi(x)},
\quad
x \in \overline{\Omega},
\end{align}
and
\begin{align}
\theta(t)=\frac{1}{(t+\delta T)(T+\delta T-t)},
\quad
t \in [0,T],
\end{align}
where $\lambda>0$ is a large parameter, $K>||\psi||_{C(\overline{\Omega})}$ and $0<\delta<\frac{1}{2}$.
Here $\delta$ is introduced to avoid singularities at time $t=0$ and $t=T$. We set 
\begin{align*}
s(t)=\tau \theta(t),\quad r(x,t)=e^{s(t)\varphi(x)},\quad \rho(x,t)=r^{-1}(x,t)
\end{align*}
with the second large parameter $\tau>0$. Finally, we also use the following notations:
\begin{align*}\begin{array}{ll}
\mathcal{P}(q):=-{\bf D}_{t}q-(\alpha-i\beta){\bf D}_{x}^{2}{\bf t}^{-}(q), &(x,t)\in\mathcal{M} \times \mathcal{N},\\
\mathcal B_{\Gamma_0}(q):={\bf D}_{t}q(0,t)+(\alpha-i\beta){\bf D}_{x}{\bf t}^{-}({\bf s}_{+}q)(0,t), &t\in\mathcal{N},\\
\mathcal B_{\Gamma_1}(q):={\bf D}_{t}q(1,t)-(\alpha-i\beta){\bf D}_{x}{\bf t}^{-}({\bf s}_{-}q)(1,t), &t\in\mathcal{N}.
\end{array}
\end{align*}

 In the following, we will use $C$ to denote generic positive constant depending on $\alpha, \beta, \delta, \omega$, $T$, $c_0$, but independent of $\lambda$ and $\tau$. Similarly, $C(\lambda)$ denotes constant also depending on $\lambda$. Moreover, we use the notation $\mathcal O_\lambda(\gamma)$, which satisfies $|\mathcal O_\lambda(\gamma)|\leq C(\lambda)|\gamma|$ with a constant $C(\lambda)$. All of these notations may vary from line to line and are independent of $\Delta x$ and $\Delta t$. Also, we denote for discrete function $q$ defined in $\mathcal M\times \mathcal N^*$ that
 \begin{align*}
 \iint\limits_{\omega\times \mathcal N^*}u=\iint\limits_{(\omega\cap\mathcal M)\times \mathcal N^*} u.
 \end{align*}

The first main result is the following discrete Carleman estimate related to the full-discrete Ginzburg-Landau operator with dynamic boundary conditions. 
 
\begin{thm} \label{carleman estimate}
Let $\omega_{0}$ be a nonempty subset of $\Omega$ such that $\omega_{0} \subset \subset \omega$,
$\psi$ be given by \eqref{psi} and define $\varphi$ according to \eqref{varphi}.
For the parameter $\lambda \ge 1$ sufficiently large, there exist positive constants $C$ depending on $\alpha, \beta, \delta, \omega$, $T, c_0$, and $C(\lambda)$, $\tau_{0}>1$, $0<\widehat{\Delta x}<1$, $0<\varepsilon_{0}<1$ also depending on $\lambda$ such that
\begin{align} \label{CE}
&\iint\limits_{\mathcal{M} \times \mathcal{N}^{*}} s^{-1}r^{2} |{\bf D}_{x}^{2}q|^{2}+\iint\limits_{\mathcal{M} \times \mathcal{N}} {\bf t}^{-}(s^{-1}r^{2}) |{\bf D}_{t}q|^{2}+
\iint\limits_{\mathcal{M} \times \mathcal{N}^{*}} sr^{2} |{\bf A}_{x}{\bf D}_{x}q|^{2}+
\nonumber \\
&\iint\limits_{\mathcal{M}^{*} \times \mathcal{N}^{*}} sr^{2} |{\bf D}_{x}q|^{2}+\iint\limits_{\mathcal{M} \times \mathcal{N}^{*}} s^{3}r^{2} |q|^{2}
+\iint\limits_{\partial\mathcal{M} \times \mathcal{N}} {\bf t}^{-}(s^{-1}r^{2}) |{\bf D}_{t}q|^{2}+\nonumber \\
&\iint\limits_{\partial\mathcal{M} \times \mathcal{N}^{*}} sr^{2} {\bf tr}(|{\bf D}_{x}q|^{2})
+\iint\limits_{\partial\mathcal{M} \times \mathcal{N}^{*}} s^{3}r^{2} |q|^{2} \nonumber \\
\le\ &C\left(\ \iint\limits_{\mathcal{M} \times \mathcal{N}} {\bf t}^{-}(r^{2})|\mathcal{P}(q)|^{2}+\iint\limits_{\{0\}\times\mathcal{N}} {\bf t}^{-}(r^{2}) |\mathcal B_{\Gamma_{0}}(q)|^{2}
+\iint\limits_{\{1\}\times\mathcal{N}} {\bf t}^{-}(r^{2}) |\mathcal B_{\Gamma_1}(q)|^{2}\right)+\nonumber\\
&C(\lambda) \iint\limits_{\omega\times \mathcal{N}^{*}} s^{3}r^{2} |q|^{2}+C(\Delta x)^{-2} \left(\ \iint\limits_{\mathcal{M} \times \partial\mathcal{N}} {\bf t}^{+}(r^2){\bf t}^{+}(|q|^{2})
+\iint\limits_{\partial\mathcal{M} \times \partial\mathcal{N}} {\bf t}^{+}(r^2){\bf t}^{+}(|q|^2)\right)
\end{align}
for all $\tau \ge \tau_{0}(T+T^{2})$, $0<\Delta x\le \widehat {\Delta x}$, $0<\Delta t\leq 1$ and $0<\delta \le \frac{1}{2}$ satisfying the conditions
\begin{align}\label{1.10}
\frac{\tau \Delta x}{\delta T^{2}} \le \varepsilon_{0},
\quad
\frac{\tau^{4}\Delta t}{\delta^{4}T^{6}} \le \varepsilon_{0}.
\end{align}
\end{thm}

\noindent{\bf Remark 1.1.}\ 
In this discrete Carleman estimate, we first choose $\lambda_0\geq C$ with some sufficiently large constant $C$. Then for  fixed sufficiently large $\lambda$, we can choose $\widehat {\Delta x}$ sufficiently small to satisfy \eqref{1.10}. This leads to that we could absorb some terms on the right-hand side of Carleman estimate.

\vspace{2mm}

\noindent{\bf Remark 1.2.}\  The class of discrete Ginzburg-Landau operator $-{\bf D}_{t}q-(\alpha-i\beta){\bf D}_{x}^{2}{\bf t}^{-}(q)$ contains both the discrete heat operator $-{\bf D}_{t}q-\alpha{\bf D}_{x}^{2}{\bf t}^{-}(q)$ in the limit $\beta\rightarrow 0$ and the Schr\"{o}dinger operator $-i{\bf D}_{t}q-{\bf D}_{x}^{2}{\bf t}^{-}(q)$ in the limit $\alpha\rightarrow 0$. From the proof of our Carleman estimate, by selecting $\beta\rightarrow 0$ we still obtain the estimate (\ref{CE}), which is consistent with the Carleman estimate for the fully-discrete parabolic equation with dynamic boundary conditions in [\ref{Lecaros-JDE-2023}]. Unfortunately, our method fails to establish a Carleman estimate for the full-discrete Schr\"{o}dinger equation as $\alpha$ approaches zero. The reason is that we need an estimate (\ref{dh2z}) below for $\alpha{\bf D}_x^2 z$ to handle errors arising from discretization, e.g. $$-\alpha\iint\limits_{\mathcal{M} \times \mathcal{N}^{*}} s^{-3}(s\Delta x)^{6}\mathcal O_{\lambda}(1) |{\bf D}_{x}^{2}z|^{2},\quad -\beta\iint\limits_{\mathcal{M} \times \mathcal{N}^{*}} s^{-1}(s\Delta x)^{2}\mathcal O_{\lambda}(1) |{\bf D}_{x}^{2}z|^{2}$$  in (\ref{A.7}) and (\ref{A.15}) respectively. However, there exists no estimate for the operator ${\bf D}_x^2$ in the Carleman estimate for Schr\"{o}dinger equation. To address this difficulty, it is imperative to either reconfigure the decompositions of the operators $\mathcal A$ and $\mathcal B$ in (\ref{3.1}) below, or adopt an alternative approach to validate the estimates. It seems that even if we apply the decomposition of the weighted Schrödinger operator analogous to the continuous case, we still can not handle the second-order difference operator induced by discretization. Therefore, Carleman estimate for full-discrete Schr\"{o}dinger equation is an interesting topic, which will be  our further work.

\vspace{2mm}

\noindent{\bf Remark 1.3.}\  In the Carleman estimate, we decouple the second parameter $\lambda$ from the related constants
appearing in the proof, which leads to that positive constants $C$ before the operators ${\mathcal P}, {\mathcal B}_{\Gamma_0}$ and $ {\mathcal B}_{\Gamma_1}$  in (\ref{CE}) are independent of $\lambda$. This is essential to deal with controllability and inverse problems for strongly coupled system, e.g. thermoelastic model [\ref{Albano-EJDE-2000}, \ref{Bellassoued-IP-2010}], where the second large parameter $\lambda$ plays a crucial role in handling the strongly coupled terms. The decoupling of the second large parameter $\lambda$ provides a potential strategy for deriving Carleman estimates for discrete strongly coupled systems.

\vspace{2mm}

Based on Carleman estimate (\ref{CE}), we can establish the following observability inequality for the adjoint system (\ref{dofc}).

\begin{thm} \label{gc}
Let $T>0$ and $\widehat {\Delta x}$ given by Theorem \ref{carleman estimate}.
Then for any $\vartheta \ge 1$, there exist positive constants $C$ and $C_{obs}$ depending on $\alpha, \beta, \delta, \omega, T, c_0, c$ and $\gamma$ such that 
\begin{align} \label{gcbds}
\|{\bf t}^{+}(q)(0)\|^2_{{L}^{2}(\mathcal{M})}
\le C_{obs} \left(\ \ \iint\limits_{\omega \times \mathcal{N}^{*}} |q|^{2} + e^{-\frac{C}{h^{\rm {min}\left\{{\vartheta}/{4},1\right\}}}} \|q_{T}\|_{{L}^{2}(\overline{\mathcal{M}})}^{2}\right),
\end{align}
for all $0< \Delta x \le \min\{\widehat{\Delta x},\widetilde{{\Delta x}} \}$ and $0<\Delta t\leq \min\{T^{-2}h^{\vartheta}, (4\rho)^{-1}\}$, where $\rho=\max\{|c|, |\gamma|\}$ and
$$\widetilde{{\Delta x}}=C\left(1+\frac{1}{T}+\rho^{2/3}\right)^{-\max\{1,4/\vartheta\}},\quad C_{obs}=e^{C \left(1+\frac{1}{T}+\rho^{{2}/{3}}+T\rho\right)}.$$
\end{thm}

By the penalized HUM method [\ref{Boyer-ESAIM-2013}, \ref{Lions-1988}], we immediately deduce the following relaxed controllability result from observability inequality (\ref{gcbds}).  

\begin{thm} \label{kz}
Let $T>0$, $\vartheta \ge 1$ and $\Delta x, \Delta t>0$ as in Theorem \ref{gc}.
Then, for any initial data $g \in L^2(\overline{\mathcal{M}})$, there exists a control function $v$ satisfying
$$\|v\|_{{L}^{2}(\omega \times \mathcal{N}^{*})} \le C\|g\|_{{L}^{2}(\overline{\mathcal{M}})}$$
and such that the corresponding controlled solution $y$ to \eqref{1.4} satisfies 
$$\|y(x,T)\|_{{L}^{2}(\overline{\mathcal{M}})} \le C\sqrt{\phi(\Delta x)}\|g\|_{{L}^{2}(\overline{\mathcal{M}})}$$
with $\phi(\Delta x):=e^{-\frac{C}{h^{\rm {min}\left\{{\vartheta}/{4},1\right\}}}} $, where  $C$ is depending on $\alpha, \beta, \delta, \omega, T, c_0, c$ and $\gamma$.
\end{thm}

The proof of Theorem \ref{kz} is standard. Thus we omit it in this paper and refer to [\ref{Carreno-ACM-2023}] for details.

The article is organized as follows.
In section 2, we present discrete settings and the preliminaries of discrete estimates to several applications of the discrete operator on the Carleman weight function.
In section 3, we prove our Carleman estimate (\ref{CE}), i.e Theorem 1.1.
In section 4, we apply this Carleman estimate to prove the observability inequality (\ref{gcbds}), i.e. Theorem 1.2.


\section{Preliminaries}

\subsection{Discrete calculus formulas}

In this section, we state the elementary notions concerning discrete calculus formulas, including some useful identities for the average and difference operators and discrete integration by parts formulas with respect to discrete space and time variables. The detailed proofs of these results can be found in [\ref{Boyer-JMPA-2010}, \ref{Boyer-POINCARE-AN-2014}, \ref{Carreno-ACM-2023}] or [{\ref{Lecaros-JDE-2023}}]. 

\begin{lem}
Let $u, v\in \mathbb{C}(\overline{\mathcal M})$.
Then, for the space difference and the average operators, we have
\begin{align} \label{ds1}
&{\bf D}_{x}(uv)
={\bf D}_{x}u{\bf A}_{x}v+{\bf A}_{x}u{\bf D}_{x}v,\\
\label{ds2}
&{\bf A}_{x}(uv)
={\bf A}_{x}u{\bf A}_{x}v+\frac{1}{4}(\Delta x)^{2}{\bf D}_{x}u{\bf D}_{x}v,\quad {\bf A}_{x}^{2}u
=u+\frac{1}{4}(\Delta x)^{2}{\bf D}_{x}^{2}u.
\end{align}
\end{lem}

The following results represent discrete integration by parts formulas with respect to spatial variable for the difference and average operators [\ref{Lecaros-ESAIM:COCV-2021}].

\begin{lem}
Let $u \in \mathbb{C}(\overline{\mathcal{M}})$ and $v \in \mathbb{C}(\mathcal{M}^{*})$.
Then, the following identities hold for the difference and average operators, respectively
\begin{align} \label{fbjf1}
&\int\limits_{\mathcal{M}}u{\bf D}_{x}v
=-\int\limits_{\mathcal{M}^{*}}v{\bf D}_{x}u+\int\limits_{\partial\mathcal{M}}u{\bf tr}(v)n_x,\\
\label{fbjf2}
&\int\limits_{\mathcal{M}}u{\bf A}_{x}v
=\int\limits_{\mathcal{M}^{*}}v{\bf A}_{x}u-\frac{h}{2}\int\limits_{\partial\mathcal{M}}u{\bf tr}(v),
\end{align}
where $n_x$ denotes the outward normal on $\partial \mathcal M$, and ${\bf tr}$ denotes the trace operator, defined by (\ref{1-1.4}) and (\ref{1-1.5}) respectively.
\end{lem}

As for time discrete calculus, we have the following formulas [\ref{Carreno-ACM-2023}] or [\ref{Lecaros-JDE-2023}].
 
\begin{lem}
For $f, g\in \mathbb{C}(\overline {\mathcal N})$, the following identities holds
\begin{align}\label{2.5}
{\bf D}_{t}(fg)={\bf D}_{t}f{\bf t}^{-}(g)+{\bf t}^{+}(f){\bf D}_{t}g,\quad
{\bf D}_{t}(fg)={\bf D}_{t}f{\bf t}^{+}(g)+{\bf t}^{-}(f){\bf D}_{t}g.
\end{align}
\end{lem}

Obviously, from (\ref{2.5}), we can obtain that
\begin{align}\label{2.6}
\begin{aligned}
&{\bf t}^{+}(f){\bf D}_{t}f^{*}+{\bf t}^{+}(f^{*}){\bf D}_{t}f={\bf D}_{t}(|f|^{2})+\Delta t |{\bf D}_{t}f|^{2}, \\
&{\bf t}^{-}(f){\bf D}_{t}f^{*}+{\bf t}^{-}(f^{*}){\bf D}_{t}f={\bf D}_{t}(|f|^{2})-\Delta t |{\bf D}_{t}f|^{2},
\end{aligned}
\end{align}
if we choose $g$ as $f^*$, the conjugate of $f$.

\begin{lem}
Let us consider $f \in \mathbb{C}(\overline{\mathcal{N}})$ and $g \in \mathbb{C}(\overline{\mathcal{N}^{*}})$.
Then, the following identities hold
\begin{align}\label{fbjf3}
&\int\limits_{\mathcal{N}}\hspace{-1.05em}-\ f{\bf t}^{-}(g)
=\int\limits_{\mathcal{N}^*}\hspace{-1.05em}-\ {\bf t}^{+}(f)g,\\
\label{fbjf4}
&\int\limits_{\mathcal{N}}\hspace{-1.05em}-\ f{\bf D}_{t}g
=-\int\limits_{\mathcal{N}^*}\hspace{-1.05em}-\ g{\bf D}_{t}f+\int\limits_{\partial\mathcal{N}}f{\bf t}^{+}(g)n_{t},
\end{align}
where $\partial\mathcal{N}=\{0,T\}$.
\end{lem}

Furthermore, we have the following useful identities, which can be found in [\ref{Lecaros-JDE-2023}].
\begin{align}
\label{2.10}&\int\limits_{\mathcal{N}}\hspace{-1.05em}-\ {\bf t}^{-}(f){\bf D}_{t}g
=-\int\limits_{\mathcal{N}}\hspace{-1.05em}-\ {\bf D}_{t}f{\bf t}^{+}(g)+\int\limits_{\partial\mathcal{N}}\hspace{-1.05em}-\ {\bf t}^{+}(fg)n_{t}.
\end{align}
for $f,g \in \mathbb{C}(\overline{\mathcal{N}^{*}})$, and
\begin{align}
&\int\limits_{\mathcal{N}^*}\hspace{-1.05em}-\ {\bf t}^{+}(f){\bf D}_{t}g
=-\int\limits_{\mathcal{N}^*}\hspace{-1.05em}-\ {\bf t}^{-}(g){\bf D}_{t}f+\int\limits_{\partial \mathcal{N}}\hspace{-1.05em}-\ fgn_{t}
\end{align}
for $f,g \in \mathbb{C}(\overline{\mathcal{N}})$.

\subsection{Discrete calculus results on the weight function}
This section states some results related to discrete operations performed on the weight functions in the discrete Carleman estimate \eqref{CE}.
We recall the weight functions
\begin{align*}s(t)=\tau \theta(t),\quad r(t,x)=e^{s(t)\varphi(x)},\quad \rho(t,x)=r^{-1}(x,t).
\end{align*}
We presents these results here without a proof and refer to [\ref{Boyer-ESAIM:COCV-2005}, \ref{Boyer-JMPA-2010}] and [\ref{Lecaros-JDE-2023}] for a complete discussion.

\begin{lem} \label{gj1}
Let $\alpha,m,n \in \mathbb{N}$.
Provided ${\tau \Delta x}{(\delta T^{2})}^{-1} \le 1$, we have
\begin{align*}r{\bf A}_{x}^{m}{\bf D}_{x}^{n}\partial_{x}^{\alpha}\rho
=r\partial_{x}^{n+\alpha}\rho
+s^{n+\alpha}{\mathcal O}_{\lambda}\left((s\Delta x)^{2}\right)
=s^{n+\alpha}\mathcal O_{\lambda}(1).
\end{align*}
\end{lem}

\begin{lem} \label{gj2}
Let $\alpha,k,l,m,n \in \mathbb{N}$.
Provided ${\tau \Delta x}(\delta T^{2})^{-1} \le 1$, we have
\begin{align*}{\bf A}_{x}^{k}{\bf D}_{x}^{l}\partial_{x}^{\alpha}(r{\bf A}_{x}^{m}{\bf D}_{x}^{n}\rho)
=\partial_{x}^{l+\alpha}(r\partial_{x}^{n}\rho)+s^{n}\mathcal O_{\lambda}\left((s\Delta x)^{2}\right)
=s^{n}{\mathcal O}_{\lambda}(1).
\end{align*}
\end{lem}

\begin{lem} \label{gj3}
Let $\alpha,\beta,i,j,k,l,m,n \in \mathbb{N}$.
Provided ${\tau \Delta x}{(\delta T^{2})}^{-1} \le 1$, we have
\begin{align*}
&{\bf A}_{x}^{m}{\bf D}_{x}^{n}\partial_{x}^{\beta}(r^{2}{\bf A}_{x}^{i}{\bf D}_{x}^{j}(\partial_{x}^{\alpha}\rho){\bf A}_{x}^{k}{\bf D}_{x}^{l}\rho)\\
=\ &\partial_{x}^{n+\beta}(r^{2}\partial_{x}^{j+\alpha}\rho\partial_{x}^{l}\rho)
+s^{j+l+\alpha}\mathcal O_{\lambda}\left((s\Delta x)^{2}\right)=s^{j+l+\alpha}\mathcal O_{\lambda}(1).
\end{align*}
\end{lem}

The following results show the effect of the time-discrete operator over
some discrete operations in the space variable, whose proof can be found in  [\ref{Boyer-ESAIM:COCV-2005}].

\begin{lem} \label{sjgj1}
Provided ${\tau \Delta t}{\left(\delta^{2} T^{3}\right)}^{-1} \le 1$, we have
\begin{align*}{\bf t}^{-}(r) {\bf D}_{t}\rho
=-\tau {\bf t}^{-}\left(\theta'\right)\varphi
+\Delta t \left(\frac{\tau}{\delta^{3}T^{4}}+\frac{\tau^{2}}{\delta^{4}T^{6}}\right)\mathcal O_{\lambda}(1).
\end{align*}
\end{lem}

\begin{lem} \label{sjgj2}
There exists a constant $C>0$ independent of $\Delta t$, $\delta$ and $T$ such that 

\vspace{2mm}

\noindent$\begin{aligned}
({\rm i})&\  |{\bf D}_{t}(\theta^{l})| \le l T {\bf t}^{-}(\theta^{l+1})
+C\frac{\Delta t}{\delta^{l+2}T^{2l+2}},
\quad l=1,2,\cdots \\
({\rm ii})& \ {\bf D}_{t}(\theta')
\le CT^{2}{\bf t}^{-}(\theta^{3})
+C\frac{\Delta t}{\delta^{4}T^{5}}.
\end{aligned}$
\end{lem}

The next result is related to time-discrete operations performed to the Carleman weights.

\begin{lem} \label{sjgj3} Provided ${\tau \Delta x}{\left(\delta T^{2}\right)^{-1}} \le 1$ and $\tau \Delta t\left(\delta^{2}T^{3}\right)^{-1} \le {1}/{2}$, we have
\begin{align*}
{\bf D}_{t}(r{\bf D}_{x}^{2}\rho)
={\bf t}^{-}(\sigma_1),\quad {\bf D}_{t}(r{\bf A}_{x}^{2}\rho)
={\bf t}^-(\sigma_2), \quad {\bf D}_{t}(r{\bf A}_{x}{\bf D}_{x}\rho)={\bf t}^{-}(\sigma_3),
\end{align*}
where
\begin{align*}
&\sigma_1=\left(Ts^{2}\theta +\frac{\tau^{2}\triangle t}{\delta^{4}T^{6}}+\frac{\tau\triangle t}{\delta^{3}T^{4}}\left(\frac{\tau \Delta x}{\delta T^{2}}\right)^{3}\right)\mathcal O_{\lambda}(1),\\
&\sigma_2=\left(T(s\Delta x)^{2}\theta+\left(\frac{\tau \Delta t}{\delta^{3}T^{4}}\right)\left(\frac{\tau \Delta x}{\delta T^{2}}\right)\right)\mathcal O_{\lambda}(1),\\
&\sigma_3=\left(Ts\theta+Ts(s\Delta x)^{2}\theta
+\Delta tT^{2}s\theta^{2}+\Delta tT^{2}s(s\Delta x)^{2}\theta^{2}\right)\mathcal O_{\lambda}(1).
\end{align*}

\end{lem}


\section{Proof of Theorem 1.1}

This section is devoted to proving Carleman estimate for the full-discrete Ginzburg-Landau operator,
i.e. Theorem \ref{carleman estimate}. The proof follows as close as possible the ideas presented
in the classical continuous setting (see e.g. [\ref{Carreno-arXiv-2023}, \ref{Fu-JFA-2009}, \ref{Rosier-Poincare-2009}]), where Carleman estimates for the complex Ginzburg-Landau equations are obtained in the continuous setting. 

In comparison with the continuous setting, additional terms on $\partial \mathcal N$ appear in the discrete Carleman estimate, due to the regular time component in the weight function is introduced in this context. Moreover, in accordance with the discrete framework of our problem, we will particularly focus on the dependence concerning the discrete parameters $\Delta x$ and $\Delta t$.

\vspace{2mm}

\noindent{\bf Proof of Theorem 1.1.}\ We split the proof into the following several steps.

{\em Step 1.}\ {\em The change of variable.}

For $q \in C(\overline{\mathcal M} \times \overline{\mathcal{N}^{*}})$, we introduce the change of variables $z=rq=\rho^{-1} q$.
By a direct calculation, we can easily see that
\begin{align*}&
{\bf D}_{x}^{2}q = {\bf D}_{x}^{2}\rho {\bf A}_{x}^{2}z + 2{\bf A}_{x}{\bf D}_{x}\rho {\bf A}_{x}{\bf D}_{x}z + {\bf A}_{x}^{2}\rho {\bf D}_{x}^{2}z,\\
&{\bf D}_{t}q = {\bf t}^{-}(\rho){\bf D}_{t}z + {\bf D}_{t}\rho{\bf t}^{+}(z).
\end{align*}
Then, we obtain
\begin{align*}
{\bf t}^- (r){\mathcal P}(q)=&-{\bf D}_t z-{\bf t}^-(r){\bf D}_t\rho {\bf t}^+(z)-\\
&(\alpha-i\beta){\bf t}^-(r{\bf D}_{x}^{2}\rho {\bf A}_{x}^{2}z + 2r{\bf A}_{x}{\bf D}_{x}\rho {\bf A}_{x}{\bf D}_{x}z + r{\bf A}_{x}^{2}\rho {\bf D}_{x}^{2}z).
\end{align*}
Together with ${\bf t}^+ (z)={\bf t}^-(z)+\Delta t {\bf D}_t (z)$ and Lemma \ref{sjgj1}, we have
\begin{align}\label{3.1}\mathcal A(z)+\mathcal B(z)=\mathcal R(z),\quad (x,t)\in \mathcal M\times \mathcal N,
\end{align}
where 
\begin{align*} {\mathcal A}(z)=\sum_{k=1}^3 {\mathcal A}_{k}(z), \quad {\mathcal B}(z)=\sum_{k=1}^3{\mathcal B}_{k}(z)
\end{align*}
with
\begin{align*}&{\mathcal A}_{1}(z)=-\alpha {\bf t}^{-}(r{\bf D}_{x}^{2}\rho {\bf A}_{x}^{2}z + r{\bf A}_{x}^{2}\rho {\bf D}_{x}^{2}z),\\
&\mathcal A_{2}(z)=i \beta {\bf t}^{-}(2r{\bf A}_{x}{\bf D}_{x}\rho {\bf A}_{x}{\bf D}_{x}z - s\partial_{x}^2\phi z),\\
&{\mathcal A}_{3}(z)=\tau \varphi  {\bf t}^{-}(\theta' z),\\
& {\mathcal B}_{1}(z)=-\alpha {\bf t}^{-}(2r{\bf A}_{x}{\bf D}_{x}\rho {\bf A}_{x}{\bf D}_{x}z - s\partial_{x}^2\phi z),\\
&{\mathcal B}_{2}(z)=i \beta {\bf t}^{-}(r{\bf D}_{x}^{2}\rho {\bf A}_{x}^{2}z + r{\bf A}_{x}^{2}\rho {\bf D}_{x}^{2}z),\\
&{\mathcal B}_{3}(z)=-{\bf D}_{t}z
\end{align*}
and
\begin{align*}\mathcal R (z)
= &\ {\bf t}^{-}(r)\mathcal{P}(q) + \Delta t\left(\frac{\tau}{\delta^{3}T^{4}}+\frac{\tau^{2}}{\delta^{4}T^{6}}\right)\mathcal O_{\lambda}(1){\bf t}^{+}(z)-\\
&\Delta t\tau{\bf t}^{-}(\theta')\varphi {\bf D}_{t}z + 2(\alpha-i\beta){\bf t}^{-}(s\partial_{x}^2\phi z).
\end{align*}

On the other hand,  under the change of variable $z=rq$, following the similar steps to obtain (\ref{3.1}) we rewrite the dynamic boundary as
\begin{align}&\mathcal C (z)+\mathcal D(z)={\mathcal R}_{\Gamma_{0}}(z), \quad x=0, t\in \mathcal{N},\\
\label{3.3}&\mathcal E (z)+\mathcal F(z)={\mathcal R}_{\Gamma_{1}}(z), \quad x=1, t\in\mathcal{N},
\end{align}
where
\begin{align*}
&\mathcal C(z)
=-\alpha {\bf t}^{-}\left(r{\bf s}_{+}({\bf D}_{x}\rho {\bf A}_{x}z)\right)
+i\beta {\bf t}^{-}(r{\bf s}_{+}({\bf A}_{x}\rho {\bf D}_{x}z))
-{\bf D}_{t}z,\\
& \mathcal D(z)
=i\beta {\bf t}^{-}(r{\bf s}_{+}({\bf D}_{x}\rho {\bf A}_{x}z))
-\alpha {\bf t}^{-}(r{\bf s}_{+}({\bf A}_{x}\rho {\bf D}_{x}z))
+\tau\varphi {\bf t}^{-}(\theta'z),\\
&\mathcal E (z)
=\alpha {\bf t}^{-}(r{\bf s}_{-}({\bf D}_{x}\rho {\bf A}_{x}z))
-i\beta {\bf t}^{-}(r{\bf s}_{-}({\bf A}_{x}\rho {\bf D}_{x}z))
-{\bf D}_{t}z,\\
&\mathcal F(z)
=-i\beta {\bf t}^{-}(r{\bf s}_{-}({\bf D}_{x}\rho {\bf A}_{x}z))
+\alpha {\bf t}^{-}(r{\bf s}_{-}({\bf A}_{x}\rho {\bf D}_{x}z))
+\tau\varphi{\bf t}^{-}(\theta' z)
\end{align*}
and
\begin{align*}
&{\mathcal R}_{\Gamma_{0}}(z)
=-{\bf t}^{-}(r){\mathcal B}_{\Gamma_{0}}(q)
+ \Delta t\left(\frac{\tau}{\delta^{3}T^{4}}+\frac{\tau^{2}}{\delta^{4}T^{6}}\right)\mathcal O_{\lambda}(1){\bf t}^{+}(z)
- \Delta t\tau\varphi {\bf t}^{-}(\theta'){\bf D}_{t}z,\\
& {\mathcal R}_{\Gamma_{1}}(z)
=-{\bf t}^{-}(r){\mathcal B}_{\Gamma_{1}}(q)
+\Delta t\left(\frac{\tau}{\delta^{3}T^{4}}+\frac{\tau^{2}}{\delta^{4}T^{6}}\right)\mathcal O_{\lambda}(1) {\bf t}^{+}(z)
-\Delta t\tau\varphi{\bf t}^{-}(\theta') {\bf D}_{t}z.
\end{align*}

From (\ref{3.1})-(\ref{3.3}), we conclude that
\begin{align} \label{dengshi}
&||{\mathcal A}(z)||^{2}_{L^{2}(\mathcal{M} \times \mathcal{N})} + ||{\mathcal B}(z)||^{2}_{L^{2}(\mathcal{M} \times \mathcal{N})}
+2 \text{Re} ( {\mathcal A}(z),{\mathcal B}(z) )_{L^{2}(\mathcal{M} \times \mathcal{N})} +\nonumber \\
&||{\mathcal C}(z)||^{2}_{L^{2}(\{0\} \times \mathcal{N})}+||{\mathcal D}(z)||^{2}_{L^{2}(\{0\} \times \mathcal{N})}
+2 \text{Re} ( {\mathcal C}(z),{\mathcal D}(z) )_{L^{2}(\{0\} \times \mathcal{N})} +\nonumber \\
&||{\mathcal E}(z)||^{2}_{L^{2}(\{1\} \times \mathcal{N})}+||{\mathcal F}(z)||^{2}_{L^{2}(\{1\} \times \mathcal{N})}
+2 \text{Re} ( {\mathcal E}(z),{\mathcal F}(z) )_{L^{2}(\{1\} \times \mathcal{N})} \nonumber \\
=\ &||{\mathcal R}(z)||^{2}_{L^{2}(\mathcal{M} \times \mathcal{N})}
+||{\mathcal R}_{\Gamma_{0}}(z)||^{2}_{L^{2}(\{0\} \times \mathcal{N})}
+||{\mathcal R}_{\Gamma_{1}}(z)||^{2}_{L^{2}(\{1\} \times \mathcal{N})}.
\end{align}


{\em Step 2. Computations of integral terms in $\mathcal M\times \mathcal N$.}

For clarity, we set
\begin{align}\label{3.5}2 \text{Re} ( {\mathcal A}(z),{\mathcal B}(z) )_{L^{2}(\mathcal{M} \times \mathcal{N})}=\sum_{i,j=1}^32\text{Re} \iint\limits_{\mathcal{M} \times \mathcal{N}} {\mathcal A}_{i}(z){\mathcal B}^{*}_{j}(z)=:\sum_{i,j=1}^3I_{ij},
\end{align}
where ${\mathcal B}^{*}_{j}(z)$ denotes the conjugate of ${\mathcal B}_{j}(z)$. Also, in the following, we use $BT_{ij}$ to denote the boundary term generated by $I_{ij}$.

Obviously, we see that $I_{12}=I_{21}=0$. The other terms in (\ref{3.5}) are estimated one by one in Appendix A, see Lemma A.1-A.7. For fixed $\lambda\geq \mathcal O(1)$,  we can choose sufficiently large $\tau_0\geq 1$  and sufficiently small $0<\varepsilon_0<1$
depending on $\alpha, \beta, \delta, \omega$, $T$, $c_0$ and $\lambda$, but independent of $\Delta x$ and $\Delta t$ to satisfy
\begin{align*}&s\geq\mathcal O_{\lambda}(1), \quad (s\Delta x)\mathcal O_{\lambda}(1)\leq \varepsilon_0 \mathcal O_{\lambda}(1)\leq \epsilon, \quad \left(s^4\Delta t\right)\mathcal O_{\lambda}(1)\leq\varepsilon_0 \mathcal O_{\lambda}(1)\leq \epsilon, \\
&\sigma_1\leq (s^2+\tau \varepsilon_0+\Delta x\varepsilon_0^3)\mathcal O_{\lambda}(1)\leq s^3, \quad \sigma_2\leq(\varepsilon_0^2+\Delta x\varepsilon_0) \mathcal O_{\lambda}(1)\leq \epsilon,\\ 
&\sigma_3\leq (s+\varepsilon_0^2s+s\Delta t+\varepsilon_0^2s\Delta t)\mathcal O_{\lambda}(1)\leq s\mathcal O_{\lambda}(1)
\end{align*}
with sufficiently small $\epsilon>0$, if
$$\tau \ge \tau_0(T+T^{2}), \quad \frac{\tau \Delta x}{\delta T^{2}} \le \varepsilon_0,\quad \frac{\tau^4 \Delta  t}{\delta^{4}T^{6}} \le \varepsilon_0.$$ 
Then we have 
\begin{align*}
&\Delta t-{\bf t}^{-}\left(\Delta t(s\Delta x)^{2}+(\Delta t)^{2}\sigma_{2}+(\Delta x)^2(\Delta t)^3\sigma_1\right)\geq (1-3\epsilon)\Delta t,\\
&s^{2}\Delta t+(\Delta t)^{2}\sigma_{1}+s^{-1}(s\Delta x)^{2}\leq s^{-1},\quad (\Delta x)^2\sigma_1 +s(s\Delta x)^{2}\mathcal O_{\lambda}(1)+\sigma_{2}\leq s,\\
&s^{-1}(s\Delta x)^{2}\mathcal O_{\lambda}(1)+(\Delta t)^{2}s\sigma_{3}+C(\epsilon)s^{2}\Delta  t\mathcal O_{\lambda}(1)\leq s^{-1},\\
&\epsilon \Delta  t+s\Delta t\Delta x\mathcal O_{\lambda}(1)+s^{-1}(\Delta  t)^{2}\sigma_{3}\leq 3\epsilon\Delta t.
\end{align*}
Then choosing $\epsilon$ such that  $\alpha(1-3\epsilon)-3\epsilon\beta\geq 0$, combining Lemma A.1-A.7 and (\ref{3.5}), and using the property of $\psi$ (\ref{psi}), we find that
\begin{align}\label{1-3.6}
&2 \text{Re} ( \mathcal A(z),\mathcal B(z) )_{L^{2}(\mathcal{M} \times \mathcal{N})} \nonumber \\
\ge\ &2(\alpha^{2}+\beta^{2}) \iint\limits_{\mathcal{M} \times \mathcal{N}^{*}} s^{3}\lambda^4\phi^3|z|^{2}
+2(\alpha^{2}+\beta^{2})\iint\limits_{\mathcal{M}^{*} \times \mathcal{N}^{*}} s \lambda^2\phi|{\bf D}_{x}z|^{2}- \iint\limits_{\mathcal{\omega}_0 \times \mathcal{N}^{*}}s^{3} \mathcal O_\lambda(1) |z|^{2}
-\nonumber \\
&\iint\limits_{\mathcal{\omega}_0^{*} \times \mathcal{N}^{*}} s \mathcal O_\lambda(1)  |{\bf D}_{x}z|^{2} - \iint\limits_{\mathcal{M} \times \mathcal{N}^{*}} s\mathcal O(1)|{\bf A}_{x}{\bf D}_{x}z|^{2}
-\epsilon\iint\limits_{\mathcal{M} \times \mathcal{N}^{*}} s^{-1}\mathcal |{\bf D}_{x}^{2}z|^{2}-
\nonumber \\
&\epsilon\iint\limits_{\mathcal{M} \times \mathcal{N}} {\bf t}^{-}(s^{-1}) |{\bf D}_{t}z|^{2}
+BT,
\end{align}
where
\begin{align*}
BT
\ge& \frac{1}{2}(\alpha^2+\beta^2) c_1 \iint\limits_{\partial\mathcal{M} \times \mathcal{N}^{*}} s^{3}\lambda^{3}\phi^3|z|^{2}
+ \frac{1}{2}(\alpha^2+\beta^2) c_1 \iint\limits_{\partial\mathcal{M} \times \mathcal{N}^{*}} s\lambda \phi{\bf tr}(|{\bf D}_{x}z|^{2})
- \nonumber \\
&\epsilon\iint\limits_{\partial\mathcal{M} \times \mathcal{N}} {\bf t}^{-}(s^{-1}) |{\bf D}_{t}z|^{2}-\iint\limits_{\mathcal{M} \times \partial\mathcal{N}}{\bf t}^{+}(s^2)\mathcal O_{\lambda}(1) {\bf t}^{+}(|z|^{2})
- \nonumber \\
&\iint\limits_{\mathcal{M}^{*} \times \partial\mathcal{N}} {\bf t}^{+}(s^2) \mathcal O_{\lambda}(1) {\bf t}^{+}(|{\bf A}_{x}z|^{2})-\iint\limits_{\mathcal{M}^{*} \times \partial\mathcal{N}} \mathcal O_{\lambda}(1) {\bf t}^+(|{\bf D}_{x}z|^{2}).
\end{align*}

{\em Step 3. Computations of boundary terms on $\partial \mathcal M\times \mathcal N$.}

We write
$$2 \text{Re} ( \mathcal C(z),\mathcal D(z) )_{L^{2}(\{0\} \times \mathcal{N})}=
\sum_{i,j=1}^32\text{Re} \iint\limits_{\{0\} \times \mathcal{N}} {\mathcal C}_{i}(z){\mathcal D}^{*}_{j}(z)=:\sum_{i,j=1}^{3}J_{ij}.$$

Obviously, we have $J_{11}=J_{22}=0$. Using
$${\bf s}_{+}({\bf A}_{x}z)
=\frac{1}{2}\Delta x{\bf s}_{+}({\bf D}_{x}z)+z
\quad
\text{in} \; \{0\} \times \mathcal{N}^{*},$$
we obtain
\begin{align*}
|J_{12}|
\leq \ &\alpha^{2}\left|\ \iint\limits_{\{0\} \times \mathcal{N}^{*}} r^{2}{\bf s}_{+}({\bf A}_{x}\rho {\bf D}_{x}\rho)\left( \frac{1}{2}\Delta x{\bf s}_{+}({\bf D}_{x}z)+z\right){\bf s}_+({\bf D}_{x}z^*)\right|.
\end{align*}
Since $$r^2(x,t)=e^{2s(t)\varphi(x)}=e^{2s(t)\varphi(x+\frac{h}{2})+\mathcal O_{\lambda}(s\Delta x)}=e^{(s\Delta x)\mathcal O_{\lambda}(1)}{\bf s}_+(r^2),$$
we have
\begin{align*}
r^{2}{\bf s}_{+}({\bf A}_{x}\rho {\bf D}_{x}\rho)=e^{\mathcal O_{\lambda}(s\Delta x)}{\bf s}_{+}(r^2{\bf A}_{x}\rho {\bf D}_{x}\rho)=e^{\mathcal O_{\lambda}(s\Delta x)}s\mathcal O_{\lambda}(1)=s\mathcal O_{\lambda}(1).
\end{align*}
Then we further obtain the estimate for $J_{12}$
\begin{align}\label{3.6}
J_{12}
\ge&-\alpha^2\iint\limits_{\{0\} \times \mathcal{N}^{*}} s^{2}\mathcal O_{\lambda}(1) |z|^{2}
-\alpha^2\iint\limits_{\{0\} \times \mathcal{N}^{*}} (1+(s\Delta x)) \mathcal O_{\lambda}(1){\bf tr}(|{\bf D}_{x}z|^{2}).
\end{align}
We easily see that $J_{21}=\beta^2/\alpha^2 J_{12}$. Therefore,  
\begin{align}
J_{21}
\ge\ &-\beta^2\iint\limits_{\{0\} \times \mathcal{N}^{*}} s^{2}\mathcal O_{\lambda}(1) |z|^{2}
-\beta^2\iint\limits_{\{0\} \times \mathcal{N}^{*}} (1+(s\Delta x))\mathcal O_{\lambda}(1)) {\bf tr}(|{\bf D}_{x}z|^{2}).
\end{align}
By a process similar to (\ref{3.6}), we have
\begin{align}
J_{13}
=&\ -2\alpha \text{Re} \iint\limits_{\{0\} \times \mathcal{N}^{*}} \tau\theta{'}\varphi r{\bf s}_{+}({\bf D}_{h}\rho)  {\bf s}_{+}({\bf A}_{x}z)z^{*} \nonumber \\
=&\ -2\alpha \text{Re} \iint\limits_{\{0\} \times \mathcal{N}^{*}} \tau\theta{'}s\mathcal O_{\lambda}(1)  \left( \frac{1}{2}\Delta x{\bf s}_{+}({\bf D}_{x}z)+z\right)z^{*}\nonumber \\
\ge&-\alpha\iint\limits_{\{0\} \times \mathcal{N}^{*}} Ts^{2}\theta \mathcal O_{\lambda}(1) |z|^{2}
-\alpha\iint\limits_{\{0\} \times \mathcal{N}^{*}} T(s\Delta x)^{2}\theta \mathcal O_{\lambda}(1) {\bf tr}(|{\bf D}_{x}z|^{2})
\end{align}
and
\begin{align}
J_{23}
=\ & 2\beta \text{Re} \iint\limits_{\{0\} \times \mathcal{N}^{*}} i\tau\theta{'}\varphi r {\bf s}_{+} ({\bf A}_{x}\rho) {\bf s}_{+} ({\bf D}_{x}z) z^{*} \nonumber \\
\ge \ &-\beta\iint\limits_{\{0\} \times \mathcal{N}^{*}} T^{2}s^{2}\theta^{2}\mathcal O_{\lambda}(1) |z|^{2}
-\beta\iint\limits_{\{0\} \times \mathcal{N}^{*}} \mathcal O_{\lambda}(1) {\bf tr}(|{\bf D}_{x}z|^{2}).
\end{align}
Using 
\begin{align*}
r{\bf s}_{+}({\bf D}_x\rho)=e^{\mathcal O_{\lambda}(s\Delta x)}{\bf s}_{+}(r{\bf D}_{x}\rho)=s\lambda\phi\mathcal O(1)+s(s\Delta x)^2\mathcal O_{\lambda}(1)
\end{align*}
and Young's inequality with $\epsilon$, we deduce that
\begin{align}
J_{31}
=\ &2\beta \text{Re} \iint\limits_{\{0\} \times \mathcal{N}} i {\bf t}^{-}(r{\bf s}_{+}({\bf D}_{x}\rho)) {\bf t}^{-}({\bf s}_{+}({\bf A}_{x}z^{*})) {\bf D}_{t}z \nonumber \\
=\ &2\beta \text{Re} \iint\limits_{\{0\} \times \mathcal{N}} i {\bf t}^{-}(s\lambda\phi\mathcal O(1)+s(s\Delta x)^2\mathcal O_{\lambda}(1)) {\bf t}^{-}\left(\frac{1}{2}\Delta x{\bf s}_{+}({\bf D}_{x}z)+z\right) {\bf D}_{t}z\nonumber\\
\geq \ &-C(\epsilon)\beta\iint\limits_{\{0\} \times \mathcal{N}^{*}} s^{3}\lambda^{2}\phi^{2}|z|^{2}
-\beta\iint\limits_{\{0\} \times \mathcal{N}^{*}} s(s\Delta x) \mathcal O_{\lambda}(1) {\bf tr}(|{\bf D}_{x}z|^{2}) -\nonumber \\
&\beta\iint\limits_{\{0\} \times \mathcal{N}} {\bf t}^{-}(\epsilon s^{-1}+s^{-1}(s\Delta x)\mathcal O_{\lambda}(1))|{\bf D}_{t}z|^{2}.
\end{align}
Moreover, by using $r{\bf s}_+({\bf A}_x \rho)=(1+\mathcal O_{\lambda}((s\Delta x)^2))\mathcal O(1)$ we obtain for $J_{32}$ that
\begin{align}
J_{32}
=\ &2\alpha \text{Re} \iint\limits_{\{0\} \times \mathcal{N}} {\bf t}^{-}(r{\bf s}_{+}({\bf A}_{x}\rho)) {\bf t}^{-}({\bf s}_{+}({\bf D}_{x}z^{*})) {\bf D}_{t}z \nonumber \\
\ge\ &-\alpha\iint\limits_{\{0\} \times \mathcal{N}^{*}} s\mathcal O(1) {\bf tr}(|{\bf D}_{x}z|^{2})
-\epsilon \alpha\iint\limits_{\{0\} \times \mathcal{N}} {\bf t}^{-}(s^{-1}) |{\bf D}_{t}z|^{2}.
\end{align}
Finally, we integrate by parts wit respect to ${\bf D}_t$ to obtain
\begin{align}\label{3.12}
J_{33}=\ & -{\rm Re}\iint\limits_{\{0\}\times \mathcal N} \tau \varphi {\bf t}^{-}(\theta'){\bf t}^{-}(z^*){\bf D}_t z\nonumber\\
=\ & \iint\limits_{\{0\} \times \mathcal{N}} \tau\varphi {\bf D}_{t}(\theta') {\bf t}^{+}(|z|^{2})
+\Delta t \iint\limits_{\{0\} \times \mathcal{N}} \tau  {\bf t}^{-}(\theta')\varphi |{\bf D}_{t}z|^{2}
-\iint\limits_{\{0\} \times \partial\mathcal{N}} \tau{\bf t}^{+}(\theta')\varphi {\bf t}^+(|z|^{2})n_{t} \nonumber \\
\ge\ &-\iint\limits_{\{0\} \times \mathcal{N}^{*}} \left(T^{2}s\theta^{2} + \left(\frac{\tau\Delta  t}{\delta^{4}T^{5}}\right)\right)\mathcal O_{\lambda}(1) |z|^{2}
-\iint\limits_{\{0\} \times \mathcal{N}} {\bf t}^{-}(T(s\Delta t)\theta ) \mathcal O_{\lambda}(1) |{\bf D}_{t}z|^{2},
\end{align}
where we have used 
$$-\iint\limits_{\{0\} \times \partial\mathcal{N}} \tau{\bf t}^{+}(\theta')\varphi {\bf t}^+(|z|^{2})n_{t}\geq 0.$$

By the above estimates (\ref{3.6})-(\ref{3.12}) for $J_{ij}$, we obtain the following estimate for boundary terms on $\{0\} \times \mathcal{N}$:
\begin{align}\label{3.13}
&2 \text{Re} ( \mathcal C(z),\mathcal D(z) )_{L^{2}(\{0\} \times \mathcal{N})} \nonumber \\
\ge\ &-C(\epsilon)\iint\limits_{\{0\} \times \mathcal{N}^{*}} s^{3}\lambda^{2}\phi^{2}\mathcal O(1) |z|^{2}
-\iint\limits_{\{0\} \times \mathcal{N}^{*}} s\mathcal O(1) {\bf tr}(|{\bf D}_{x}z|^{2})-\epsilon \iint\limits_{\{0\} \times \mathcal{N}} {\bf t}^{-}(s^{-1}) |{\bf D}_{t}z|^{2}.
\end{align}
By a process similar to (\ref{3.13}),  we can obtain a lower bound for boundary terms on $\{1\} \times \mathcal{N}$
\begin{align}
&2 \text{Re} ( \mathcal E(z),\mathcal F(z) )_{L^{2}(\{1\} \times \mathcal{N})} \nonumber \\
\ge\ &-C(\epsilon)\iint\limits_{\{1\} \times \mathcal{N}^{*}} s^{3}\lambda^{2}\phi^{2}\mathcal O(1) |z|^{2}
-\iint\limits_{\{1\} \times \mathcal{N}^{*}}s\mathcal O (1) {\bf tr}(|{\bf D}_{x}z|^{2})-\epsilon \iint\limits_{\{1\} \times \mathcal{N}} {\bf t}^{-}(s^{-1} )|{\bf D}_{t}z|^{2}.
\end{align}

{\em Step 4. Computations of remainder terms.}

Using H\"{o}lder inequality and ${\bf t}^+(|z|^2)\leq C{\bf t}^-(|z|^2)+C(\Delta t)^2|{\bf D}_t z|^2$, we can easily obtain that
\begin{align} \label{rz}
&\|\mathcal R(z)\|^{2}_{L^{2}(\mathcal{M} \times \mathcal{N})}\nonumber\\
\le\ & \iint\limits_{\mathcal{M} \times \mathcal{N}} {\bf t}^{-}(r^{2})|\mathcal{P}(q)|^{2}
+\iint\limits_{\mathcal{M} \times \mathcal{N}^{*}} \left(s^{2}+\left(\frac{\tau\Delta t}{\delta^{3}T^{4}}\right)^{2}+\left(\frac{\tau^{2}\Delta t}{\delta^{4}T^{6}}\right)^{2}\right)\mathcal O_{\lambda}(1) |z|^{2} +\nonumber \\
&\iint\limits_{\mathcal{M} \times \mathcal{N}} {\bf t}^{-}\left((s\Delta t )^{2}\right)\mathcal O_{\lambda}(1)  |{\bf D}_{t}z|^{2}\nonumber\\
\le\ & \iint\limits_{\mathcal{M} \times \mathcal{N}} {\bf t}^{-}(r^{2})|\mathcal{P}(q)|^{2}
+\iint\limits_{\mathcal{M} \times \mathcal{N}^{*}} s^{2}\mathcal O_{\lambda}(1) |z|^{2} +\epsilon\iint\limits_{\mathcal{M} \times \mathcal{N}} {\bf t}^{-}(s^{-1}) |{\bf D}_{t}z|^{2}.
\end{align}
Similarly, 
\begin{align} \label{rz0}
&\|\mathcal R_{\Gamma_{0}}(z)\|^{2}_{L^{2}(\{0\} \times \mathcal{N})}\nonumber\\
\le\ & \iint\limits_{\{0\} \times \mathcal{N}} {\bf t}^{-}(r^{2})|\mathcal B_{\Gamma_{0}}(q)|^{2}
+\iint\limits_{\{0\} \times \mathcal{N}^{*}}\mathcal O(1) |z|^{2}+\epsilon\iint\limits_{\{0\} \times \mathcal{N}} {\bf t}^{-}(s^{-1}) |{\bf D}_{t}z|^{2},
\end{align}
and
\begin{align} \label{rz1}
&\|{\mathcal R}_{\Gamma_{1}}(z)\|^{2}_{L^{2}(\{M+1\} \times \mathcal{N})}\nonumber\\
\le\ & \iint\limits_{\{1\} \times \mathcal{N}} {\bf t}^{-}(r^{2})|\mathcal B_{\Gamma_{1}}(q)|^{2}
+\iint\limits_{\{1\} \times \mathcal{N}^{*}}\mathcal O(1) |z|^{2}+\epsilon\iint\limits_{\{1\} \times \mathcal{N}} {\bf t}^{-}(s^{-1}) |{\bf D}_{t}z|^{2}.
\end{align}

From (\ref{dengshi}), (\ref{1-3.6}) and (\ref{3.13})-(\ref{rz1}), it follows that 
\begin{align}\label{3.19}
&\|\mathcal A(z)\|^{2}_{L^{2}(\mathcal{M} \times \mathcal{N})} + \|\mathcal B(z)\|^{2}_{L^{2}(\mathcal{M} \times \mathcal{N})}+\|\mathcal C(z)\|^{2}_{L^{2}(\{0\} \times \mathcal{N})}+\|\mathcal D(z)\|^{2}_{L^{2}(\{0\} \times \mathcal{N})}+ 
\nonumber\\
&\|\mathcal E(z)\|^{2}_{L^{2}(\{1\} \times \mathcal{N})}+\|\mathcal F(z)\|^{2}_{L^{2}(\{1\} \times \mathcal{N})} +(\alpha^2+\beta^2)\times\nonumber \\
& \left(\ \iint\limits_{\mathcal{M} \times \mathcal{N}^{*}} s^{3}\lambda^{4}\phi^{3} |z|^{2}
+\iint\limits_{\mathcal{M}^{*} \times \mathcal{N}^{*}} s\lambda^{2}\phi |{\bf D}_{x}z|^{2} 
+ \iint\limits_{\partial\mathcal{M} \times \mathcal{N}^{*}} s^{3}\lambda^{3}\phi^{3} |z|^{2}
+ \iint\limits_{\partial\mathcal{M} \times \mathcal{N}^{*}} s\lambda\phi {\bf tr}(|{\bf D}_{x}z|^{2})\right) \nonumber \\
\le\ &\iint\limits_{\mathcal{M} \times \mathcal{N}} {\bf t}^{-}(r^{2})|\mathcal{P}(q)|^{2}
+\iint\limits_{\{0\} \times \mathcal{N}} {\bf t}^{-}(r^{2})|\mathcal B_{\Gamma_{0}}(q)|^{2}
+\iint\limits_{\{1\} \times \mathcal{N}} {\bf t}^{-}(r^{2})|\mathcal B_{\Gamma_{1}}(q)|^{2} +\nonumber \\
&\epsilon \left(\ \iint\limits_{\mathcal{M} \times \mathcal{N}} {\bf t}^{-}(s^{-1}) |{\bf D}_{t}z|^{2}
+\iint\limits_{\mathcal{M} \times \mathcal{N}^{*}}s^{-1} |{\bf D}_{x}^{2}z|^{2}+\iint\limits_{\partial\mathcal{M} \times \mathcal{N}} {\bf t}^{-}(s^{-1})|{\bf D}_{t}z|^{2}\right)+ \nonumber \\
&C({\lambda})\left(\ \iint\limits_{\mathcal{\omega}_{0} \times \mathcal{N}^{*}} s^{3}|z|^{2}
+\iint\limits_{\mathcal{\omega}_{0}^{*} \times \mathcal{N}^{*}} s|{\bf D}_{x}z|^{2}\right)
+C({\lambda})\iint\limits_{\mathcal{M} \times \mathcal{N}^{*}} |{\bf A}_{x}{\bf D}_{x}z|^{2}
+\nonumber \\
&C({\lambda}) \left(\ \iint\limits_{\mathcal{M} \times \partial\mathcal{N}} {\bf t}^{+}(s^{2}){\bf t}^{+}(|z|^{2})+\iint\limits_{\mathcal{M}^{*} \times \partial\mathcal{N}}{\bf t}^{+} (s^{2}) {\bf t}^{+}(|{\bf A}_{x}z|^{2})
+\iint\limits_{\mathcal{M}^{*} \times \partial\mathcal{N}} {\bf t}^{+}(|{\bf D}_{x}z|^{2})\right).
\end{align}

{\em Step 5. Estimates for $|{\bf A}_{x}{\bf D}_{x}z|^{2}$, $|{\bf D}_{x}^{2}z|^{2}$ and $|{\bf D}_{t}z|^{2}$.}

By Lemma 2.12 in [\ref{Carreno-ACM-2023}], we know that there exists sufficiently small $\widehat {\Delta x}<0$ such that 
\begin{align*}
&\iint\limits_{\mathcal{M} \times \mathcal{N}^{*}} s\lambda^{2} \phi |{\bf A}_{x}{\bf D}_{x}z|^{2}
+\frac{1}{4}(\Delta x)^{2} \iint\limits_{\mathcal{M} \times \mathcal{N}^{*}} s\lambda^{2}\phi |{\bf D}_{x}^{2}z|^{2}\nonumber\\
\leq\ 
&\iint\limits_{\mathcal{M}^{*} \times \mathcal{N}^{*}} s\lambda^{2} \phi |{\bf D}_{x}z|^{2}+(\Delta x)^{2} \iint\limits_{\mathcal{M}^{*} \times \mathcal{N}^{*}} s\mathcal O_{\lambda}(1) |{\bf D}_{x}z|^{2}
+\iint\limits_{\mathcal{M} \times \mathcal{N}^{*}} (s\Delta x)\mathcal O_{\lambda}(1) |{\bf A}_{x}{\bf D}_{x}z|^{2}+ \nonumber \\
&(\Delta x)^{4} \iint\limits_{\mathcal{M} \times \mathcal{N}^{*}} s\mathcal O_{\lambda}(1) |{\bf D}_{x}^{2}z|^{2},
\end{align*}
for all $0<\Delta x\leq \widehat {\Delta x}$, which implies 
\begin{align}\label{3.20}
&\iint\limits_{\mathcal{M} \times \mathcal{N}^{*}}s\lambda^{2} \phi |{\bf A}_{x}{\bf D}_{x}z|^{2}
\leq\ \iint\limits_{\mathcal{M}^{*} \times \mathcal{N}^{*}} s\lambda^{2} \phi |{\bf D}_{x}z|^{2}+\epsilon \iint\limits_{\mathcal{M} \times \mathcal{N}^{*}} s^{-1}|{\bf D}_{x}^{2}z|^{2}.
\end{align}

Using the definition of ${\mathcal A}(z)$, we have
\begin{align*}
&\alpha {\bf t}^-(1+(s\Delta x )^2\mathcal O_{\lambda}(1)){\bf t}^{-} ({\bf D}_{x}^{2}z)\nonumber\\
=&-\mathcal A (z)-\alpha {\bf t}^{-}(s^2\lambda^2\phi^2\mathcal O(1) ){\bf t}^{-}({\bf A}^2_x z)+i \beta {\bf t}^-(s)\left(\lambda\phi\mathcal O(1){\bf t}^-({\bf A}_x{\bf D}_xz)-\mathcal O_{\lambda}(1){\bf t}^-(z)\right) +\nonumber \\
&\tau \mathcal O_{\lambda}(1) {\bf t}^-(\theta'){\bf t}^-(z).
\end{align*}
Multiplying the above equation by ${\bf t}^{-}(s^{-1/2})\phi^{-1/2}$, integrating over $\mathcal M\times \mathcal N$ and using ${\bf A}_x^2 z=z+1/4(\Delta x)^2{\bf D}_x^2 z$, we find that
\begin{align*}
\alpha^2\iint\limits_{\mathcal{M} \times \mathcal{N}^{*}} s^{-1}\phi^{-1} |{\bf D}_{x}^{2}z|^{2}
\le&\iint\limits_{\mathcal{M} \times \mathcal{N}} {\bf t}^{-}(s^{-1}) |{\mathcal A}(z)|^{2}
+\alpha^2\iint\limits_{\mathcal{M} \times \mathcal{N}^{*}} s^{-1}(s\Delta x)^{4}\mathcal O_{\lambda}(1)|{\bf D}_{x}^{2}z|^{2}+\nonumber \\
&\iint\limits_{\mathcal{M} \times \mathcal{N}^{*}} s^{3}\lambda^4\phi^3\mathcal O(1) |z|^{2}
+\beta^2\iint\limits_{\mathcal{M} \times \mathcal{N}^{*}} s \lambda^2\phi\mathcal O(1)|{\bf A}_{x}{\bf D}_{x}z|^{2},
\end{align*}
which implies
\begin{align} \label{dh2z}
&\alpha^2\iint\limits_{\mathcal{M} \times \mathcal{N}^{*}} s^{-1}\phi^{-1} |{\bf D}_{x}^{2}z|^{2}\nonumber\\
\le& \ \|{\mathcal A}(z)\|^{2}_{L^2(\mathcal{M} \times \mathcal{N})}+\iint\limits_{\mathcal{M} \times \mathcal{N}^{*}} s^{3}\lambda^4\phi^3\mathcal O(1) |z|^{2}
+\beta^2\iint\limits_{\mathcal{M} \times \mathcal{N}^{*}} s\lambda^2\phi\mathcal O(1) |{\bf A}_{x}{\bf D}_{x}z|^{2},
\end{align}
if $s^{-1}(s\Delta x)^{4}\leq 1/2$. By the definition of $\mathcal B(z)$, we also obtain
\begin{align} \label{dtz}
\iint\limits_{\mathcal{M} \times \mathcal{N}} {\bf t}^{-}(s^{-1})\phi^{-1} |{\bf D}_{t}z|^{2}
\le\ & ||\mathcal B (z)||_{L^2(\mathcal{M} \times \mathcal{N})}^{2}
+\iint\limits_{\mathcal{M} \times \mathcal{N}^{*}} s^{3}\lambda^4\phi^3\mathcal O (1) |z|^{2}
+\nonumber\\
&\iint\limits_{\mathcal{M} \times \mathcal{N}^{*}} s\lambda^2\phi\mathcal O(1)  |{\bf A}_{x}{\bf D}_{x}z|^{2}
+\beta^2\iint\limits_{\mathcal{M} \times \mathcal{N}^{*}} s^{-1} \phi^{-1} \mathcal O(1)|{\bf D}_{x}^{2}z|^{2}.
\end{align}
Combining \eqref{dh2z} and \eqref{dtz}, we obtain 
\begin{align}\label{3.23}
&\alpha^2\iint\limits_{\mathcal{M} \times \mathcal{N}^{*}} s^{-1}\phi^{-1} |{\bf D}_{x}^{2}z|^{2}
+\iint\limits_{\mathcal{M} \times \mathcal{N}} {\bf t}^{-}(s^{-1})\phi^{-1} |{\bf D}_{t}z|^{2}\nonumber \\
\le\ &||\mathcal A(z)||_{L^{2}(\mathcal{M} \times \mathcal{N})}^{2}
+||\mathcal B(z)||_{L^{2}(\mathcal{M} \times \mathcal{N})}^{2}
+C\left(\ \iint\limits_{\mathcal{M} \times \mathcal{N}^{*}} s^{3}\lambda^4\phi^3 \mathcal |z|^{2}
+\iint\limits_{\mathcal{M} \times \mathcal{N}^{*}} s\lambda^2\phi |{\bf A}_{x}{\bf D}_{x}z|^{2}\right).
\end{align}
Similarly, by the definitions of $\mathcal C(z)$ and $\mathcal E(z)$, we obtain
%
\begin{align}\label{3.24}
\iint\limits_{\partial\mathcal{M} \times \mathcal{N}} {\bf t}^{-}(s^{-1})|{\bf D}_{t}z|^{2}
\le\ & \|\mathcal C(z)\|_{L^{2}(\{0\} \times \mathcal{N})}^{2}
+\|\mathcal E(z)\|_{L^{2}(\{1\} \times \mathcal{N})}^{2} +\nonumber \\
&C\left(\ \ \iint\limits_{\partial\mathcal{M} \times \mathcal{N^*}} s \mathcal O_{\lambda}(1)|z|^{2}
+\iint\limits_{\partial\mathcal{M} \times \mathcal{N^*}} s^{-1}\mathcal O_{\lambda}(1) {\bf tr}(|{\bf D}_{x}z|^{2})\right).
\end{align}

{\em Step 6. The complete of the proof.}

 In order to eliminate the local term of $|{\bf D}_{x}z|^{2}$ on the right-hand side of (\ref{3.19}), we need a local estimate for $z$ on $\mathcal{\omega}_{0}^{*} \times \mathcal{N}^{*}$, which can be found in [\ref{Carreno-ACM-2023}]. By Lemma 2.16 in [\ref{Carreno-ACM-2023}], we obtain that 
\begin{align} \label{jbdhz}
\iint\limits_{\mathcal{\omega}_{0}^{*} \times \mathcal{N}^{*}} s |{\bf D}_{x}z|^{2}
\le&C(\epsilon) \iint\limits_{\mathcal{\omega}\times \mathcal{N}^{*}} s^{3} |z|^{2}
+C \iint\limits_{\mathcal{M} \times \mathcal{N}^{*}} |{\bf A}_{x}{\bf D}_{x}z|^{2}
+C \iint\limits_{\mathcal{M} \times \mathcal{N}^{*}} s |z|^{2}+ \nonumber \\
&\epsilon \iint\limits_{\mathcal{M} \times \mathcal{N}^{*}} s^{-1} |{\bf D}_{x}^{2}z|^{2}
+C \iint\limits_{\mathcal{M} \times \mathcal{N}^{*}} s^{-1}(sh)^{2} |{\bf D}_{x}^{2}z|^{2}.
\end{align}
for any $\epsilon > 0$.

Now we choose $\epsilon$ such that $\epsilon\mathcal O_{\lambda}(1)\leq \alpha^2$ sufficiently small and substitute (\ref{3.20}), (\ref{3.23}) and (\ref{3.24}) into (\ref{3.19}) to eliminate the terms of ${\bf D}_x^2 z, {\bf D}_t z$ and ${\bf A}_x{\bf D}_x z$ on the right-hand side of (\ref{3.19}). The we find that
\begin{align}\label{3.26}
&\alpha^2\iint\limits_{\mathcal{M} \times \mathcal{N}^{*}} s^{-1}\phi^{-1} |{\bf D}_{x}^{2}z|^{2}
+\iint\limits_{\mathcal{M} \times \mathcal{N}} {\bf t}^{-}(s^{-1})\phi^{-1} |{\bf D}_{t}z|^{2} +\iint\limits_{\partial\mathcal{M} \times \mathcal{N}} {\bf t}^{-}(s^{-1}) |{\bf D}_{t}z|^{2}+\nonumber \\
& (\alpha^2+\beta^2)\times\left(\ \iint\limits_{\mathcal{M} \times \mathcal{N}^{*}} s^{3}\lambda^{4}\phi^{3} |z|^{2}
+\iint\limits_{\mathcal{M}^{*} \times \mathcal{N}^{*}} s\lambda^{2}\phi |{\bf D}_{x}z|^{2}+\iint\limits_{\mathcal{M} \times \mathcal{N}^{*}}s\lambda^{2} \phi |{\bf A}_{x}{\bf D}_{x}z|^{2}\right.+\nonumber\\
&\left.\iint\limits_{\partial\mathcal{M} \times \mathcal{N}^{*}} s^{3}\lambda^{3}\phi^{3} |z|^{2}
+ \iint\limits_{\partial\mathcal{M} \times \mathcal{N}^{*}} s\lambda\phi {\bf tr}(|{\bf D}_{x}z|^{2})\right) \nonumber \\
\le\ &C\left(\ \iint\limits_{\mathcal{M} \times \mathcal{N}} {\bf t}^{-}(r^{2})|\mathcal{P}(q)|^{2}
+\iint\limits_{\{0\} \times \mathcal{N}} {\bf t}^{-}(r^{2})|\mathcal B_{\Gamma_{0}}(q)|^{2}
+\iint\limits_{\{1\} \times \mathcal{N}} {\bf t}^{-}(r^{2})|\mathcal B_{\Gamma_{1}}(q)|^{2}\right) +\nonumber \\
&C({\lambda}) \left(\ \  \iint\limits_{\mathcal{\omega}\times \mathcal{N}^{*}} s^{3}|z|^{2}
+\iint\limits_{\mathcal{M} \times \partial\mathcal{N}} {\bf t}^{+}(s^{2}){\bf t}^{+}(|z|^{2})+\iint\limits_{\mathcal{M}^{*} \times \partial\mathcal{N}}{\bf t}^{+} (s^{2}) {\bf t}^{+}(|{\bf A}_{x}z|^{2})
+\right.\nonumber \\
&\left.\ \iint\limits_{\mathcal{M}^{*} \times \partial\mathcal{N}} {\bf t}^{+}(|{\bf D}_{x}z|^{2})\right).
\end{align}

Moreover, since $\text{max}_{t \in [0,T+\triangle t]} \theta(t) \le {2}({\delta T^{2}})^{-1}$ and ${\tau \Delta x}({\delta T^{2}})^{-1} \le \varepsilon_0$, we obtain
\begin{align}\label{3.27}\iint\limits_{\mathcal{M} \times \partial\mathcal{N}}{\bf t}^{+} (s^2){\bf t}^{+}(|z|^{2})
\le 4 (\Delta x)^{-2} \iint\limits_{\mathcal{M} \times \partial\mathcal{N}} \left(\frac{\tau \Delta x}{\delta T^{2}}\right)^{2}{\bf t}^{+}(|z|^{2})
\le C (\Delta x)^{-2} \iint\limits_{\mathcal{M} \times \partial\mathcal{N}} {\bf t}^{+}(|z|^{2}).
\end{align}
Noticing that ${\bf s}_{+}(\mathcal M^*)\cup {\bf s}_{-}(\mathcal M^*)=\mathcal M \cup \partial \mathcal M$ and $$|{\bf A}_{x}z|^{2} \leq \frac{1}{4}\left({\bf s}_{+}(|z^{2}|) + {\bf s}_{-}(|z|)^2\right) \ {\rm and} \ |{\bf D}_{x}z|^{2} \le \frac{1}{2} (\Delta x)^{-2}\left({\bf s}_{+}(|z|^{2})+{\bf s}_{-}(|z|^{2})\right),$$
we deduce that 
\begin{align}
\iint\limits_{\mathcal{M} \times \partial\mathcal{N}} {\bf t}^{+}(s^{2}){\bf t}^{+}(|{\bf A}_{x} z|^{2})
\le& \ C \left(\ \iint\limits_{\mathcal{M} \times \partial\mathcal{N}} {\bf t}^{+} (s^2){\bf t}^{+}(|z|^{2})
+\iint\limits_{\partial\mathcal{M} \times \partial\mathcal{N}} {\bf t}^{+} (s^2){\bf t}^{+}(|z|^{2})\right) \nonumber \\
\le& \ C  (\Delta x)^{-2} \left(\ \iint\limits_{\mathcal{M} \times \partial\mathcal{N}} {\bf t}^{+}(|z|^{2})
+\iint\limits_{\partial\mathcal{M} \times \partial\mathcal{N}}{\bf t}^{+}( |z|^{2})\right)
\end{align}
and
\begin{align}\label{3.29}
\iint\limits_{\mathcal{M}^{*} \times \partial\mathcal{N}} {\bf t}^{+}(|{\bf D}_{x}z|^{2})
\le C (\Delta x)^{-2} \left(\ \iint\limits_{\mathcal{M} \times \partial\mathcal{N}} {\bf t}^{+}(|z|^{2})
+\iint\limits_{\partial\mathcal{M} \times \partial\mathcal{N}} {\bf t}^{+}(|z|^{2})\right).
\end{align}

Substituting (\ref{3.27})-(\ref{3.29}) in (\ref{3.26}), we obtain
\begin{align}\label{3.30}
&\alpha^2\iint\limits_{\mathcal{M} \times \mathcal{N}^{*}} s^{-1} |{\bf D}_{x}^{2}z|^{2}
+\iint\limits_{\mathcal{M} \times \mathcal{N}} {\bf t}^{-}(s^{-1}) |{\bf D}_{t}z|^{2} +\iint\limits_{\partial\mathcal{M} \times \mathcal{N}} {\bf t}^{-}(s^{-1}) |{\bf D}_{t}z|^{2}+\nonumber \\
& (\alpha^2+\beta^2)\times\left(\ \iint\limits_{\mathcal{M} \times \mathcal{N}^{*}} s^{3}\lambda^{4}\phi^{3} |z|^{2}
+\iint\limits_{\mathcal{M}^{*} \times \mathcal{N}^{*}} s\lambda^{2}\phi |{\bf D}_{x}z|^{2}+\iint\limits_{\mathcal{M} \times \mathcal{N}^{*}}s\lambda^{2} \phi |{\bf A}_{x}{\bf D}_{x}z|^{2}\right.+\nonumber\\
&\left.\iint\limits_{\partial\mathcal{M} \times \mathcal{N}^{*}} s^{3}\lambda^{3}\phi^{3} |z|^{2}
+ \iint\limits_{\partial\mathcal{M} \times \mathcal{N}^{*}} s\lambda\phi {\bf tr}(|{\bf D}_{x}z|^{2})\right) \nonumber \\
\le\ &C\left(\ \iint\limits_{\mathcal{M} \times \mathcal{N}} {\bf t}^{-}(r^{2})|\mathcal{P}(q)|^{2}
+\iint\limits_{\{0\} \times \mathcal{N}} {\bf t}^{-}(r^{2})|\mathcal B_{\Gamma_{0}}(q)|^{2}
+\iint\limits_{\{1\} \times \mathcal{N}} {\bf t}^{-}(r^{2})|\mathcal B_{\Gamma_{1}}(q)|^{2}\right) +\nonumber \\
&C({\lambda}) \iint\limits_{\mathcal{\omega} \times \mathcal{N}^{*}} s^{3}|z|^{2}
+C(\Delta x)^{-2}\left(\ \iint\limits_{\mathcal{M} \times \partial\mathcal{N}} {\bf t}^{+}(|z|^{2})+\iint\limits_{\partial\mathcal{M} \times \partial\mathcal{N}} {\bf t}^{+}(|z|^{2})\right).
\end{align}
Finally, by standard argument we can return to the original variable and then obtain (\ref{CE}). This completes the proof of  Theorem \ref{carleman estimate}. \hfill$\Box$


\section{Proof of Theorem 1.2}

In this section, we will apply Carleman estimate (\ref{CE}) to prove observability inequality \eqref{gcbds} for the adjoint system \eqref{dofc}, i.e. Theorem 1.2. To do this, we also need an energy estimate for the full-discrete Ginzburg-Landau equation with dynamic boundary conditions.

\begin{lem} \label{ditui} For all $t^n\in {\mathcal{N}}$, there exists positive constant $C$ such that
\begin{align}\label{4.1}
\iint\limits_{\overline {\mathcal M} \times \{0\}} {\bf t}^{+}(|q|^{2})
\le C\iint\limits_{\overline {\mathcal M} \times \{t^n\}} {\bf t}^{+}(|q|^{2})
\end{align}
for all $9<\Delta t\leq 1/(4|c|)$, where $q\in \mathbb{C}(\overline {\mathcal M}\times\overline{\mathcal N ^*})$ is the solution to (\ref{dofc}).
\end{lem}

\noindent{\bf Proof.} \ Multiplying by ${\bf t}^-(q^*)$   the equation of $q$ in \eqref{dofc}, integrating over $\mathcal{M} $ and extracting the real part, we obtain for all $t^n\in \mathcal N$ that
\begin{align*}
-\text{Re} \int\limits_{\mathcal{M} } {\bf t}^{-}(q^{*}) {\bf D}_{t}q
-\text{Re} (\alpha-i\beta) \int\limits_{\mathcal{M}} {\bf t}^{-}(q^{*}) {\bf D}_{x}^{2}{\bf t}^{-}(q)
+\int\limits_{\mathcal{M}} c |{\bf t}^{-}(q)|^{2}
=0.
\end{align*}
Applying 
\begin{align*}
{\rm Re}({\bf t}^-(q^*){\bf D}_tq)=\frac{1}{2}{\bf D}_t(|q|^2)-\frac{1}{2}\Delta t|{\bf D}_t q|^2
\end{align*}
and discrete integration by parts \eqref{fbjf1}, we further obtain
\begin{align} \label{fcnlgj}
&-\frac{1}{2} \int\limits_{\mathcal{M} }{\bf D}_t(|q|^{2})
+\frac{1}{2}\Delta t \ \int\limits_{\mathcal{M} } |{\bf D}_{t}q|^{2}
+\int\limits_{\mathcal{M} } c {\bf t}^{-}(|q|^{2})
+\alpha \int\limits_{\mathcal{M}^{*} } {\bf t}^{-}(|{\bf D}_{x}  q|^{2}) \nonumber \\
=&\ \text{Re} (\alpha-i\beta) \int\limits_{\partial\mathcal{M} } {\bf t}^{-}(q^{*}) {\bf t}^{-}({\bf tr}({\bf D}_{x}q)) n_{x}.
\end{align}
Similarly, by the dynamic boundary conditions in  \eqref{dofc}, we can derive
\begin{align} \label{tjnlgj}
&\text{Re} (\alpha-i\beta) \int\limits_{\partial \mathcal M} {\bf t}^{-}(q^{*}){\bf t}^{-}( {\bf tr}({\bf D}_{x}q)) n_x\nonumber\\
=\ &\frac{1}{2}\  \int\limits_{\partial \mathcal M } {\bf D}_t(|q|^{2})
-\frac{1}{2}\Delta t \ \int\limits_{\partial \mathcal M} |{\bf D}_{t}q|^{2} -\int\limits_{\partial \mathcal M} c {\bf t}^{-}(|q|^{2}).
\end{align}
Next, substituting  \eqref{tjnlgj} into \eqref{fcnlgj} yields
\begin{align*}
&-\frac{1}{2} \int\limits_{\overline{\mathcal{M}} }{\bf D}_t(|q|^{2})
+\frac{1}{2}\Delta t \ \int\limits_{ \overline{\mathcal{M}} } |{\bf D}_{t}q|^{2}
+\alpha \int\limits_{\mathcal{M}^{*} } {\bf t}^{-}(|{\bf D}_{x}  q|^{2}) \leq |c| \int\limits_{\overline{\mathcal{M}}  } {\bf t}^{-}(|q|^{2}),
\end{align*}
which implies
\begin{align}
&\iint\limits_{\overline {\mathcal M} \times \{t^n\}} {\bf t}^{-}(|q|^{2})\leq \iint\limits_{\overline {\mathcal M} \times \{t^{n+1}\}} {\bf t}^{-}(|q|^{2})+2|c|\Delta t\iint\limits_{\overline {\mathcal M} \times \{t^n\}} {\bf t}^{-}(|q|^{2}) ,\quad t^n\in \mathcal N.
\end{align}
Therefore, we obtain
\begin{align}\label{4.5}
\iint\limits_{\overline {\mathcal M} \times \{t^n\}} {\bf t}^{-}(|q|^{2})\leq \frac{1}{1-2|c|\Delta t}\iint\limits_{\overline {\mathcal M} \times \{t^{n+1}\}} {\bf t}^{-}(|q|^{2})\leq e^{4|c|\Delta t}\iint\limits_{\overline {\mathcal M} \times \{t^{n+1}\}} {\bf t}^{-}(|q|^{2})
\end{align}
if we choose $\Delta t$ sufficiently small such that $2|c|\Delta t\leq 1/2$. From (\ref{4.5}), we deduce
\begin{align}
\iint\limits_{\overline {\mathcal M} \times \{0\}} {\bf t}^{+}(|q|^{2})
\le e^{4|c|n\Delta t}\iint\limits_{\overline {\mathcal M} \times \{t^n\}} {\bf t}^{+}(|q|^{2})\leq C\iint\limits_{\overline {\mathcal M} \times \{t^n\}} {\bf t}^{+}(|q|^{2}).
\end{align}
 This completes the proof of Lemma 4.1.\hfill$\Box$

\vspace{2mm}

Now we prove Theorem 1.2.

\vspace{2mm}

\noindent{\bf Proof of Theorem 1.2.}\ Firstly, applying Carleman estimate \eqref{CE} to $q$, we obtain
\begin{align} \label{CEyy}
&\iint\limits_{\overline{\mathcal{M}} \times \mathcal{N}^{*}} s^{3}r^{2} |q|^{2}\nonumber \\
\le\ &C\rho^2\ \iint\limits_{\overline{\mathcal{M}}\times \mathcal{N}} {\bf t}^{-}(r^2){\bf t}^{-}(|q|^{2})+C(\lambda) \iint\limits_{\omega \times \mathcal{N}^{*}} s^{3}r^{2} |q|^{2}+C (\Delta x)^{-2}\ \iint\limits_{\overline{\mathcal{M}}\times \partial\mathcal{N}} {\bf t}^{+}(r^2){\bf t}^{+}(|q|^{2}),
\end{align}
where $\rho := \text{max}\{|c|,|\gamma|\}$.
Then, we choose $\tau$ sufficiently large such that $s^3\geq C\rho^2$, i.e. $\tau \ge CT^{2} \rho^{{2}/{3}}$ 
to absorb the first term on the right-hand side. Therefore, for all $\tau \ge \tau_{1}(T+T^{2}+T^{2}\rho^{\frac{2}{3}})$ with $\tau_{1} \ge \tau_{0}$, we have
\begin{align} \label{4.8}
&\iint\limits_{\overline{\mathcal{M}} \times \mathcal{N}^{*}} s^{3}r^{2} |q|^{2}\le\ C\iint\limits_{\omega \times \mathcal{N}^{*}} s^{3}r^{2} |q|^{2}+C (\Delta x)^{-2}\ \iint\limits_{\overline{\mathcal{M}}\times \partial\mathcal{N}} {\bf t}^{+}(r^2){\bf t}^{+}(|q|^{2}).
\end{align}

Let $\widehat {\mathcal N}=(T/4,3T/4) \cap \mathcal N$. Using ${\bf t}^+(r^2) \ge e^{-\frac{2^{5}\tau K_{0}}{3T^{2}}}$ in $\overline {\mathcal M}\times \widehat {\mathcal N}$, where $K_{0}:=\text{max}_{\overline{\Omega}}\{-\varphi\}$, we obtain that
\begin{align}\label{4.9}
\iint\limits_{\overline{\mathcal{M}} \times \mathcal{N}^{*}} s^{3}r^{2} |q|^{2}
\ge \frac{\tau^{3}}{T^{6}} e^{-\frac{2^{5}\tau K_{0}}{3T^{2}}} \sum\limits_{t^n\in \widehat {\mathcal N}}\ \ \iint\limits_{\overline {\mathcal{M}} \times \{t^n\}}{\bf t}^+( |q|^{2})\Delta t.
\end{align}
On the other hand, we apply (\ref{4.1}) to obtain
\begin{align}\label{4.10}
\iint\limits_{\overline {\mathcal M} \times \{0\}} {\bf t}^{+}(|q|^{2})\leq C\sum\limits_{t^n\in \widehat {\mathcal N}}\ \ \iint\limits_{\overline {\mathcal{M}} \times \{t^n\}}{\bf t}^+( |q|^{2})\Delta t.
\end{align}
From (\ref{4.9}) and (\ref{4.10}), it follows 
\begin{align}\label{4.11}
\iint\limits_{\overline{\mathcal{M}} \times \mathcal{N}^{*}} s^{3}r^{2} |q|^{2}
\ge C\frac{\tau^{3}}{T^{6}} e^{-\frac{2^{5}\tau K_{0}}{3T^{2}}}\iint\limits_{\overline {\mathcal M} \times \{0\}} {\bf t}^{+}(|q|^{2})\geq CT e^{-\frac{C \tau}{T^{2}}} \iint\limits_{\overline {\mathcal M} \times \{0\}} {\bf t}^{+}(|q|^{2}).
\end{align}
Using ${\bf t}^+(r^2)\leq e^{-\frac{2\tau k_0}{\delta T^2}}$ with $k_{0}:=\text{min}_{\overline{\Omega}}\{-\varphi\}$ in $\overline{\mathcal M}\times \partial \mathcal N$ and (\ref{4.1}), we have
\begin{align}\label{4.12}
&(\Delta x)^{-2}\ \iint\limits_{\overline{\mathcal{M}}\times \partial\mathcal{N}} {\bf t}^{+}(r^2){\bf t}^{+}(|q|^{2})
\le C(\Delta x)^{-2}\ \iint\limits_{\overline{\mathcal{M}}\times \{T\}} {\bf t}^{+}(r^2){\bf t}^{+}(|q|^{2})\nonumber\\
\leq \ & C(\Delta x)^{-2} e^{-\frac{2\tau k_{0}}{\delta T^{2}}}  \iint\limits_{\overline{\mathcal{M}}\times \{T\}}{\bf t}^{+}(|q|^{2})\le C (\Delta x)^{-2} e^{-\frac{2\tau k_{0}}{\delta T^{2}}} \|q_{T}\|_{{L}^{2}(\overline{\mathcal{M}})}^{2}.
\end{align}

Therefore, from (\ref{4.8}), (\ref{4.11}) and (\ref{4.12}), it follows that
\begin{align}
\|{\bf t}^{+}q(0)\|_{{L}^{2}(\overline{\mathcal{M}})}^{2}
\le C T^{-1} e^{\frac{C \tau}{T^{2}} } \iint\limits_{\omega \times \mathcal{N}^{*}} |q|^{2}
+C T^{-1} (\Delta x)^{-2} e^{\frac{\tau}{T^{2}}(C - \frac{C'}{\delta}) } ||q_{T}||_{{L}^{2}(\overline{\mathcal{M}})}^{2}
\end{align}
for $\tau \ge \tau_{2}(T+T^{2}+T^{2}\rho^{\frac{2}{3}})$ with $\tau_{2}=\text{max}\{\tau_{1},\frac{3T^{2}}{2k_{0}}\}$.
Then, for $0<\delta<\delta_{1}<\frac{1}{2}$ small enough we obtain by the above inequality that
\begin{align}
\|{\bf t}^{+}q(0)\|_{{L}^{2}(\overline{\mathcal{M}})}^{2}
\le C T^{-1} e^{\frac{C \tau}{T^{2}} } \iint\limits_{\omega \times \mathcal{N}^{*}} |q|^{2}
+C T^{-1} (\Delta x)^{-2} e^{-\frac{C'\tau}{\delta T^{2}} } ||q_{T}||_{{L}^{2}(\overline{\mathcal{M}})}^{2}.
\end{align}
Finally following Step 2 in Proposition 3.1 in [\ref{Carreno-ACM-2023}], we can deduce (\ref{gcbds}) and then complete the proof of Theorem 1.2.\hfill$\Box$

\section{Appendix A}
\renewcommand{\theequation}{A.\arabic{equation}}

\setcounter{equation}{0}
 {In this appendix, we will present some lemmas to show the  estimates for $I_{ij}$  one by one in $2{\rm Re} ( \mathcal A (z), \mathcal B(z) )_{L^{2}(\mathcal{M} \times \mathcal{N})}$.
 
 \vspace{2mm}
 
 {\noindent}{\bf Lemma A.1.}\ {\em Provided ${\tau \Delta x}{(\delta T^{2})}^{-1} \le 1$ and $\tau\geq T^2$, we have
	\begin{align}\label{A.1}
I_{11}
\ge\ &4\alpha^{2} \iint\limits_{\mathcal{M} \times \mathcal{N}^{*}} s^{3}\lambda^{4}\phi^{3}(\partial_{x}\psi)^{4} |z|^{2}
+4\alpha^{2} \iint\limits_{\mathcal{M}^{*} \times \mathcal{N}^{*}} s\lambda^{2}\phi(\partial_{x}\psi)^{2} |{\bf D}_{x}z|^{2} -\nonumber\\
&\alpha^2\iint\limits_{\mathcal{M} \times \mathcal{N}^{*}} (s^{3}\lambda^{3}\phi^{3}{\mathcal O}(1) + s^{2}{\mathcal O}_{\lambda}(1) + s^{3}(s\Delta x)^{2}{\mathcal O}_{\lambda}(1) ) |z|^{2}-\nonumber \\
&\alpha^2\iint\limits_{\mathcal{M}^{*} \times \mathcal{N}^{*}} (s\lambda\phi {\mathcal O}(1) + s(s\Delta x)^{2}{\mathcal O}_{\lambda}(1)) |{\bf D}_{x}z|^{2}
+BT_{11},
\end{align}
where
\begin{align*}
BT_{11}
\ge\ & \alpha^2 c_1 \iint\limits_{\partial\mathcal{M} \times \mathcal{N}^{*}} s^{3}\lambda^{3}\phi^3|z|^{2}+\alpha^2 c_1\iint\limits_{{\partial \mathcal M}\times{\mathcal N}^*} s\lambda\phi{\bf tr}(|{\bf D}_x z|^2)-\nonumber \\
&\alpha^2\iint\limits_{\partial\mathcal{M} \times \mathcal{N}^{*}} (s^{2}+s^{3}(s\Delta x)){\mathcal O}_{\lambda}(1) |z|^{2}-\alpha^2 \iint\limits_{\partial\mathcal{M} \times \mathcal{N}^{*}} s(s\Delta x)^2{\mathcal O}_{\lambda}(1) {\bf tr}(|{\bf D}_{x}z|^{2})
\end{align*}
with a constant $c_1$ such that $c_1\geq \min\{|\psi_x(0)|, |\psi_x(1)|\}>0$.
}

\vspace{2mm}

\noindent{\bf Proof.}\ By the definition of $I_{11}$, we see that
\begin{align}\label{A.2}
I_{11}
=\ &2\alpha^{2} \text{Re} \iint\limits_{\mathcal{M} \times \mathcal{N}^{*}} (r{\bf D}_{x}^{2}\rho {\bf A}_{x}^{2}z + r{\bf A}_{x}^{2}\rho {\bf D}_{x}^{2}z) (2r{\bf A}_{x}{\bf D}_{x}\rho {\bf A}_{x}{\bf D}_{x}z^{*} - s\partial_{x}^2\phi z^{*}) \nonumber\\
=\ &4\alpha^{2} {\rm Re} \iint\limits_{\mathcal{M} \times \mathcal{N}^{*}} r^{2}{\bf A}_x{\bf D}_x\rho {\bf D}_x^{2}\rho {\bf A}_x^{2}z {\bf A}_x{\bf D}_xz^{*} -2\alpha^{2} \text{Re} \iint\limits_{\mathcal{M} \times \mathcal{N}^{*}} s\partial^2_{x}\phi r{\bf D}_x^{2}\rho {\bf A}_x^{2}z z^{*}+\nonumber\\
&4\alpha^{2} \text{Re} \iint\limits_{\mathcal{M} \times \mathcal{N}^{*}} r^{2}{\bf A}_x{\bf D}_x\rho {\bf A}_x^{2}\rho {\bf D}_x^{2}z {\bf A}_x{\bf D}_xz^{*}-2\alpha^{2} \text{Re} \iint\limits_{\mathcal{M} \times \mathcal{N}^{*}} s\partial_{x}^2\phi r{\bf A}_x^{2}\rho {\bf D}_x^{2}zz^{*} \nonumber\\
=\ &X_1+X_2+X_3+X_4.
\end{align}

Obviously, by (\ref{ds1}) we have
\begin{align}\label{A.3}
{\rm Re} ({\bf A}_x{\bf D}_x z^* {\bf D}_x^2 z)=\frac{1}{2}({\bf A}_x{\bf D}_x z^* {\bf D}_x^2 z+{\bf A}_x{\bf D}_x z {\bf D}_x^2 z^*)=\frac{1}{2}{\bf D}_x(|{\bf D}_x z|^2).
\end{align}
Using (\ref{ds2}), (\ref{A.3}) and integration by parts in space, i.e Lemma 2.2, we can obtain that
\begin{align*}
& X_1\nonumber\\ 
=\ &4\alpha^{2} {\rm Re} \iint\limits_{\mathcal{M} \times \mathcal{N}^{*}} r^{2}{\bf A}_x{\bf D}_x\rho {\bf D}_x^{2}\rho z {\bf A}_x{\bf D}_xz^{*} +\frac{1}{2}\alpha^{2}(\Delta x)^2\iint\limits_{\mathcal{M} \times \mathcal{N}^{*}} r^{2}{\bf A}_x{\bf D}_x\rho {\bf D}_x^{2}\rho {\bf D}_x(|{\bf D}_x z|^2) \nonumber\\
=\ &4\alpha^{2} {\rm Re} \iint\limits_{\mathcal{M}^* \times \mathcal{N}^{*}} {\bf A}_x(r^{2}{\bf A}_x{\bf D}_x\rho {\bf D}_x^{2}\rho z) {\bf D}_xz^{*} -\frac{1}{2}\alpha^{2}(\Delta x)^2\iint\limits_{\mathcal{M}^* \times \mathcal{N}^{*}} {\bf D}_x(r^{2}{\bf A}_x{\bf D}_x\rho {\bf D}_x^{2}\rho) |{\bf D}_x z|^2+\nonumber\\
&B_1+B_2\nonumber\\
=\ &-2\alpha^{2} \iint\limits_{\mathcal{M} \times \mathcal{N}^{*}} {\bf A}_{x}{\bf D}_{x}(r^{2}{\bf A}_{x}{\bf D}_{x}\rho {\bf D}_{x}^{2}\rho) |z|^{2}+\frac{1}{2}\alpha^{2} (\Delta x)^{2} \iint\limits_{\mathcal{M}^{*} \times \mathcal{N}^{*}} {\bf D}_{x}(r^{2}{\bf A}_{x}{\bf D}_{x}\rho {\bf D}_{x}^{2}\rho) |{\bf D}_{x}z|^{2}+\nonumber\\
&B_{1}+B_2+B_3,
\end{align*}
where
\begin{align*}
&B_{1}=-2\alpha^{2} \Delta x {\rm Re} \iint\limits_{\partial\mathcal{M} \times \mathcal{N}^{*}} r^{2}{\bf A}_{x}{\bf D}_{x}\rho {\bf D}_{x}^{2}\rho z {\bf tr}({\bf D}_{x}z^{*}),\\
&B_2=\frac{1}{2}\alpha^{2} (\Delta x)^{2} \iint\limits_{\partial\mathcal{M} \times \mathcal{N}^{*}} r^{2}{\bf A}_{x}{\bf D}_{x}\rho {\bf D}_{x}^{2}\rho {\bf tr}(|{\bf D}_{x}z|^{2}) n_{x},\\
&B_3=2\alpha^{2} \iint\limits_{\partial\mathcal{M} \times \mathcal{N}^{*}}  {\bf tr}({\bf A}_{x}(r^{2}{\bf A}_{x}{\bf D}_{x}\rho {\bf D}_{x}^{2}\rho))|z|^{2} n_{x}.
\end{align*}
From Lemma \ref{gj3}, we can deduce for $\tau \Delta x(\delta  T^2)^{-1}\leq 1$ that 
\begin{align*}
 &\hspace{-3.9cm}{\bf D}_{x}(r^{2}{\bf A}_{x}{\bf D}_{x}\rho {\bf D}_{x}^{2}\rho)
=s^{3}\mathcal O_{\lambda}(1),\\
 {\bf A}_{x}{\bf D}_{x}(r^{2}{\bf A}_{x}{\bf D}_{x}\rho {\bf D}_{x}^{2}\rho)=\ &\partial_x(r^2\partial_x\rho\partial_x^2\rho)+s^3\mathcal O_{\lambda}((s\Delta x)^2)\\
=\ &-3s^{3}\lambda^{4}\phi^{3}(\partial_{x}\psi)^{4} + s^{3}\lambda^{3}\phi^{3}\mathcal O(1) + s^{2}\mathcal O_{\lambda}(1) + s^{3}(s\Delta x)^{2}\mathcal O_{\lambda}(1).
\end{align*}
Therefore,
\begin{align}\label{A.4}
&X_1\ge\  6\alpha^{2} \iint\limits_{\mathcal{M} \times \mathcal{N}^{*}} s^{3}\lambda^{4}\phi^{3}(\partial_{x}\psi)^{4} |z|^{2}
 +Y_1+\sum_{j=1}^3B_j,
\end{align}
where
\begin{align*}
Y_1=\ &-\alpha^{2}\iint\limits_{\mathcal{M} \times \mathcal{N}^{*}} \left(s^{3}\lambda^{3}\phi^{3}{\mathcal O}(1)+s^{2}{\mathcal O}_{\lambda}(1)+s^{3}(s\Delta x)^{2}{\mathcal O}_{\lambda}(1)\right) |z|^{2}-\\
&\ \alpha^2\iint\limits_{\mathcal{M}^{*} \times \mathcal{N}^{*}} s(s\Delta x)^{2}{\mathcal O}_{\lambda}(1) |{\bf D}_{x}z|^{2}.
\end{align*}
From Lemma \ref{gj3}, we obtain that
\begin{align*}{\bf A}_{x}(r^{2}{\bf A}_{x}{\bf D}_{x}\rho {\bf D}_{x}^{2}\rho)
=\ &r^2\partial_x\rho\partial_x^2\rho+s^3\mathcal O_{\lambda}((s\Delta x)^2)\\
=\ &-s^{3}\lambda^{3}\phi^{3}(\partial_{x}\psi)^{3}
+s^{2}O_{\lambda}(1)+s^{3}(s\Delta x)^{2}\mathcal O_{\lambda}(1).
\end{align*}
Since $\psi_{x}(0)>0$ and $\psi_{x}(1)<0$, there exists $c_1\geq \min\{|\psi_x(0)|, |\psi_x(1)|\}>0$ such that 
\begin{align*}
-s^{3}\lambda^{3}\phi^{3}(\partial_{x}\psi)^{3} n_x\geq c_1s^3\lambda^3\phi^3,\quad (x,t)\in \partial\mathcal M\times {\mathcal N}^*.
\end{align*}
Furthermore, we have the following estimate for boundary terms $\sum_{j=1}^3B_j$ 
\begin{align*}
&\sum_{j=1}^3B_j\\
\ge\ & \alpha^2 c_1 \iint\limits_{\partial\mathcal{M} \times \mathcal{N}^{*}} s^{3}\lambda^{3}\phi^3|z|^{2}
-\alpha^2\iint\limits_{\partial\mathcal{M} \times \mathcal{N}^{*}} \left(s(s\Delta x)^2+(s\Delta x)^2+s(s\Delta x)^{4}\right)\mathcal O_{\lambda}(1) {\bf tr}(|{\bf D}_{x}z|^{2})- \nonumber \\
&\alpha^2\iint\limits_{\partial\mathcal{M} \times \mathcal{N}^{*}} \left(s^{2}+s^{2}(s\Delta x)+s(s\Delta x)^{2}+s^{3}(s\Delta x)^{2}\right)\mathcal O_{\lambda}(1) |z|^{2}.
\end{align*}

For $X_2$, integration by parts together with (\ref{ds2}) and
${\rm Re} ({\bf D}_x z{\bf A}_x z^*)=\frac{1}{2}{\bf D}_x(|z|^2)$
yields
\begin{align*}
&X_2 \\
=\ &-2\alpha^{2} {\rm Re}\iint\limits_{\mathcal{M} \times \mathcal{N}^{*}} s\partial^2_{x}\phi r{\bf D}_{x}^{2}\rho|z|^{2}
-\frac{1}{2}{\alpha^{2}}(\Delta x)^{2} {\rm Re}\iint\limits_{\mathcal{M} \times \mathcal{N}^{*}} s\partial_x^2\phi r {\bf D}_{x}^{2}\rho {\bf D}_{x}^{2}z z^*\nonumber\\
=\ &-2\alpha^{2} {\rm Re}\iint\limits_{\mathcal{M} \times \mathcal{N}^{*}} s\partial^2_{x}\phi r{\bf D}_{x}^{2}\rho|z|^{2}
+\frac{1}{4}{\alpha^{2}}(\Delta x)^{2} \iint\limits_{{\mathcal M}^* \times \mathcal{N}^{*}} {\bf D}_x( s\partial_x^2\phi r {\bf D}_{x}^{2}\rho) {\bf D}_{x}(|z|^2)+\nonumber\\  
&\frac{1}{2}\alpha^2 (\Delta x)^{2} \iint\limits_{{\mathcal M}^* \times \mathcal{N}^{*}} {\bf A}_x( s\partial_x^2\phi r {\bf D}_{x}^{2}\rho) |{\bf D}_{x}z|^2+B_4\nonumber\\
=\ &-2\alpha^{2} {\rm Re}\iint\limits_{\mathcal{M} \times \mathcal{N}^{*}} s\partial^2_{x}\phi r{\bf D}_{x}^{2}\rho|z|^{2}
-\frac{1}{4}\alpha^{2}(\Delta x)^{2} \iint\limits_{\mathcal{M} \times \mathcal{N}^{*}} {\bf D}_{x}^{2}(s\partial_{x}^2\phi r{\bf D}_{x}^{2}\rho) |z|^{2} +\nonumber \\
&\frac{1}{2}\alpha^{2}(\Delta x)^{2} \iint\limits_{\mathcal{M}^{*} \times \mathcal{N}^{*}} {\bf A}_{x}(s\partial_{x}^2\phi r{\bf D}_{x}^{2}\rho) |{\bf D}_{x}z|^{2}
+B_4+B_5,
\end{align*}
where
\begin{align*}
&B_4=-\frac{1}{2}\alpha^{2} (\Delta x)^{2} {\rm Re} \iint\limits_{\partial\mathcal{M} \times \mathcal{N}^{*}} s\partial^2_{x}\phi r{\bf D}_{x}^{2}\rho z^{*} {\bf tr}({\bf D}_{x}z) n_x, \\
&B_5=\frac{1}{4} \alpha^{2} (\Delta x)^{2} \iint\limits_{\partial\mathcal{M}\times \mathcal{N}^{*}} {\bf tr}({\bf D}_{x}(s\partial^2_{x}\phi r{\bf D}_{x}^{2}\rho)) |z|^{2}  n_x.
\end{align*}
From Lemma \ref{gj2}, we can see that
\begin{align*}
&{\bf A}_{x}(s\partial^2_{x}\phi r{\bf D}_{x}^{2}\rho)
=s^{3}\mathcal O_{\lambda}(1),\quad {\bf D}_{x}^{k}(s\partial^2_{x}\phi r{\bf D}_{x}^{2}\rho)=s^{3}\mathcal O_{\lambda}(1), \ k=0,1,2\\
&s\partial^2_{x}\phi r{\bf D}_{x}^{2}\rho
=s^{3}\lambda^{4}\phi^{3}|\partial_{x}\psi|^{4}
+s^{3}\lambda^{3}\phi^{3}|\partial_x\psi|^2\partial_x^2\psi +s^{2}\mathcal O_{\lambda}(1)+s^{3}(s\Delta x)^{2}\mathcal O_{\lambda}(1).
\end{align*}
Then, we obtain the following estimate for $X_2$
\begin{align}\label{A.5}
X_2\ge\ -2\alpha^{2} \iint\limits_{\mathcal{M} \times \mathcal{N}^{*}} s^{3}\lambda^{4}\phi^{3}(\partial_{x}\psi)^{4} |z|^{2} +Y_2+\sum_{j=4}^5{B_j},
\end{align}
where
\begin{align*}
Y_2=\ &-\alpha^2\iint\limits_{\mathcal{M} \times \mathcal{N}^{*}} (s^{3}\lambda^{3}\phi^{3}\mathcal O(1)+s^{2}\mathcal O_{\lambda}(1)+s^{3}(s\Delta x)^{2}\mathcal O_{\lambda}(1)+s(s\Delta x)^{2}\mathcal O_{\lambda}(1)) |z|^{2}-\\
&\ \alpha^2\iint\limits_{\mathcal{M}^{*} \times \mathcal{N}^{*}} s(s\Delta x)^{2}\mathcal O_{\lambda}(1) |{\bf D}_{x}z|^{2}
\end{align*}
and
\begin{align*}
& \sum_{j=4}^5{B_j}
\ge -\alpha^2\iint\limits_{\partial\mathcal{M} \times \mathcal{N}^{*}} s(s\Delta x)^{2}\mathcal O_{\lambda}(1) |z|^{2}
-\alpha^2\iint\limits_{\partial\mathcal{M} \times \mathcal{N}^{*}} s(s\Delta x)^{2}\mathcal O_{\lambda}(1) {\bf tr}(|{\bf D}_{x}z|^{2}).
\end{align*}

Similarly, by using integration by parts in space yields we obtain
\begin{align*}
&X_3+X_4\nonumber\\
=\ &2\alpha^{2} \iint\limits_{\mathcal{M} \times \mathcal{N}^{*}} r^{2}{\bf A}_x{\bf D}_x\rho {\bf A}_x^{2}\rho {\bf D}_x(|{\bf D}_x z|^{2})+2\alpha^{2} \text{Re} \iint\limits_{\mathcal{M}^* \times \mathcal{N}^{*}} {\bf D}_x(s\partial_{x}^2\phi r{\bf A}_x^{2}\rho)  {\bf A}_x z^*{\bf D}_xz+\nonumber\\
&2\alpha^{2} \iint\limits_{\mathcal{M}^* \times \mathcal{N}^{*}} {\bf A}_x(s\partial_{x}^2\phi r{\bf A}_x^{2}\rho)  |{\bf D}_x z|^2+B_6 \nonumber\\
=\ &-2\alpha^{2} \iint\limits_{\mathcal{M}^{*} \times \mathcal{N}^{*}} {\bf D}_{x}(r^{2}{\bf A}_{x}{\bf D}_{x}\rho {\bf A}_{x}^{2}\rho) |{\bf D}_{x}z|^{2}-\alpha^{2} \iint\limits_{\mathcal{M} \times \mathcal{N}^{*}} {\bf D}_{x}^{2}(s\partial_{x}^2\phi r{\bf A}_{x}^{2}\rho) |z|^{2}+\nonumber\\
&2\alpha^{2} \iint\limits_{\mathcal{M}^{*} \times \mathcal{N}^{*}} {\bf A}_{x}(s\partial^2_{x}\phi r{\bf A}_{x}^{2}\rho) |{\bf D}_{x}z|^{2}+B_6+B_7+B_8,
\end{align*} 
where
\begin{align*}
&B_6=-2\alpha^2 {\rm Re}\iint\limits_{\partial\mathcal{M}\times \mathcal{N}^{*}} s\partial^2_{x}\phi r{\bf A}_{x}^{2}\rho{\bf tr}({\bf D}_x z)z^*n_x,\\
&B_7=2\alpha^{2} \iint\limits_{\partial\mathcal{M}\times \mathcal{N}^{*}}  r^2{\bf A}_x{\bf D}_{x}\rho {\bf A}_{x}^{2}\rho {\bf tr}(|{\bf D}_xz|^{2})  n_x,\\
&B_8=\alpha^{2} \iint\limits_{\partial\mathcal{M}\times \mathcal{N}^{*}} {\bf tr}({\bf D}_x( s\partial_x^2\phi r{\bf A}_x^2\rho)) |z|^{2}  n_x.
\end{align*}
From Lemma \ref{gj2} and Lemma \ref{gj3}, we see 
\begin{align*}
& {\bf D}_{x}(r^{2}{\bf A}_{x}{\bf D}_{x}\rho {\bf A}_{x}^{2}\rho)
=-s\lambda^{2}\phi(\partial_{x}\psi)^{2} + s\lambda\phi {\mathcal O}(1) + s(s\Delta x)^{2}{\mathcal O}_{\lambda}(1),\\
&{\bf A}_{x}(s\partial^2_{x}\phi r{\bf A}_{x}^{2}\rho)={\bf A}_{x}(s\partial^2_{x}\phi){\bf A}_x( r{\bf A}_{x}^{2}\rho)+\frac{(\Delta x)^2}{4}{\bf D}_{x}(s\partial^2_{x}\phi){\bf D}_x (r{\bf A}_{x}^{2}\rho)\\
&\hspace{2.3cm}= s\lambda^2\phi(\partial_x\psi)^2+s\lambda\phi\mathcal O(1)+s(s\Delta x)^2\mathcal O_{\lambda}(1),\\
& r^{2}{\bf A}_{x}{\bf D}_{x}\rho {\bf A}_{x}^{2}\rho=-s\lambda\phi\partial_{x}\psi + s(s\Delta x)^{2}{\mathcal O}_{\lambda}(1),\\
& {\bf D}_x^j(s\partial_x^2\phi r {\bf A}_x^2\rho)=s{\mathcal O}_{\lambda}(1),\quad j=0,1,2,
\end{align*}
and then obtain the following estimate for $X_3+X_4$
\begin{align}\label{A.6}
& X_3+X_4\geq 4\alpha^{2} \iint\limits_{\mathcal{M}^{*} \times \mathcal{N}^{*}} s\lambda^{2}\phi(\partial_{x}\psi)^{2} |{\bf D}_{x}z|^{2}+Y_3+\sum_{j=6}^8{B_j},
\end{align}
where
\begin{align*}
Y_3=\ &-\alpha^2\iint\limits_{\mathcal{M} \times \mathcal{N}^{*}} s{\mathcal O}_{\lambda}(1)|z|^{2}-\alpha^2\iint\limits_{\mathcal{M}^{*} \times \mathcal{N}^{*}} (s\lambda\phi \mathcal O(1)+s(s\Delta x)^{2}{\mathcal O}_{\lambda}(1)) |{\bf D}_{x}z|^{2}
\end{align*}
and
\begin{align*}
\sum_{j=6}^8{B_j}
\ge&\alpha^2 c_1\iint\limits_{{\partial \mathcal M}\times{\mathcal N}^*} s\lambda\phi{\bf tr}(|{\bf D}_x z|^2)-\alpha^2\iint\limits_{\partial\mathcal{M} \times \mathcal{N}^{*}} s{\mathcal O}_{\lambda}(1) |z|^{2}.
\end{align*}

Finally, substituting (\ref{A.4})-(\ref{A.6}) into (\ref{A.2}) and noticing that $s\Delta x\leq 1$ and $s=\tau\theta(t)\geq 1$ provided ${\tau \Delta x}{(\delta T^{2})}^{-1} \le 1$ and $\tau \geq T^2$, we obtain (\ref{A.1}) and then complete the proof of Lemma A.1. \hfill$\Box$

\vspace{2mm}
 
 {\noindent}{\bf Lemma A.2.}\ {\em Provided ${\tau \Delta x}{(\delta T^{2})}^{-1} \le 1$, $\tau \Delta t\left(\delta^{2}T^{3}\right)^{-1} \le {1}/{2}$ and $\tau\geq T^2$, we have
 \begin{align}\label{A.7}I_{13}\geq \ &\alpha\iint\limits_{\mathcal{M}^{*} \times \mathcal{N}}\left(\Delta t-{\bf t}^{-}\left(\Delta t(s\Delta x)^{2}+(\Delta t)^{2}\sigma_{2}+(\Delta x)^2(\Delta t)^3\sigma_1\right)\right) \mathcal O_{\lambda}(1)|{\bf D}_{t}{\bf D}_{x}z|^{2}-\nonumber\\
&\alpha\iint\limits_{\mathcal{M} \times \mathcal{N}^{*}} s(s\Delta x)^{2}\mathcal O_{\lambda}(1) |{\bf A}_{x}{\bf D}_{x}z|^{2}-\alpha\iint\limits_{\mathcal{M} \times \mathcal{N}^{*}} s^{-3}(s\Delta x)^{6}\mathcal O_{\lambda}(1) |{\bf D}_{x}^{2}z|^{2}-\nonumber\\
&
\alpha\iint\limits_{\mathcal{M} \times \mathcal{N}} {\bf t}^{-}\left(s^{2}\Delta t+(\Delta t)^{2}\sigma_{1}+s^{-1}(s\Delta x)^{2}\right)\mathcal O_{\lambda}(1) |{\bf D}_{t}z|^{2}-\nonumber \\
&\alpha\iint\limits_{\mathcal{M} \times \mathcal{N}^{*}} \sigma_1 |z|^{2}-\alpha\iint\limits_{\mathcal{M}^{*} \times \mathcal{N}^{*}} \left((\Delta x)^2\sigma_1 +s(s\Delta x)^{2}\mathcal O_{\lambda}(1)+\sigma_{2}\right) |{\bf D}_{x}z|^{2}+BT_{13},
\end{align}
where
\begin{align*}
&BT_{13}\\
\ge\ & -\iint\limits_{\mathcal{M} \times \partial\mathcal{N}} {\bf t}^{+}(s)^{2}\mathcal O_{\lambda}(1){\bf t}^{+}(|z|^{2})-\alpha\iint\limits_{\mathcal{M}^{*} \times \partial\mathcal{N}} {\bf t}^{+}(1+(s\Delta x)^{2}\mathcal O_{\lambda}(1)) {\bf t}^{+}(|{\bf D}_{x}z|^{2})
- \nonumber \\
&\alpha\iint\limits_{\partial\mathcal{M} \times \mathcal{N}} {\bf t}^{-}(\epsilon s^{-1}+s^{-1}(s\Delta x)^{2}\mathcal O_{\lambda}(1)) |{\bf D}_{t}z|^{2}-\alpha\iint\limits_{\partial\mathcal{M} \times \mathcal{N}^{*}} (C(\epsilon)s+s(s\Delta x)^{2}\mathcal O_{\lambda}(1)) {\bf tr}(|{\bf D}_{x}z|^{2})
\end{align*}
for any $\epsilon>0$.
 }

\vspace{2mm}

{\noindent\bf Proof.}\ We write $I_{13}$ as
\begin{align}\label{A.8}
I_{13}
=&\ 2\alpha \text{Re} \iint\limits_{\mathcal{M} \times \mathcal{N}} {\bf t}^{-}(r{\bf D}_{x}^{2}\rho {\bf A}_{x}^{2}z) {\bf D}_{t}z^{*} + 2\alpha \text{Re} \iint\limits_{\mathcal{M} \times \mathcal{N}} {\bf t}^{-} (r{\bf A}_{x}^{2}\rho {\bf D}_{x}^{2}z) {\bf D}_{t}z^{*} = X_{5}+X_{6}.
\end{align}
Using the equality (\ref{ds2}), (\ref{2.6}) and integration by parts in time, we have
\begin{align}\label{1-A.9}
X_5 =\ &2\alpha \text{Re} \iint\limits_{\mathcal{M} \times \mathcal{N}} {\bf t}^{-}(r{\bf D}_{x}^{2}\rho) {\bf t}^{-}\left(z+\frac{1}{4}(\Delta x)^{2}{\bf D}_{x}^{2}z\right) {\bf D}_{t}z^{*} \nonumber\\
=\ &\alpha \iint\limits_{\mathcal{M} \times \mathcal{N}^*} r{\bf D}_{x}^{2}\rho {\bf D}_t{\bf t}^{+}(|z|^2)-\alpha \Delta t \iint\limits_{\mathcal{M} \times \mathcal{N}} {\bf t}^{-}(r{\bf D}_{x}^{2}\rho) |{\bf D}_{t}z|^{2}+\nonumber\\
&\frac{1}{2}\alpha (\Delta x)^{2} \text{Re} \iint\limits_{\mathcal{M} \times \mathcal{N}} {\bf t}^{-}(r{\bf D}_{x}^{2}\rho)  {\bf t}^{-}({\bf D}_{x}^{2}z) {\bf D}_{t}z^{*}\nonumber\\
=\ &-\alpha \iint\limits_{\mathcal{M} \times \mathcal{N}} {\bf D}_{t}(r{\bf D}_{x}^{2}\rho) {\bf t}^{+}(|z|^{2})
-\alpha \Delta t \iint\limits_{\mathcal{M} \times \mathcal{N}} {\bf t}^{-}(r{\bf D}_{x}^{2}\rho) |{\bf D}_{t}z|^{2} +\nonumber \\
&\frac{1}{2}\alpha (\Delta x)^{2} \text{Re} \iint\limits_{\mathcal{M} \times \mathcal{N}} {\bf t}^{-}(r{\bf D}_{x}^{2}\rho)  {\bf t}^{-}({\bf D}_{x}^{2}z) {\bf D}_{t}z^{*}+B_9,
\end{align}
where
\begin{align*}
B_9=\alpha \iint\limits_{\mathcal{M} \times \partial\mathcal{N}} {\bf t}^{+}(r{\bf D}_{x}^{2}\rho){\bf t}^{+}( |z|^{2}) n_{t}.
\end{align*}
For the third term in (\ref{1-A.9}), we use integrations by parts with respect to ${\bf D}_x$ and ${\bf D}_t$ given by (\ref{fbjf1}) and (\ref{2.10}) respectively, and  
${\rm Re}({\bf t}^-({\bf D}_xz){\bf D}_t{\bf D}_x z^*)=1/2{\bf D}_t(|{\bf D}_xz|^2)-\Delta t/2|{\bf D}_t{\bf D}_xz|^2$ to obtain
\begin{align*}
&\frac{1}{2}\alpha (\Delta x)^{2} \text{Re} \iint\limits_{\mathcal{M} \times \mathcal{N}} {\bf t}^{-}(r{\bf D}_{x}^{2}\rho)  {\bf t}^{-}({\bf D}_{x}^{2}z) {\bf D}_{t}z^{*}\nonumber\\
=\ &-\frac{1}{2}\alpha (\Delta x)^{2} \text{Re} \iint\limits_{\mathcal{M}^* \times \mathcal{N}} {\bf t}^{-}({\bf D}_x(r{\bf D}_{x}^{2}\rho))  {\bf t}^{-}({\bf D}_{x}z) {\bf A}_x{\bf D}_{t}z^{*}-\nonumber\\
&\frac{1}{4}\alpha (\Delta x)^{2} \iint\limits_{\mathcal{M}^* \times \mathcal{N}} {\bf t}^{-}({\bf A}_x(r{\bf D}_{x}^{2}\rho))  {\bf D}_t(|{\bf D}_xz|^2)+\nonumber\\
&\frac{1}{4}\alpha (\Delta x)^{2} \Delta t\iint\limits_{\mathcal{M}^* \times \mathcal{N}} {\bf t}^{-}({\bf A}_x(r{\bf D}_{x}^{2}\rho))  |{\bf D}_t{\bf D}_xz|^2+B_{10}\nonumber\\
=\ &-\frac{1}{2}\alpha (\Delta x)^{2} \text{Re} \iint\limits_{\mathcal{M}^* \times \mathcal{N}} {\bf t}^{-}({\bf D}_x(r{\bf D}_{x}^{2}\rho))  {\bf t}^{-}({\bf D}_{x}z) {\bf A}_x{\bf D}_{t}z^{*}+\nonumber\\
&\frac{1}{4}\alpha (\Delta x)^{2} \iint\limits_{\mathcal{M}^* \times \mathcal{N}} {\bf D}_t({\bf A}_x(r{\bf D}_{x}^{2}\rho))  {\bf t}^+(|{\bf D}_xz|^2)+\nonumber\\
&\frac{1}{4}\alpha (\Delta x)^{2} \Delta t\iint\limits_{\mathcal{M}^* \times \mathcal{N}} {\bf t}^{-}({\bf A}_x(r{\bf D}_{x}^{2}\rho))  |{\bf D}_t{\bf D}_xz|^2+B_{10}+B_{11}
\end{align*}
with
\begin{align*}
&B_{10}=\frac{1}{2}\alpha (\Delta x)^{2} \text{Re} \iint\limits_{\partial\mathcal{M} \times \mathcal{N}} {\bf t}^{-}(r{\bf D}_{x}^{2}\rho)  {\bf t}^{-}({\bf tr}({\bf D}_{x}z)) {\bf D}_{t}z^{*}n_x,\\
&B_{11}=-\frac{1}{4}\alpha (\Delta x)^{2} \iint\limits_{\mathcal{M}^* \times \partial\mathcal{N}} {\bf t}^{+}({\bf A}_x(r{\bf D}_{x}^{2}\rho))  {\bf t}^+(|{\bf D}_xz|^2)n_t.
\end{align*}
By Lemma 2.10, we have ${\bf D}_t({\bf A}_x(r{\bf D}_{x}^{2}\rho))={\bf t}^-(\sigma_1)$. Combining this estimate with \begin{align*}
&{\bf D}_x(r{\bf D}_{x}^{2}\rho)=s^2\mathcal O_{\lambda}(1),\quad{\bf A}_x(r{\bf D}_{x}^{2}\rho)=s^2\mathcal O_{\lambda}(1),\\
&{\bf t}^+(|{\bf D}_xz|^2)\leq C{\bf t}^-(|{\bf D}_xz|^2)+C(\Delta t)^2|{\bf D}_t{\bf D}_xz|^2|,
\end{align*} we have 
\begin{align}\label{1-A.10}
&\frac{1}{2}\alpha (\Delta x)^{2} \text{Re} \iint\limits_{\mathcal{M} \times \mathcal{N}} {\bf t}^{-}(r{\bf D}_{x}^{2}\rho)  {\bf t}^{-}({\bf D}_{x}^{2}z) {\bf D}_{t}z^{*}\nonumber\\
\geq\ &-\alpha  \iint\limits_{\mathcal{M}^* \times \mathcal{N}} {\bf t}^{-}(s^{-1}(s\Delta x)^{2})\mathcal O_{\lambda}(1)|{\bf A}_x{\bf D}_{t}z|^2-\alpha\iint\limits_{\mathcal{M}^* \times \mathcal{N}^*} s(s\Delta x)^{2}\mathcal O_{\lambda}(1) |{\bf D}_xz|^2-\nonumber\\
&\alpha \iint\limits_{\mathcal{M}^* \times \mathcal{N}^*}(\Delta x)^2\sigma_1 |{\bf D}_xz|^2-\alpha \iint\limits_{\mathcal{M}^* \times \mathcal{N}}{\bf t}^-((\Delta x)^2(\Delta t)^2\sigma_1) |{\bf D}_t{\bf D}_xz|^2-\nonumber\\
&\alpha \Delta t\iint\limits_{\mathcal{M}^* \times \mathcal{N}} {\bf t}^{-}((s\Delta x)^{2} )\mathcal O_{\lambda}(1)  |{\bf D}_t{\bf D}_xz|^2+B_{10}+B_{11}.
\end{align}
Now we estimate the first term in (\ref{1-A.10}). By $|{\bf A}_x{\bf D}_{t}z|^2\leq {\bf A}_x(|{\bf D}_{t}z|^2)$ and discrete integration by part for ${\bf A}_x$, the first term in (\ref{1-A.10}) can be estimated as
\begin{align}\label{1-A.11}
&\alpha  \iint\limits_{\mathcal{M}^* \times \mathcal{N}} {\bf t}^{-}(s^{-1}(s\Delta x)^{2})\mathcal O_{\lambda}(1)|{\bf A}_x{\bf D}_{t}z|^2\leq  \alpha  \iint\limits_{\mathcal{M}^* \times \mathcal{N}} {\bf t}^{-}(s^{-1}(s\Delta x)^{2})\mathcal O_{\lambda}(1){\bf A}_x|{\bf D}_{t}z|^2\nonumber\\
\leq \ & \alpha  \iint\limits_{\mathcal{M} \times \mathcal{N}} {\bf t}^{-}(s^{-1}(s\Delta x)^{2})\mathcal O_{\lambda}(1)|{\bf D}_{t}z|^2-B_{12}
\end{align}
with
\begin{align*}
B_{12}=-\alpha \Delta x \iint\limits_{\partial\mathcal{M} \times \mathcal{N}} {\bf t}^{-}(s^{-1}(s\Delta x)^{2})\mathcal O_{\lambda}(1)|{\bf D}_{t}z|^2.
\end{align*}
Substituting (\ref{1-A.10}) and (\ref{1-A.11}) into (\ref{1-A.9}) and using ${\bf t}^+(|z|^2)\leq C{\bf t}^-(|z|^2)+C(\Delta t)^2|{\bf D}_t z|^2$, we obtain the following estimate for $X_5$
\begin{align}\label{A.9}
X_5
\ge\ & -\alpha\iint\limits_{\mathcal{M} \times \mathcal{N}^{*}} \sigma_1 |z|^{2}-\alpha \iint\limits_{\mathcal{M}^* \times \mathcal{N}^*}\left((\Delta x)^2\sigma_1 +s(s\Delta x)^{2}\mathcal O_{\lambda}(1)\right)|{\bf D}_xz|^2-\nonumber\\
&-\alpha\iint\limits_{\mathcal{M} \times \mathcal{N}} {\bf t}^{-}\left(s^{2}\Delta t+(\Delta t)^{2}\sigma_{1}+s^{-1}(s\Delta x)^2\right)\mathcal O_{\lambda}(1) |{\bf D}_{t}z|^{2}-\nonumber\\
&\alpha \Delta t\iint\limits_{\mathcal{M}^* \times \mathcal{N}} {\bf t}^{-}((s\Delta x)^{2} \mathcal O_{\lambda}(1)+(\Delta x)^2(\Delta t)^2\sigma_1)  |{\bf D}_t{\bf D}_xz|^2+\sum_{j=9}^{12}B_{j},
\end{align}
where  
\begin{align*}
\sum_{j=9}^{12}B_{j}
\ge \ &-\alpha \iint\limits_{\mathcal{M} \times \partial\mathcal{N}} {\bf t}^{+}(s)^{2}\mathcal O_{\lambda}(1){\bf t}^{+}(|z|^{2})-\alpha \iint\limits_{\partial\mathcal{M} \times \mathcal{N}} {\bf t}^{-}(1+\Delta x)\left(s^{-1}(s\Delta x)^{2}\right)\mathcal O_{\lambda}(1)|{\bf D}_{t}z|^2\\
&\ -\alpha\iint\limits_{\partial\mathcal{M} \times \mathcal{N}^{*}} s(s\Delta x)^{2}\mathcal O_{\lambda}(1) {\bf tr}(|{\bf D}_{x}z|^{2})-\alpha \iint\limits_{\mathcal{M}^* \times \partial\mathcal{N}} {\bf t}^{+}((s\Delta x)^{2} )\mathcal O_{\lambda}(1)  {\bf t}^+(|{\bf D}_xz|^2).
\end{align*}

Applying integration by parts (\ref{fbjf1}) and 
\begin{align*}
{\rm Re}\left({\bf t}^{-}({\bf D}_x z){\bf D}_t({\bf D}_x z^*)\right)=\frac{1}{2}{\bf D}_t(|{\bf D}_x z|^2)-\frac{1}{2}\Delta t|{\bf D}_{t}{\bf D}_{x}z|^{2}
\end{align*}
 leads to
\begin{align*}
X_6
=\ &-2\alpha \text{Re} \iint\limits_{\mathcal{M}^* \times \mathcal{N}} {\bf t}^{-}({\bf D}_{x}(r{\bf A}_{x}^{2}\rho)){\bf t}^{-} ({\bf D}_{x}z) {\bf A}_x{\bf D}_{t}z^{*}-\\
&\ {\alpha} \iint\limits_{\mathcal{M}^* \times \mathcal{N}} {\bf t}^{-}\left({\bf A}_{x}(r{\bf A}_{x}^{2}\rho)\right) \left({\bf D}_t(|{\bf D}_x z|^2)-\Delta t|{\bf D}_{t}{\bf D}_{x}z|^{2}\right)+B_{13},
\end{align*}
where
\begin{align*}
B_{13}=2\alpha{\rm Re}\iint\limits_{\partial \mathcal M\times \mathcal N}{\bf t}^-(r{\bf A}_x^2\rho){\bf D}_{t} z^*{\bf t}^-\left({\bf tr}({\bf D }_x z)\right) n_x.
\end{align*}
Furthermore, by integration by parts  (\ref{fbjf2})  and (\ref{2.10}) we obtain
\begin{align*}
X_6
=\ &-2\alpha \text{Re} \iint\limits_{\mathcal{M} \times \mathcal{N}} {\bf t}^{-}({\bf A}_{x}{\bf D}_{x}(r{\bf A}_{x}^{2}\rho)) {\bf t}^-({\bf A}_{x}{\bf D}_{x}z) {\bf D}_{t}z^{*}-\\
&\ \frac{1}{2} \alpha (\Delta x)^{2} \text{Re} \iint\limits_{\mathcal{M} \times \mathcal{N}} {\bf t}^{-}({\bf D}_{x}^{2}(r{\bf A}_{x}^{2}\rho)) {\bf t}^{-}({\bf D}_{x}^{2}z) {\bf D}_{t}z^{*}+ \nonumber \\
&\  \alpha \Delta t \iint\limits_{\mathcal{M}^{*} \times \mathcal{N}} {\bf t}^{-}({\bf A}_{x}(r{\bf A}_{x}^{2}\rho)) |{\bf D}_{t}{\bf D}_{x}z|^{2}
+\\
&\  \alpha \iint\limits_{\mathcal{M}^{*} \times \mathcal{N}} {\bf D}_{t}({\bf A}_{x}(r{\bf A}_{x}^{2}\rho)) {\bf t}^{+}(|{\bf D}_{x}z|^{2})+\sum_{j=13}^{15}B_j,
\end{align*}
where
\begin{align*}
&B_{14}=-\alpha\Delta x{\rm Re}\iint\limits_{\partial \mathcal M\times \mathcal N}{\bf t}^-\left({\bf tr}({\bf D}_x(r{\bf A}_x^2\rho)\right){\bf t}^-({\bf D}_x z){\bf D}_{t} z^*,\\
&B_{15}=-\alpha\iint\limits_{\mathcal M^*\times \partial\mathcal N}{\bf t}^+\left({\bf A}_x(r{\bf A}_x^2\rho)\right){\bf t}^+(|{\bf D}_{z}|^2) n_t.
\end{align*}
We note that provided ${\tau \Delta x}{(\delta T^{2})}^{-1} \le 1$ and $\tau \Delta t\left(\delta^{2}T^{3}\right)^{-1} \le {1}/{2}$, from Lemma 2.10 and Lemma 2.6 we deduce 
\begin{align*}
&{\bf A}_{x}{\bf D}_{x}(r{\bf A}_{x}^{2}\rho)= s^2(s\Delta x)^2 \mathcal O_{\lambda}(1),\quad {\bf D}_{x}^{2}(r{\bf A}_{x}^{2}\rho)=(s\Delta x)^2\mathcal O_{\lambda}(1),\\
& {\bf A}_{x}(r{\bf A}_{x}^{2}\rho)=1+(s\Delta x)^2\mathcal O_{\lambda}(1),\quad {\bf D}_{t}({\bf A}_{x}(r{\bf A}_{x}^{2}\rho))={\bf t}^{-}(\sigma_2).
\end{align*}
Therefore, together with ${\bf t}^+(|{\bf D}_xz|^2)\leq C{\bf t}^-(|{\bf D}_x z|^2)+C(\Delta t)^2|{\bf D}_t{\bf D}_x z|^2$, for $X_6$ we obtain the following estimate
\begin{align}\label{A.10}
X_6
\ge\ &\alpha \Delta t \iint\limits_{\mathcal{M}^{*} \times \mathcal{N}} |{\bf D}_{t}{\bf D}_{x}z|^{2}
-\alpha\iint\limits_{\mathcal{M}^{*} \times \mathcal{N}^{*}} \sigma_{2} |{\bf D}_{x}z|^{2}
-\alpha\iint\limits_{\mathcal{M} \times \mathcal{N}^{*}} s(s\Delta x)^{2}\mathcal O_{\lambda}(1) |{\bf A}_{x}{\bf D}_{x}z|^{2}- \nonumber \\
&\alpha\iint\limits_{\mathcal{M} \times \mathcal{N}^{*}} s^{-3}(s\Delta x)^{6}\mathcal O_{\lambda}(1) |{\bf D}_{x}^{2}z|^{2}
-\alpha\iint\limits_{\mathcal{M} \times \mathcal{N}} {\bf t}^{-}(s^{-1}(s\Delta x)^{2})\mathcal O_{\lambda}(1) |{\bf D}_{t}z|^{2} - \nonumber \\
&\alpha\iint\limits_{\mathcal{M}^{*} \times \mathcal{N}} {\bf t}^{-}\left(\Delta t(s\Delta x)^{2}+(\Delta t)^{2}\sigma_{2}\right) \mathcal O_{\lambda}(1)|{\bf D}_{t}{\bf D}_{x}z|^{2}+
\sum_{j=13}^{15}B_j,
\end{align}
where
\begin{align*}
&\sum_{j=13}^{15}B_j\\
\ge\ & -\alpha\iint\limits_{\mathcal{M}^{*} \times \partial\mathcal{N}} {\bf t}^{+}(1+(s\Delta x)^{2}\mathcal O_{\lambda}(1)) {\bf t}^{+}(|{\bf D}_{x}z|^{2})
-\alpha\iint\limits_{\partial\mathcal{M} \times \mathcal{N}} \tau^{-}(\epsilon s^{-1}+s^{-1}(s\Delta x)^{2}\mathcal O_{\lambda}(1)) |{\bf D}_{t}z|^{2} \nonumber \\
&-\alpha\iint\limits_{\partial\mathcal{M} \times \mathcal{N}^{*}} (C(\epsilon)s+s(s\Delta x)^{4}\mathcal O_{\lambda}(1)) {\bf tr}(|{\bf D}_{x}z|^{2}).
\end{align*}
Here we have used Young's inequality with any $\epsilon>0$.

Substituting (\ref{A.9}) and (\ref{A.10}) into (\ref{A.8}) and using $s\Delta x \le 1$ due to ${\tau \Delta x}{(\delta T^{2})}^{-1} \le 1$, we obtain the desired estimate (\ref{A.7}) for $I_{13}$ and complete the proof of Lemma A.2.\hfill$\Box$

\vspace{2mm}

Notice that $I_{22}={\beta^2}/{\alpha^2} I_{11}$. Therefore, from Lemma A.1 we deduce the following estimate for $I_{22}$.

\vspace{2mm}

{\noindent}{\bf Lemma A.3.}\ {\em Provided ${\tau \Delta x}{(\delta T^{2})}^{-1} \le 1$ and $\tau\geq T^2$, we have

	\begin{align}\label{A.1}
I_{22}
\ge\ &4\beta^{2} \iint\limits_{\mathcal{M} \times \mathcal{N}^{*}} s^{3}\lambda^{4}\phi^{3}(\partial_{x}\psi)^{4} |z|^{2}
+4\beta^{2} \iint\limits_{\mathcal{M}^{*} \times \mathcal{N}^{*}} s\lambda^{2}\phi(\partial_{x}\psi)^{2} |{\bf D}_{x}z|^{2} -\nonumber\\
&\beta^2\iint\limits_{\mathcal{M} \times \mathcal{N}^{*}} (s^{3}\lambda^{3}\phi^{3}{\mathcal O}(1) + s^{2}{\mathcal O}_{\lambda}(1) + s^{3}(s\Delta x)^{2}{\mathcal O}_{\lambda}(1) ) |z|^{2}-\nonumber \\
&\beta^2\iint\limits_{\mathcal{M}^{*} \times \mathcal{N}^{*}} (s\lambda\phi {\mathcal O}(1) + s(s\Delta x)^{2}{\mathcal O}_{\lambda}(1)) |{\bf D}_{x}z|^{2}
+BT_{22},
\end{align}
where
\begin{align*}
BT_{22}
\ge\ & \beta^2 c_1 \iint\limits_{\partial\mathcal{M} \times \mathcal{N}^{*}} s^{3}\lambda^{3}\phi^3|z|^{2}+\beta^2 c_1\iint\limits_{{\partial \mathcal M}\times{\mathcal N}^*} s\lambda\phi{\bf tr}(|{\bf D}_x z|^2)-\nonumber \\
&\beta^2\iint\limits_{\partial\mathcal{M} \times \mathcal{N}^{*}} (s^{2}+s^{3}(s\Delta x)){\mathcal O}_{\lambda}(1) |z|^{2}- \beta^2 \iint\limits_{\partial\mathcal{M} \times \mathcal{N}^{*}}s(s\Delta x)^2{\mathcal O}_{\lambda}(1) {\bf tr}(|{\bf D}_{x}z|^{2})
\end{align*}
with a constant $c_1$ such that $c_1\geq \min\{|\psi_x(0)|, |\psi_x(1)|\}>0$.
}

\vspace{2mm}

Next we give the estimate for $I_{23}$.

\vspace{2mm}

 {\noindent}{\bf Lemma A.4.}\ {\em Provided ${\tau \Delta x}{(\delta T^{2})}^{-1} \le 1$ and $\tau\geq T^2$, we have
\begin{align}\label{A.12}
I_{23}
\ge&-\beta\iint\limits_{\mathcal{M} \times \mathcal{N}^{*}} \left(s^{3}(s\Delta x)^{2}\mathcal O_{\lambda}(1)+s\sigma_{3}\right) |z|^{2}
-\beta\iint\limits_{\mathcal{M}^{*} \times \mathcal{N}^{*}} s^{-1}\sigma_{3} |{\bf D}_{x}z|^{2}- \nonumber \\
&\beta\iint\limits_{\mathcal{M} \times \mathcal{N}^{*}} s^{-1}(s\Delta x)^{2}\mathcal O_{\lambda}(1) |{\bf A}_{x}{\bf D}_{x}z|^{2}
-\beta\iint\limits_{\mathcal{M} \times \mathcal{N}^{*}} s^{-1}(s\Delta x)^{2}\mathcal O_{\lambda}(1) |{\bf D}_{x}^{2}z|^{2} -\nonumber \\
&\beta\iint\limits_{\mathcal{M} \times \mathcal{N}} {\bf t}^{-}\left(s^{-1}(s\Delta x)^{2}\mathcal O_{\lambda}(1)+(\Delta t)^{2}s\sigma_{3}+C(\epsilon)s^{2}\Delta  t\mathcal O_{\lambda}(1)\right) |{\bf D}_{t}z|^{2} -\nonumber \\
&\beta\iint\limits_{\mathcal{M}^{*} \times \mathcal{N}} {\bf t}^{-}(\epsilon \Delta  t+s\Delta t\Delta x\mathcal O_{\lambda}(1)+s^{-1}(\Delta  t)^{2}\sigma_{3}) |{\bf D}_{t}{\bf D}_{x}z|^{2}
+BT_{23},
\end{align}
where
\begin{align*}
&BT_{23}\\
\ge\ &-\beta\iint\limits_{\partial\mathcal{M} \times \mathcal{N}^{*}} \left(C(\epsilon)s^3\lambda^2\phi^2+(s\Delta x)\mathcal O_{\lambda}(1)+s\sigma_3\right) |z|^{2}
-\beta\iint\limits_{\partial\mathcal{M} \times \mathcal{N}^{*}} (s\Delta x)\mathcal O_{\lambda}(1) {\bf tr}(|{\bf D}_{x}z|^{2})- \nonumber \\
&\beta\iint\limits_{\partial\mathcal{M} \times \mathcal{N}} {\bf t}^{-}\left(\epsilon s^{-1}+(s\Delta x)\mathcal O_{\lambda}(1)+s\sigma_3(\Delta t)^2+C(\epsilon)s^{2}\Delta t\mathcal O_{\lambda}(1)\right) |{\bf D}_{t}z|^{2}-\\
&\beta\iint\limits_{\mathcal{M}^{*} \times \partial\mathcal{N}}{\bf t}^{+}(s^{2})\mathcal O_{\lambda}(1) {\bf t}^+( |{\bf A}_{x}z|^{2})
-\beta\iint\limits_{\mathcal{M}^{*} \times \partial\mathcal{N}}{\bf t}^{+}\left( |{\bf D}_{x}z|^{2}\right)
\end{align*}
for any $\epsilon>0$.
}

\vspace{2mm}

{\noindent\bf Proof.}\ 
To compute the term $I_{23}$, we borrow some ideas in [\ref{Rosier-Poincare-2009}], where the corresponding continuous case is investigated. Firstly we write
\begin{align} \label{23}
I_{23}
=\ &-\beta \iint\limits_{\mathcal{M} \times \mathcal{N}} i {\bf t}^{-}(2r{\bf A}_{x}{\bf D}_{x}\rho {\bf A}_{x}{\bf D}_{x}z - s\partial_{x}^2\phi z) {\bf D}_{t}z^{*}+\nonumber \\
&\beta \iint\limits_{\mathcal{M} \times \mathcal{N}} i {\bf t}^{-}(2r{\bf A}_{x}{\bf D}_{x}\rho {\bf A}_{x}{\bf D}_{x}z^{*} - s\partial_{x}^2\phi z^{*}) {\bf D}_{t}z.
\end{align}
On one hand, by integration by parts in space, for the first term in (\ref{23}) we have
\begin{align} \label{A.13}
&-\beta \iint\limits_{\mathcal{M} \times \mathcal{N}} i {\bf t}^{-}(2r{\bf A}_{x}{\bf D}_{x}\rho {\bf A}_{x}{\bf D}_{x}z - s\partial_{x}^2\phi z) {\bf D}_{t}z^{*}\nonumber\\
=\ &\beta \iint\limits_{\mathcal{M} \times \mathcal{N}} i {\bf t}^{-}(s\partial^2_{x}\phi z) {\bf D}_{t}z^{*}+2\beta \iint\limits_{\mathcal{M}^* \times \mathcal{N}} i{\bf t}^{-}({\bf A}_{x}(r{\bf A}_{x}{\bf D}_{x}\rho)) {\bf t}^{-}({\bf A}_xz) {\bf D}_{t}{\bf D}_xz^{*}+\nonumber\\
&2\beta \iint\limits_{\mathcal{M}^* \times \mathcal{N}} i{\bf t}^{-}({\bf D}_{x}(r{\bf A}_{x}{\bf D}_{x}\rho)){\bf t}^{-}({\bf A}_xz) {\bf D}_{t}{\bf A}_xz^{*}+B_{16}\nonumber\\
=\ &\beta \iint\limits_{\mathcal{M} \times \mathcal{N}} i {\bf t}^{-}(s\partial^2_{x}\phi z) {\bf D}_{t}z^{*}+2\beta \iint\limits_{\mathcal{M}^* \times \mathcal{N}} i{\bf t}^{-}({\bf A}_{x}(r{\bf A}_{x}{\bf D}_{x}\rho)) {\bf t}^{-}({\bf A}_xz) {\bf D}_{t}{\bf D}_xz^{*}+\nonumber\\
& 2\beta \iint\limits_{\mathcal{M} \times \mathcal{N}} i{\bf t}^{-}({\bf A}_x{\bf D}_{x}(r{\bf A}_{x}{\bf D}_{x}\rho)){\bf t}^{-}\left(z+\frac{1}{4}(\Delta x)^2{\bf D}_x^2 z\right){\bf D}_t z^*+\nonumber\\
& \frac{1}{2}\beta (\Delta x)^2\iint\limits_{\mathcal{M} \times \mathcal{N}} i{\bf t}^{-}({\bf D}^2_{x}(r{\bf A}_{x}{\bf D}_{x}\rho)){\bf t}^{-}\left({\bf A}_x{\bf D}_x z\right){\bf D}_t z^*+B_{16}+B_{17}\nonumber\\
=\ &\beta \iint\limits_{\mathcal{M} \times \mathcal{N}} i {\bf t}^{-}(s\partial^2_{x}\phi z) {\bf D}_{t}z^{*}+2\beta \iint\limits_{\mathcal{M}^* \times \mathcal{N}} i{\bf t}^{-}({\bf A}_{x}(r{\bf A}_{x}{\bf D}_{x}\rho)) {\bf t}^{-}({\bf A}_xz) {\bf D}_{t}{\bf D}_xz^{*}+\nonumber\\
& 2\beta \iint\limits_{\mathcal{M} \times \mathcal{N}} i{\bf t}^{-}({\bf A}_x{\bf D}_{x}(r{\bf A}_{x}{\bf D}_{x}\rho)){\bf t}^{-}(z){\bf D}_t z^*-\nonumber\\
& \frac{1}{2}\beta (\Delta x)^2\iint\limits_{\mathcal{M}^* \times \mathcal{N}} i{\bf t}^{-}({\bf D}_{x}(r{\bf A}_{x}{\bf D}_{x}\rho)){\bf t}^{-}\left({\bf D}_x z\right){\bf D}_t{\bf D}_x z^*+B_{16}+B_{17}+B_{18},
\end{align}
where
\begin{align*}
&B_{16}
=-2\beta \iint\limits_{\partial\mathcal{M} \times \mathcal{N}} i {\bf t}^{-}(r{\bf A}_{x}{\bf D}_{x}\rho) {\bf t}^{-}({\bf tr}({\bf A}_{x}z)) {\bf D}_{t}z^{*}  n_{x}, \nonumber \\
&B_{17}=\beta \Delta x \iint\limits_{\partial\mathcal{M} \times \mathcal{N}} i {\bf t}^-\left({\bf tr}({\bf D}_x(r{\bf A}_x{\bf D}_x\rho))\right){\bf t}^-({\bf tr}({\bf A}_x z)) {\bf D}_{t}z^{*},\\
&B_{18}=\frac{1}{2}\beta(\Delta x)^2\iint\limits_{\partial \mathcal M\times \mathcal N}i {\bf t}^{-}\left({\bf tr}({\bf D}_x(r{\bf A}_x{\bf D}_x\rho))\right){\bf t}^-({\bf tr}({\bf D}_x z)){\bf D}_t z^* n_x.
\end{align*}
We see that
\begin{align}\label{1-A.15} {\bf tr}({\bf A}_{x}z)=
\left\{\begin{array}{ll}
\frac{\Delta x}{2}\left({\bf tr}({\bf D}_{x}z)\right)+z
&{\rm in} \ \{0\} \times \mathcal{N},\\
-\frac{\Delta x}{2}\left( {\bf tr}({\bf D}_{x}z)\right)+z
&{\rm in}\  \{1\} \times \mathcal{N},
\end{array}
\right.
\end{align}
which together with
\begin{align*}
r{\bf A}_x{\bf D}_x \rho=s\lambda\phi\partial_x\psi+s(s\Delta x)^2\mathcal O_{\lambda}(1)
\end{align*}
implies
\begin{align*}
|B_{16}|\leq \ &\beta \iint\limits_{\partial\mathcal{M} \times \mathcal{N}} {\bf t}^{-}(s\Delta x\mathcal ) O_{\lambda}(1) {\bf t}^{-}(|{\bf tr}({\bf D}_{x}z)|) |{\bf D}_{t}z^{*}|+\\
&\beta\iint\limits_{\partial\mathcal{M} \times \mathcal{N}} {\bf t}^{-}(s\lambda\phi\partial_x\psi+s(s\Delta x)^2\mathcal O_{\lambda}(1)) {\bf t}^{-}(|z|) |{\bf D}_{t}z^{*}|\\
\leq \ &  \beta\iint\limits_{\partial\mathcal{M} \times \mathcal{N}^{*}} (s\Delta x)\mathcal O_{\lambda}(1) {\bf tr}(|{\bf D}_{x}z|^{2})+\beta\iint\limits_{\partial\mathcal{M} \times \mathcal{N}} {\bf t}^{-}\left(\epsilon s^{-1}+(s\Delta x)\mathcal O_{\lambda}(1)\right) |{\bf D}_{t}z|^{2} +\nonumber \\
&\beta\iint\limits_{\partial\mathcal{M} \times \mathcal{N}^{*}} C(\epsilon)s^3\lambda^2\phi^2|z|^{2}
\end{align*}
for any $\epsilon>0$. Thus, we can obtain
\begin{align*}
\sum_{j=16}^{18}{\rm Re} (B_j)\ge\ &-\beta\iint\limits_{\partial\mathcal{M} \times \mathcal{N}^{*}} \left(C(\epsilon)s^3\lambda^2\phi^2+(s\Delta x)\mathcal O_{\lambda}(1)\right) |z|^{2}
-\beta\iint\limits_{\partial\mathcal{M} \times \mathcal{N}^{*}} (s\Delta x)\mathcal O_{\lambda}(1) {\bf tr}(|{\bf D}_{x}z|^{2}) \nonumber \\
&-\beta\iint\limits_{\partial\mathcal{M} \times \mathcal{N}} {\bf t}^{-}(\epsilon s^{-1}+(s\Delta x)\mathcal O_{\lambda}(1)) |{\bf D}_{t}z|^{2}.
\end{align*}
On the other hand, using integration by parts with respect to ${\bf A}_x$ and  integration by parts with respect to ${\bf D}_t$ given by (\ref{2.10}), we find for the second term in (\ref{23}) that
\begin{align}\label{A.14}
&\beta \iint\limits_{\mathcal{M} \times \mathcal{N}} i {\bf t}^{-}(2r{\bf A}_{x}{\bf D}_{x}\rho {\bf A}_{x}{\bf D}_{x}z^{*} - s\partial_{x}^2\phi z^{*}) {\bf D}_{t}z\nonumber\\
=\ &2\beta \iint\limits_{\mathcal{M}^* \times \mathcal{N}} i {\bf t}^{-}({\bf A}_{x}(r{\bf A}_{x}{\bf D}_{x}\rho)){\bf t}^{-}({\bf D}_x z^*) {\bf D}_t{\bf A}_{x}z-\beta \iint\limits_{\mathcal{M} \times \mathcal{N}} i {\bf t}^{-}(s\partial_{x}^2\phi z^{*}) {\bf D}_{t}z+\nonumber\\
&\frac{1}{2}\beta(\Delta x)^2 \iint\limits_{\mathcal{M}^* \times \mathcal{N}} i {\bf t}^{-}({\bf D}_{x}(r{\bf A}_{x}{\bf D}_{x}\rho)){\bf t}^{-}({\bf D}_x z^*) {\bf D}_t{\bf D}_{x}z+B_{19}\nonumber\\
=\ & -2\beta \iint\limits_{\mathcal{M}^{*} \times \mathcal{N}} i {\bf D}_{t}({\bf A}_{x}(r{\bf A}_{x}{\bf D}_{x}\rho)) {\bf t}^{+}({\bf A}_{x}z){\bf t}^+( {\bf D}_{x}z^{*})-\nonumber\\
&2\beta \iint\limits_{\mathcal{M}^{*} \times \mathcal{N}} i {\bf t}^{-}({\bf A}_{x}(r{\bf A}_{x}{\bf D}_{x}\rho)) {\bf t}^{-}({\bf A}_{x}z) {\bf D}_{t}{\bf D}_{x}z^{*}-\nonumber\\
&2\beta \Delta t \iint\limits_{\mathcal{M}^{*} \times \mathcal{N}} i {\bf t}^{-}({\bf A}_{x}(r{\bf A}_{x}{\bf D}_{x}\rho)) {\bf D}_{t}{\bf A}_{x}z{\bf D}_{t}{\bf D}_{x}z^{*}-\beta \iint\limits_{\mathcal{M} \times \mathcal{N}} i {\bf t}^{-}(s\partial_{x}^2\phi z^{*}) {\bf D}_{t}z+\nonumber \\
&\frac{1}{2}\beta (\Delta x)^{2} \iint\limits_{\mathcal{M}^{*} \times \mathcal{N}} i {\bf t}^{-}({\bf D}_{x}(r{\bf A}_{x}{\bf D}_{x}\rho)) {\bf t}^{-}({\bf D}_{x}z^{*}) {\bf D}_{t}{\bf D}_{x}z
+B_{19}+B_{20},
\end{align}
where
\begin{align*}
&B_{19}
=-\beta \Delta x\iint\limits_{\partial\mathcal{M} \times \mathcal{N}} i {\bf t}^{-}(r{\bf A}_{x}{\bf D}_{x}\rho) {\bf t}^{-}({\bf tr}({\bf D}_{x}z^*)) {\bf D}_{t}z, \nonumber \\
&B_{20}=2\beta \iint\limits_{\mathcal{M}^* \times \partial\mathcal{N}} i {\bf t}^+\left({\bf A}_x(r{\bf A}_x{\bf D}_x\rho)\right){\bf t}^+({\bf A}_x z) {\bf t}^+({\bf D}_{x}z^{*})n_t,
\end{align*}
satisfy 
\begin{align*}
\sum_{j=19}^{20}{\rm Re}({B_j})
\ge\ &-\beta\iint\limits_{\partial\mathcal{M} \times \mathcal{N}^{*}} (s\Delta x)\mathcal O_{\lambda}(1) {\bf tr}(|{\bf D}_{x}z|^{2})
-\beta\iint\limits_{\partial\mathcal{M} \times \mathcal{N}} {\bf t}^{-}((s\Delta x)\mathcal O_{\lambda}(1)) |{\bf D}_{t}z|^{2}- \nonumber \\
&\beta\iint\limits_{\mathcal{M}^{*} \times \partial\mathcal{N}}{\bf t}^{+}(s^{2})\mathcal O_{\lambda}(1) {\bf t}^+( |{\bf A}_{x}z|^{2})
-\beta\iint\limits_{\mathcal{M}^{*} \times \partial\mathcal{N}}{\bf t}^{+}\left( |{\bf D}_{x}z|^{2}\right).
\end{align*}

Thus, substituting \eqref{A.13} and \eqref{A.14} into \eqref{23} and noticing that ${\rm Re}(I_{23})=I_{23}$, we have
\begin{align}\label{A.15}
I_{23}
=\ &\beta \text{Re} \iint\limits_{\mathcal{M} \times \mathcal{N}}i{\bf t}^{-}(s\partial_{x}^2\phi){\bf t}^{-}(z) {\bf D}_{t}z^{*} -\beta \text{Re} \iint\limits_{\mathcal{M} \times \mathcal{N}} i {\bf t}^{-}(s\partial_{x}^2\phi)  {\bf t}^{-}(z^*){\bf D}_{t}z+\nonumber \\
&2\beta \text{Re} \iint\limits_{\mathcal{M} \times \mathcal{N}} i{\bf t}^{-}({\bf A}_{x}{\bf D}_{x}(r{\bf A}_{x}{\bf D}_{x}\rho)) {\bf t}^{-}(z) {\bf D}_{t}z^{*}-\nonumber \\
&2\beta \text{Re}\iint\limits_{\mathcal{M}^{*} \times \mathcal{N}} i {\bf D}_{t}({\bf A}_{x}(r{\bf A}_{x}{\bf D}_{x}\rho)) {\bf t}^{+}({\bf A}_{x}z){\bf t}^+( {\bf D}_{x}z^{*})-\nonumber \\
&\beta (\Delta x)^{2} \text{Re}\iint\limits_{\mathcal{M}^{*} \times \mathcal{N}} i {\bf t}^{-}({\bf D}_{x}(r{\bf A}_{x}{\bf D}_{x}\rho)) {\bf t}^{-}({\bf D}_{x}z) {\bf D}_{t}{\bf D}_{x}z^*-\nonumber \\
&2\beta \Delta t \text{Re}\iint\limits_{\mathcal{M}^{*} \times \mathcal{N}} i {\bf t}^{-}({\bf A}_{x}(r{\bf A}_{x}{\bf D}_{x}\rho)) {\bf D}_{t}{\bf A}_{x}z{\bf D}_{t}{\bf D}_{x}z^{*}+\sum_{j=13}^{17} {\rm Re}(B_{j})\nonumber\\
=\ &\sum_{j=7}^{12}{X_j}+\sum_{j=16}^{20} {\rm Re}(B_{j}).
\end{align}

We can easily find that
\begin{align*}& {\bf A}_{x}{\bf D}_{x}(r{\bf A}_{x}{\bf D}_{x}\rho)
=-s\lambda^{2}\phi|\partial_{x}\psi|^{2}-s\lambda\phi\partial_{x}^{2}\psi+s(s\Delta x)^{2}\mathcal O_{\lambda}(1),\\
&s\partial_{x}^2\phi=s\lambda^{2}\phi|\partial_{x}\psi|^{2}+s\lambda\phi\partial_{x}^{2}\psi.
\end{align*}
Then we obtain 
\begin{align}\label{A.16}
\sum_{j=7}^{9}{X_j}=\ &-2\beta {\rm Im}\iint\limits_{\mathcal{M} \times \mathcal{N}}{\bf t}^{-}(s\partial_{x}^2\phi){\bf t}^{-}(z) {\bf D}_{t}z^{*} -2\beta {\rm Im}\iint\limits_{\mathcal{M} \times \mathcal{N}} {\bf t}^{-}({\bf A}_{x}{\bf D}_{x}(r{\bf A}_{x}{\bf D}_{x}\rho)) {\bf t}^{-}(z) {\bf D}_{t}z^{*}\nonumber \\
=\ &-2\beta\text{Im} \iint\limits_{\mathcal{M} \times \mathcal{N}}  {\bf t}^{-}\left(s(s\Delta x)^{2}\right)\mathcal O_{\lambda}(1) {\bf t}^{-}(z) {\bf D}_{t}z^{*} \nonumber \\
\ge\ &-\beta\iint\limits_{\mathcal{M} \times \mathcal{N}^{*}} s^{3}(s\Delta x)^{2}\mathcal O_{\lambda}(1) |z|^{2}
-\beta\iint\limits_{\mathcal{M} \times \mathcal{N}} {\bf t}^{-}\left(s^{-1}(s\Delta x)^{2}\right)\mathcal O_{\lambda}(1) |{\bf D}_{t}z|^{2}.
\end{align}
From Lemma \ref{sjgj3}, we see ${\bf D}_{t}({\bf A}_x(r{\bf A}_{x}{\bf D}_{x}\rho))={\bf t}^{-}(\sigma_3)$.
Then for $X_{10}$, using ${\bf t}^+={\bf t}^-+\Delta t {\bf D}_t$, ${\bf A}_x=\frac{1}{2}\Delta x {\bf D}_x+{\bf s}_{-}$ and ${\bf s}_-(\mathcal M^*)=\mathcal M\cup \{0\}$, we obtain
\begin{align}
&\ X_{10}\nonumber \\
=\ &-2\beta \text{Re}\iint\limits_{\mathcal{M}^{*} \times \mathcal{N}} i {\bf t}^-(\sigma_3)({\bf t}^{-}+\Delta t{\bf D}_t)\left(\frac{1}{2}\Delta x{\bf D}_xz+{\bf s}_{-}(z)\right)({\bf t}^{-}+\Delta t{\bf D}_t)\left({\bf D}_x z^*\right)\nonumber\\
\ge\ &-\beta\iint\limits_{\mathcal{M}^* \times \mathcal{N}^{*}} s\sigma_{3} \left((\Delta x)^2|{\bf D}_xz|^{2}+|{\bf s}_-(z)|^2\right)-\beta\iint\limits_{\mathcal{M}^{*} \times \mathcal{N}} {\bf t}^{-}(s^{-1}\sigma_{3})(\Delta t)^{2} |{\bf D}_{t}{\bf D}_{x}z|^{2}-\nonumber \\
&\beta\iint\limits_{\mathcal{M}^{*} \times \mathcal{N}^{*}} s^{-1}\sigma_{3} |{\bf D}_{x}z|^{2}-\beta\iint\limits_{\mathcal{M}^* \times \mathcal{N}} {\bf t}^-(s\sigma_3)(\Delta t)^{2}\left((\Delta x)^2|{\bf D}_{t}{\bf D}_xz|^{2}+|{\bf s}^{-}({\bf D}_tz)|^2\right)\nonumber\\
\geq\ & -\beta\iint\limits_{\mathcal{M} \times \mathcal{N}^{*}} s\sigma_{3} |z|^2-\beta\iint\limits_{\mathcal{M}^{*} \times \mathcal{N}} {\bf t}^{-}(s^{-1}\sigma_{3})(\Delta t)^{2} |{\bf D}_{t}{\bf D}_{x}z|^{2}-\beta\iint\limits_{\mathcal{M}^{*} \times \mathcal{N}^{*}} s^{-1}\sigma_{3} |{\bf D}_{x}z|^{2}-\nonumber\\
&\beta\iint\limits_{\mathcal{M} \times \mathcal{N}} {\bf t}^-(s\sigma_3)(\Delta t)^{2}|{\bf D}_tz|^2-\beta\iint\limits_{\{0\}\times \mathcal{N}^{*}} s\sigma_{3} |z|^2-\beta\iint\limits_{\{0\}\times \mathcal{N}} {\bf t}^-(s\sigma_3)(\Delta t)^{2}|{\bf D}_tz|^2,
\end{align}
where we have used $\Delta x\leq s^{-1}$ due to ${\tau \Delta x}{(\delta T^{2})}^{-1} \le 1$.
Integration by parts with respect to ${\bf D}_x$ given by \eqref{fbjf1}, and using 
\begin{align*}
{\bf A}_{x}{\bf D}_{x}(r{\bf A}_{x}{\bf D}_{x}\rho)=s(s\Delta x)^2\mathcal O_{\lambda}(1),\quad {\bf D}^2_{x}(r{\bf A}_{x}{\bf D}_{x}\rho)=s(s\Delta x)^2\mathcal O_{\lambda}(1)
\end{align*}
 yields the following estimate for the term  $X_{11}$
\begin{align}
X_{11} =\ &\beta (\Delta x)^{2} \text{Re}\iint\limits_{\mathcal{M}^{*} \times \mathcal{N}} i {\bf t}^{-}({\bf D}^2_{x}(r{\bf A}_{x}{\bf D}_{x}\rho)) {\bf t}^{-}({\bf A}_x{\bf D}_{x}z) {\bf D}_{t}z^*+\nonumber \\
&\beta (\Delta x)^{2} \text{Re}\iint\limits_{\mathcal{M}^{*} \times \mathcal{N}} i {\bf t}^{-}({\bf A}_x{\bf D}_{x}(r{\bf A}_{x}{\bf D}_{x}\rho)) {\bf t}^{-}({\bf D}^2_{x}z) {\bf D}_{t}{\bf D}_{x}z^*+B_{18}\nonumber\\
\ge&-\beta\iint\limits_{\mathcal{M} \times \mathcal{N}^{*}} s^{-1}(s\Delta x)^{2}\mathcal O_{\lambda}(1) |{\bf A}_{x}{\bf D}_{x}z|^{2}
-\beta\iint\limits_{\mathcal{M} \times \mathcal{N}^{*}} s^{-1}(s\Delta x)^{2}\mathcal O_{\lambda}(1) |{\bf D}_{x}^{2}z|^{2} \nonumber \\
&-\beta\iint\limits_{\mathcal{M} \times \mathcal{N}} {\bf t}^{-}\left(s^{-1}(s\Delta x)^{2}\right) \mathcal O_{\lambda}(1) |{\bf D}_{t}z|^{2}+B_{21},
\end{align}
where
\begin{align*}
B_{21}=-\beta(\Delta x)^2{\rm Re}\iint\limits_{\partial\mathcal{M} \times \mathcal{N}} i {\bf t}^-\left({\bf tr}({\bf D}_x(r{\bf A}_x{\bf D}_x\rho){\bf D}_x z)\right){\bf D}_t z^* n_x
\end{align*}
satisfies
\begin{align*}
B_{21}\geq-\beta\iint\limits_{\partial\mathcal{M} \times \mathcal{N}^{*}} s^{-1}(s\Delta x)^{2}\mathcal O_{\lambda}(1) {\bf tr}(|{\bf D}_{x}z|^{2})-\beta\iint\limits_{\partial\mathcal{M} \times \mathcal{N}} {\bf t}^{-}(s^{-1}(s\Delta x)^{2}\mathcal O_{\lambda}(1)) |{\bf D}_{t}z|^{2}.
\end{align*}
Finally,  we obtain for $X_{12}$ that
\begin{align}\label{A.19}
X_{12}=\ &-2\beta \Delta t \text{Re} \iint\limits_{\mathcal{M}^{*} \times \mathcal{N}} i {\bf t}^{-1}(s)\mathcal O_{\lambda}(1) {\bf D}_{t}\left(\frac{1}{2}\Delta x {\bf D}_x z+{\bf s}_-(z)\right) {\bf D}_{t}{\bf D}_{x}z^{*} \nonumber \\
\ge&-\beta\iint\limits_{\mathcal{M}^{*} \times \mathcal{N}} \left(\varepsilon \Delta t+{\bf t}^{-}(s)\Delta t\Delta x\mathcal O_{\lambda}(1)\right) |{\bf D}_{t}{\bf D}_{x}z|^{2}-\nonumber \\
&\beta\iint\limits_{\mathcal{M} \times \mathcal{N}} C(\epsilon){\bf t}^{-}(s^{2})\Delta t\mathcal O_{\lambda}(1) |{\bf D}_{t}z|^{2}-\beta\iint\limits_{\{0\}\times \mathcal{N}} C(\epsilon){\bf t}^{-}(s^{2})\Delta t\mathcal O_{\lambda}(1) |{\bf D}_{t}z|^{2}  
\end{align}
for any $\epsilon>0$.

Therefore, from (\ref{A.15})-(\ref{A.19}) and using $s\Delta x \le 1$ and $s\geq 1$ due to $\tau\geq T^2$, we deduce (\ref{A.12}) and then complete the proof of Lemma A.4.\hfill$\Box$

\vspace{2mm}

{\noindent}{\bf Lemma A.5.}\ {\em Provided ${\tau \Delta x}{(\delta T^{2})}^{-1} \le 1$ and $\tau\geq T^2$, we have
\begin{align}\label{A.21}
I_{31}
\ge&-\alpha\iint\limits_{\mathcal{M} \times \mathcal{N}^{*}} Ts^{2}\theta  \mathcal O_{\lambda}(1)|z|^{2}-\alpha\iint\limits_{\mathcal{M}^{*} \times \mathcal{N}^{*}} T(s\Delta x)^{2}\theta \mathcal O_{\lambda}(1)|{\bf D}_{x}z|^{2}
+BT_{31},
\end{align}
where
\begin{align*}
&BT_{31}\geq -\alpha\iint\limits_{\partial\mathcal{M} \times \mathcal{N}^{*}} s^{2}\mathcal O_{\lambda}(1) |z|^{2}
-\alpha\iint\limits_{\partial\mathcal{M} \times \mathcal{N}^{*}} (s\Delta x)^{2}\mathcal O_{\lambda}(1) {\bf tr}(|{\bf D}_{x}z|^{2}).
\end{align*}
}

\vspace{2mm}

{\noindent\bf Proof.}\ The term $I_{31}$ is given by
\begin{align}\label{A.22}I_{31}
=2\alpha \iint\limits_{\mathcal{M} \times \mathcal{N}^{*}} \tau\theta'\varphi s\partial_{x}^2\phi |z|^{2}
-4\alpha \text{Re} \iint\limits_{\mathcal{M} \times \mathcal{N}^{*}} \tau\theta'\varphi r{\bf A}_{x}{\bf D}_{x}\rho z {\bf A}_{x}{\bf D}_{x}z^{*}.
\end{align}

Using the discrete integration by parts for the operator ${\bf A}_x$ and then the discrete integration by parts for the operator ${\bf D}_x$,  we obtain for the second term in (\ref{A.22}) that
\begin{align*}
&-4\alpha \text{Re} \iint\limits_{\mathcal{M} \times \mathcal{N}^{*}} \tau\theta'\varphi r{\bf A}_{x}{\bf D}_{x}\rho z {\bf A}_{x}{\bf D}_{x}z^{*} \nonumber \\
=\ & -2\alpha \iint\limits_{\mathcal{M}^* \times \mathcal{N}^{*}} \tau\theta'{\bf A}_x \left(\varphi r{\bf A}_{x}{\bf D}_{x}\rho\right) {\bf D}_x(|z|^2)-\alpha (\Delta x)^{2} \iint\limits_{\mathcal{M}^{*} \times \mathcal{N}^{*}} \tau\theta{'}{\bf D}_{x}(\varphi r{\bf A}_{x}{\bf D}_{x}\rho) |{\bf D}_{x}z|^{2} +\nonumber \\
&\ B_{22}\nonumber\\
=\ &2\alpha \iint\limits_{\mathcal{M} \times \mathcal{N}^{*}} \tau\theta{'}{\bf A}_{x}{\bf D}_{x}(\varphi r{\bf A}_{x}{\bf D}_{x}\rho) |z|^{2}
-\alpha (\Delta x)^{2} \iint\limits_{\mathcal{M}^{*} \times \mathcal{N}^{*}} \tau\theta{'}{\bf D}_{x}(\varphi r{\bf A}_{x}{\bf D}_{x}\rho) |{\bf D}_{x}z|^{2}+ \nonumber \\
&
B_{22}+B_{23}
\end{align*}
where
\begin{align*}
&B_{22}=2\alpha \Delta x \text{Re} \iint\limits_{\partial\mathcal{M} \times \mathcal{N}^{*}} \tau\theta{'}\varphi r{\bf A}_{x}{\bf D}_{x}\rho z {\bf tr}({\bf D}_{x}z^{*}),\\
&B_{23}=-2\alpha \iint\limits_{\partial\mathcal{M} \times \mathcal{N}^{*}} \tau\theta{'}  {\bf tr}({\bf A}_{x}(\varphi r{\bf A}_{x}{\bf D}_{x}\rho)) |z|^{2} n_x.
\end{align*}

From Lemma \ref{gj2}, we can see that $\tau\theta{'}\varphi s\partial_{x}^2\phi
=Ts^{2}\theta {\mathcal O}_{\lambda}(1)$ and
\begin{align*}
\tau\theta{'}{\bf A}^k_x{\bf D}^l_{x}(\varphi r{\bf A}_{x}{\bf D}_{x}\rho)
=Ts^{2}\theta \mathcal O_{\lambda}(1) + T(s\Delta x)^{2}\theta{\mathcal  O}_{\lambda}(1)=Ts^{2}\theta \mathcal O_{\lambda}(1) ,\quad k,l=0,1,
\end{align*}
due to $s\Delta x\leq 1$ and $s\geq 1$. Then, we obtain the estimate (\ref{A.21}) for $I_{31}$ and complete the proof of Lemma A.5.\hfill$\Box$

\vspace{2mm}

{\noindent}{\bf Lemma A.6.}\ {\em Provided ${\tau \Delta x}{(\delta T^{2})}^{-1} \le 1$ and $\tau\geq T^2$, we have
\begin{align}\label{A.24}
I_{32}\ge&-\beta\iint\limits_{\mathcal{M} \times \mathcal{N}^{*}} s^{2}\mathcal O_{\lambda}(1) |z|^{2}
-\beta\iint\limits_{\mathcal{M} \times \mathcal{N}^{*}} \mathcal O_{\lambda}(1)|{\bf A}_{x}{\bf D}_{x}z|^{2}- \nonumber \\
&\ \beta\iint\limits_{\mathcal{M} \times \mathcal{N}^{*}} s^{-3}\Delta x\mathcal O_{\lambda}(1) |{\bf D}_{x}^{2}z|^{2}
+BT_{23},
\end{align}
where
\begin{align*}
BT_{32}
\ge&-\beta\iint\limits_{\partial\mathcal{M} \times \mathcal{N}^{*}} s^{2}\mathcal O_{\lambda}(1) |z|^{2}-\beta\iint\limits_{\partial\mathcal{M} \times \mathcal{N}^{*}} \mathcal O_{\lambda}(1) {\bf tr}(|{\bf D}_{x}z|^{2}).
\end{align*}
}

\vspace{2mm}

{\noindent\bf Proof.}\ For $I_{32}$, we have
\begin{align}\label{A.25}
I_{32}
=&2\beta \text{Re} \iint\limits_{\mathcal{M} \times \mathcal{N}^{*}} i \tau\theta{'}\varphi r{\bf D}_{x}^{2}\rho z{\bf A}_{x}^{2}z^{*} +2\beta \text{Re} \iint\limits_{\mathcal{M} \times \mathcal{N}^{*}} i \tau\theta{'}\varphi   r{\bf A}_{x}^{2}\rho z{\bf D}_{x}^{2}z^{*}.
\end{align}

Applying integration by parts with respect to ${\bf A}_x$ yields  
\begin{align*}
&\ 2\beta \text{Re} \iint\limits_{\mathcal{M} \times \mathcal{N}^{*}} i \tau\theta{'}\varphi r{\bf D}_{x}^{2}\rho z{\bf A}_{x}^{2}z^{*}\nonumber\\
=\ & 2\beta \text{Re} \iint\limits_{\mathcal{M}^* \times \mathcal{N}^{*}} i \tau\theta{'}\left({\bf A}_{x}(\varphi r{\bf D}_{x}^{2}\rho) {\bf A}_xz +\frac{(\Delta x)^2}{4}{\bf D}_x(\varphi r{\bf D}_{x}^{2}\rho){\bf D}_xz\right) {\bf A}_{x}z^{*}+B_{24}
\end{align*}
with
\begin{align*}
B_{24}=-\beta \Delta x\text{Re} \iint\limits_{\partial\mathcal{M} \times \mathcal{N}^{*}} i \tau\theta{'}\varphi r{\bf D}_{x}^{2}\rho z{\bf tr}({\bf A}_{x}z).
\end{align*}
Again, using integration by parts with respect to ${\bf D}_x$ and noticing that ${\rm Re}(i{\bf A}_x z{\bf A}_x z^*)=0$ and ${\bf A}_x^2 z^*=z^*+(\Delta x)^2/4 {\bf D}_x^2 z^*$ we obtain 
\begin{align}
&\ 2\beta \text{Re} \iint\limits_{\mathcal{M} \times \mathcal{N}^{*}} i \tau\theta{'}\varphi r{\bf D}_{x}^{2}\rho z{\bf A}_{x}^{2}z^{*}\nonumber\\
=\ & -\frac{1}{2}\beta (\Delta x)^2\text{Re} \iint\limits_{\mathcal{M} \times \mathcal{N}^{*}} i \tau\theta{'}\left({\bf D}^2_x(\varphi r{\bf D}_{x}^{2}\rho)z {\bf A}^2_{x}z^{*}+{\bf A}_x{\bf D}_x(\varphi r{\bf D}_{x}^{2}\rho)z {\bf A}_{x}{\bf D}_xz^{*}\right)+B_{21}+B_{22}\nonumber\\
=\ &-\frac{1}{2}\beta (\Delta x)^{2} \text{Re} \iint\limits_{\mathcal{M} \times \mathcal{N}^{*}} i \tau\theta{'}{\bf A}_{x}{\bf D}_{x}(\varphi r{\bf D}_{x}^{2}\rho) z {\bf A}_{x}{\bf D}_{x}z^{*}
-\nonumber\\
&\ \frac{1}{8}\beta(\Delta x)^{4} \text{Re} \iint\limits_{\mathcal{M} \times \mathcal{N}^{*}} i \tau\theta{'}{\bf D}_{x}^{2}(\varphi r{\bf D}_{x}^{2}\rho) z {\bf D}_{x}^{2}z^{*} +B_{24}+B_{25}
\end{align}
with
\begin{align*}
B_{25}=\frac{1}{2}\beta (\Delta x)^2\text{Re} \iint\limits_{\partial\mathcal{M} \times \mathcal{N}^{*}} i \tau\theta{'}z{\bf tr}({\bf D}_x(\varphi r{\bf D}_{x}^{2}\rho){\bf A}_x z^*)n_x.
\end{align*}
Similarly, for the second term in (\ref{A.25})
\begin{align}\label{A.27}
&\ 2\beta \text{Re} \iint\limits_{\mathcal{M} \times \mathcal{N}^{*}} i \tau\theta{'}\varphi   r{\bf A}_{x}^{2}\rho z{\bf D}_{x}^{2}z^{*}\nonumber\\
=\ &-2\beta \text{Re} \iint\limits_{\mathcal{M} \times \mathcal{N}^{*}} i \tau\theta{'}{\bf A}_{x}{\bf D}_{x}(\varphi r{\bf A}_{x}^{2}\rho) z {\bf A}_{x}{\bf D}_{x}z^{*}
-\frac{1}{2} \beta (\Delta x)^{2} \text{Re} \iint\limits_{\mathcal{M} \times \mathcal{N}^{*}} i \tau\theta{'}{\bf D}_{x}^{2}(\varphi r{\bf A}_{x}^{2}\rho) z {\bf D}_{x}^{2}z^{*}+\nonumber\\
 &\ B_{26}+B_{27} 
\end{align}
with
\begin{align*}
&B_{26}=2\beta \text{Re} \iint\limits_{\partial\mathcal{M} \times \mathcal{N}^{*}} i \tau\theta{'}\varphi r{\bf A}_{x}^{2}\rho z {\bf tr}({\bf D}_x z^*)n_x,\\
&B_{27}=-\beta \Delta x\text{Re} \iint\limits_{\partial\mathcal{M} \times \mathcal{N}^{*}} i \tau\theta{'}z{\bf tr}({\bf D}_x(\varphi r{\bf A}_{x}^{2}\rho){\bf D}_x z^*).
\end{align*}

Therefore, it follows from (\ref{A.25})-(\ref{A.27}) that
\begin{align*}
I_{32}
=\ &-\beta \text{Re} \iint\limits_{\mathcal{M} \times \mathcal{N}^{*}} i \tau\theta{'}\left(\frac{1}{2}(\Delta x)^{2} {\bf A}_{x}{\bf D}_{x}(\varphi r{\bf D}_{x}^{2}\rho)+2{\bf A}_{x}{\bf D}_{x}(\varphi r{\bf A}_{x}^{2}\rho)\right) z {\bf A}_{x}{\bf D}_{x}z^{*}
-\nonumber\\
&\  \beta\text{Re} \iint\limits_{\mathcal{M} \times \mathcal{N}^{*}} i \tau\theta{'}\left(\frac{1}{8}(\Delta x)^{4}{\bf D}_{x}^{2}(\varphi r{\bf D}_{x}^{2}\rho)+\frac{1}{2}(\Delta x)^{2}{\bf D}_{x}^{2}(\varphi r{\bf A}_{x}^{2}\rho)\right) z {\bf D}_{x}^{2}z^{*} +\sum_{j=24}^{27}B_{j}.
\end{align*}
Finally, using 
\begin{align*}
&\tau\theta{'}{\bf A}^k_{x}{\bf D}^l_{x}(\varphi r{\bf D}_{x}^{2}\rho)=Ts^3\theta\mathcal O_{\lambda}(1),\ \tau\theta{'}{\bf A}^k_{x}{\bf D}^l_{x}(\varphi r{\bf A}_{x}^{2}\rho)=Ts\theta\mathcal O_{\lambda}(1),\quad k,l=0,1,2
\end{align*}
and (\ref{1-A.15}), we obtain the estimate (\ref{A.24}) for $I_{32}$ and then complete the proof of Lemma A.6.\hfill$\Box$

\vspace{2mm}

{\noindent}{\bf Lemma A.7.}\ {\em Provided ${\tau^3 \Delta t}{(\delta^3 T^{4})}^{-1} \le 1$ and $\tau\geq 1/2(T+T^2)$, we have
\begin{align}\label{A.28}
I_{33}
\ge\ &-\iint\limits_{\mathcal{M} \times \mathcal{N}^{*}} \left(s+\frac{\tau\Delta t}{\delta^{4}T^{5}}\right)\mathcal O_{\lambda}(1)|z|^{2} -\iint\limits_{\mathcal{M} \times \mathcal{N}} {\bf t}^{-}\left(s\Delta t+\frac{\tau\Delta t}{\delta^{4}T^{5}}\right)\mathcal O_{\lambda}(1) |{\bf D}_{t}z|^{2}.
\end{align}
}

\vspace{2mm}

{\noindent\bf Proof.}\ Discrete integrating by parts with respect to ${\bf D}_t$ given by (\ref{2.10}) yields
\begin{align*}
I_{33}
=\ &-2 \text{Re} \iint\limits_{\mathcal{M} \times \mathcal{N}} \tau{\bf t}^{-}(\theta{'})\varphi {\bf t}^{-}(z) {\bf D}_{t}z^{*} \nonumber \\
=\ & 2\text{Re} \iint\limits_{\mathcal{M} \times \mathcal{N}} \tau\varphi\left({\bf D}_t(\theta'){\bf t}^{+}(z)+{\bf t}^-(\theta') {\bf D}_{t}z\right){\bf t}^+(z^*)+B_{25} \\
=\ &2\iint\limits_{\mathcal{M} \times \mathcal{N}} \tau {\bf D}_{t}(\theta^{'})\varphi {\bf t}^{+}(|z|^{2})+2\Delta t \iint\limits_{\mathcal{M} \times \mathcal{N}} \tau{\bf t}^{-}(\theta{'})\varphi |{\bf D}_{t}z|^{2}+\\
&2 \text{Re} \iint\limits_{\mathcal{M} \times \mathcal{N}} \tau{\bf t}^{-}(\theta{'})\varphi  {\bf D}_{t}z{\bf t}^{-}(z^*) +B_{28} 
\end{align*}
with 
\begin{align*}
B_{28}=-2\iint\limits_{\mathcal M\times {\partial\mathcal N}}\tau{\bf t}^+(\theta')\varphi{\bf t}^+(|z|^2)n_t\geq 0.
\end{align*}
Notice that $-I_{33}=2 \text{Re} \iint\limits_{\mathcal{M} \times \mathcal{N}} \tau{\bf t}^{-}(\theta{'})\varphi  {\bf D}_{t}z{\bf t}^{-}(z^*)$. Then we have
\begin{align}\label{A.39}
I_{33}
\geq\ &\iint\limits_{\mathcal{M} \times \mathcal{N}} \tau {\bf D}_{t}(\theta{'})\varphi {\bf t}^{+}(|z|^{2})+\Delta t \iint\limits_{\mathcal{M} \times \mathcal{N}} \tau{\bf t}^{-}(\theta{'})\varphi |{\bf D}_{t}z|^{2}\nonumber\\
\ge\ &-\iint\limits_{\mathcal{M} \times \mathcal{N}} \left(T^{2}{\bf t}^-(s\theta^{2}) + \frac{\tau\Delta t}{\delta^{4}T^{5}}\right)\mathcal O_{\lambda}(1) {\bf t}^{+}(|z|^{2}) -\Delta t\iint\limits_{\mathcal{M} \times \mathcal{N}} {\bf t}^{-}(Ts\theta) \mathcal O_{\lambda}(1) |{\bf D}_{t}z|^{2}.
\end{align}
Since ${\tau^3 \Delta t}{(\delta^3 T^{4})}^{-1} \le 1$ and $\tau\geq 1/2(T+T^2)$, we have $\Delta t\leq 1$. Then combining (\ref{A.39}) and ${\bf t}^{+}(|z|^{2})\geq C{\bf t}^-(|z|^2)+C(\Delta t)^2|{\bf D}_t z|^2$ yields (\ref{A.28}). The proof of Lemma A.7 is completed. \hfill$\Box$

}

\vspace*{0.5cm}
 
\newcounter{cankao}
\begin{list}
	{[\arabic{cankao}]}{\usecounter{cankao}\itemsep=0cm} \centerline{\bf
		References} \vspace*{0.5cm} \small

\item\label{Albano-EJDE-2000} Albano P and Tataru D. Carleman estimates and boundary observability for a coupled parabolic-hyperbolic system. Electronic Journal of Differential Equations, 2000, 22:1-15.

\item\label{Aranson-RMP-2002} Aranson I S and Kramer L.  The world of the complex Ginzburg-Landau equation. Reviews of modern physics, 2002, 74(1): 99.
    
\item\label{Baroun-JDCS-2023} Baroun  M,  Boulite  S,  Elgrou  A  and  Maniar  L.  Null  controllability  for  stochastic parabolic equations with dynamic boundary conditions.  Journal of Dynamical and Control Systems, 2023, 29(4): 1727-1756.
    
\item\label{Baroun-arXiv-2024} Baroun  M,  Boulite  S  and  Maniar  A.   Null  controllability  for  backward  stochastic  parabolic  convection-diffusion  equations  with  dynamic  boundary  conditions. preprint, arXiv: 2401.048 16, 2024.
        
\item\label{Baudouin-IP-2002}  Baudouin  L  and  Puel  J  P.  Uniqueness  and  stability  in  an  inverse  problem  for the Schr\"{o}dinger equation, Inverse Problems, 2002, 18(6): 1537-1554.

\item\label{Baudouin-SICON-2023} Baudouin L and Ervedoza S. Convergence of an inverse problem for a 1-d discrete wave equation. SIAM Journal on Control and Optimization, 2013, 51(1): 556-598.
    
\item\label{Bellassoued-IP-2010}    Bellassoued M and Yamamoto M. Carleman estimates and an inverse heat source problem for the thermoelasticity system. Inverse problems, 2010, 27(1): 015006.
    
\item\label{Bothe-SPTHA-2005} Bothe D, Pr\"{u}ss J, Simonett G. Well-posedness of a two-phase flow with soluble surfactant. Nonlinear Elliptic and Parabolic Problems:  A Special Tribute to the Work of Herbert Amann, 2005, 64: 37-61.
    
\item\label{Boyer-ESAIM-2013} Boyer F. On the penalised HUM approach and its applications to the numerical approximation of null-controls for parabolic problems. ESAIM: Proceedings. EDP Sciences, 2013, 41: 15-58.
    
\item\label{Boyer-ESAIM:COCV-2005} Boyer F and Hern\'{a}ndez-Santamar\'{i}a V. Carleman estimates for time-discrete parabolic equations and applications to controllability.   ESAIM: Control,  Optimisation and Calculus of Variations, 2020, 26: 12.
    
\item\label{Boyer-JMPA-2010}  Boyer F, Hubert F and Le Rousseau J. Discrete carleman estimates for elliptic operators and uniform controllability of semi-discretized parabolic equations.  Journal de math\'{e}matiques pures et appliqu\'{e}es, 2010, 93(3): 240-276.
    
\item\label{Boyer-POINCARE-AN-2014} Boyer F and Le Rousseau J. Carleman estimates for semi-discrete parabolic operators and application to the controllability of semi-linear semi-discrete parabolic equations. Annales de l'Institut Henri Poincar\'{e}, Analyse non lin\'{e}aire, 2014, 31(5): 1035-1078.
    
\item\label{Carreno-arXiv-2023}  Carre\~{n}o N, Mercado A and Morales R.  Local null controllability of a cubic Ginzburg-Landau equation with dynamic boundary conditions. preprint, arXiv:  2301.03429, 2023.
    
\item\label{Carreno-ACM-2023}   Casanova P G and Hern\'{a}ndez-Santamar\'{i}a V. Carleman estimates and controllability results for fully discrete approximations of 1D parabolic equations. Advances in Computational Mathematics, 2021, 47:  1-71.
    
\item\label{Cerpa-JMPA-2022}  Cerpa E, Lecaros R, Nguyen T and P\'{e}rez A. Carleman estimates and controllability for a semi-discrete fourth-order parabolic equation. Journal de Math\'{e}matiques Pures et Appliqu\'{e}es, 2022, 164: 93-130.
    
\item\label{Correa-NARWA-2018}  Corr\^{e}a W J and \"{O}zsarı T. Complex Ginzburg–Landau equations with dynamic boundary conditions[J]. Nonlinear Analysis: Real World Applications, 2018, 41: 607-641.

\item\label{Farkas-MBE-2010}  Farkas J Z and Hinow P.   Physiologically structured populations with diffusion and dynamic boundary conditions.   Mathematical  Biosciences  and  Engineering,  2010, 8(2): 503-513.
    
\item\label{Fu-CRM-2006}  Fu X.  A weighted identity for partial differential operators of second order  and its applications. Comptes Rendus Mathematique, 2006, 342(8): 579-584.
    
\item\label{Fu-JFA-2009} Fu X. Null controllability for the parabolic equation with a complex principal part. Journal of Functional Analysis, 2009, 257(5): 1333-1354.
    
\item\label{Fu-SICON-2017} Fu X and Liu X. Controllability and observability of some stochastic complex Ginzburg-Landau equations.  SIAM Journal on Control and Optimization, 2017, 55(2): 1102-1127.
    
\item\label{Fu-JDE-2017} Fu X and Liu X.  A weighted identity for stochastic partial differential operators and its applications. Journal of Differential Equations, 2017, 262(6): 3551-3582.
    
\item\label{Gal-DCDS-2008}  Gal C and Grasselli M. The non-isothermal Allen-Cahn equation with dynamic boundary conditions.   Discrete and Continuous Dynamical Systems, 2008, 22(4): 1009-1040.
    
\item\label{Hernandez-arXiv-2021}  Hern\'{a}ndez-Santamar\'{i}a  V.  Controllability  of a simplified  time-discrete  stabilized Kuramoto-Sivashinsky system. preprint, arXiv: 2103.12238, 2021.
    
\item\label{Khoutaibi-DCDS-S-2022} Khoutaibi A, Maniar L and Oukdach O. Null controllability for semilinear heat equation with dynamic boundary conditions. Discrete and Continuous Dynamical Systems-S, 2022, 15(6):  1525-1546.
    
\item\label{Lecaros-JDE-2023}  Lecaros R, Morales R, P\'{e}rez A and Zamorano S. Discrete carleman estimates and application to controllability for a fully-discrete parabolic operator with dynamic boundary conditions. Journal of Differential Equations, 2023, 365:  832-881.
      
\item\label{Lecaros-ESAIM:COCV-2021} Lecaros R, Ortega J and P\'{e}rez A.  Stability estimate for the semi-discrete linearized Benjamin-Bona-Mahony equation.  ESAIM: Control,  Optimisation and Calculus of Variations, 2021, 27: 93.
    
\item\label{Lions-SIREV-1988}    Lions J L. Exact controllability, stabilization and perturbations for distributed systems. SIAM review, 1988, 30(1): 1-68.
    
 \item\label{Lions-1988} Lions J L. Contr\^{o}labilité exacte, perturbations et stabilisation de systèmes distribués, Recherches en Mathématiques Appliquées, Tome 1, 8. Masson, Paris, 1988.

\item\label{Lu-ESAIM:COCV-2013}  Lu  Q.  Observability  estimate  for  stochastic  Schr\"{o}dinger  equations  and  its  applications. SIAM Journal on Control and Optimization, 2013, 51(1): 121-144.
    
\item\label{Mercado-SICON-2023} Mercado A and Morales R. Exact Controllability for a Schr\"{o}dinger equation with dynamic boundary conditions. SIAM Journal on Control and Optimization, 2023, 61(6): 3501-3525.

\item\label{Nguyen-MCRF-2014} Nguyen T. Carleman estimates for semi-discrete parabolic operators with a discontinuous diffusion coefficient and applications to controllability[J]. Mathematical Control And Related Fields, 2014, 4(2): 203-259.
    
\item\label{Rosier-Poincare-2009}  Rosier L and Zhang B Y.  Null controllability of the complex Ginzburg-Landau equation. Annales del'Institut Henri Poincar\'{e}, Analyse non lin\'{e}aire, 2009, 26(2):  649-673.
    
\item\label{Wang-arXiv-2024} Wang Y and  Zhao  Q.   Null  controllability  for  stochastic  fourth  order  semi-discrete parabolic equations. preprint, arXiv: 2405.03257, 2024.
    
\item\label{Wu-IP-2024} Wu B, Wang Y and Wang Z.  Carleman estimates for space semi-discrete approximations of one-dimensional stochastic parabolic equation and its applications.  Inverse Problems, 2024, 40(11): 115003.
    
\item\label{Zhao-arXiv-2024} Zhao Q. Null controllability for stochastic semi-discrete parabolic equations[J]. preprint, arXiv: 2402.12651, 2024.
    
\item\label{Zhao-JMAA-2023}  Zhao  Z  and  Zhang W.  Stability  of  a  coefficient  inverse  problem  for the  discrete Schr\"{o}dinger equation and a convergence result.  Journal of Mathematical Analysis and Applications, 2023, 518(1):126665.
    
\item\label{Zuazua-SIREV-2005} Zuazua E. Propagation, observation, and control of waves approximated by finite dif- ference methods. SIAM Review, 2005, 47(2):  197-243.

\end{list}

\end{document}